%  USING TEX
%-------------------------------------------------------------------------
%
%
%
%   
%
%
%
%
%   INSTRUCTIONS:
%   1) If you have any problem with fonts which are not available on
%   your system, search for the statement 
%   \fnts=1
%   and comment it (with a %)
%   2) Run TeX twice to resolve cross-references
%----------------------------------------------------------------------------
%
%%%%%%%%%%%%%%% FORMAT
%\magnification=\magstep1
\hoffset=0.cm
\voffset=0.5truecm
\hsize=16.5truecm 
\vsize=22.0truecm
\headline={\ifnum\pageno>1\ifodd\pageno\rightheadline\else\leftheadline\fi\fi}
\def\rightheadline{\eightrm\hfil\break  
On the homotopy invariance of higher signatures for manifolds with boundary}
\def\leftheadline
{\eightrm E. Leichtnam, J. Lott and P. Piazza
\hfil\break}
\footline={\ifodd\pageno\rightfootline\else\leftfootline\fi}
\def\rightfootline{\hfil\break\folio\hfil\break} 
\def\leftfootline{\hfil\break\folio\hfil\break}
\baselineskip=14pt plus0.1pt minus0.1pt 
\parindent=25pt
\lineskip=4pt\lineskiplimit=0.1pt      
\parskip=0.1pt plus1pt
%

%%%%%%%%%%%%%%%%%%%%%%%%%%%  FONTS

\font\eightrm=cmr8

\font\sixrm=cmr6

%%%%%%%%%%%%%%%%%%%%%  Page numbering

\def\data{\number\day/\ifcase\month\or gennaio \or febbraio \or marzo \or
aprile \or maggio \or giugno \or luglio \or agosto \or settembre
\or ottobre \or novembre \or dicembre \fi/\number\year}

%%\newcount\tempo
%%\tempo=\number\time\divide\tempo by 60}

\setbox200\hbox{$\scriptscriptstyle \data $}

\newcount\pgn 
\pgn=1
\def\foglio{\veroparagrafo:\number\pgn
\global\advance\pgn by 1}

\global\newcount\numsec
\global\newcount\numfor
\global\newcount\numfig
\global\newcount\numtheo

\gdef\profonditastruttura{\dp\strutbox}

\def\senondefinito#1{\expandafter\ifx\csname#1\endcsname\relax}

\def\SIA #1,#2,#3 {\senondefinito{#1#2}%
   \expandafter\xdef\csname #1#2\endcsname{#3}\else
   \write16{???? ma #1,#2 e' gia' stato definito !!!!}\fi}

\def\etichetta(#1){(\veroparagrafo.\veraformula)
   \SIA e,#1,(\veroparagrafo.\veraformula)
   \global\advance\numfor by 1
   \write15{\string\FU (#1){\equ(#1)}}
   \write16{ EQ \equ(#1) == #1  }}

\def\FU(#1)#2{\SIA fu,#1,#2 }

%------------------- theorems ----------------------------
%
\def\tetichetta(#1){{\veroparagrafo.\verotheo}%
   \SIA theo,#1,{\veroparagrafo.\verotheo}
   \global\advance\numtheo by 1%
   \write15{\string\FUth (#1){\thm[#1]}}%
   \write16{ TH \thm[#1] == #1  }}

\def\FUth(#1)#2{\SIA futh,#1,#2 }
%
%--------------------------------------------------------

\def\getichetta(#1){Fig. \verafigura
 \SIA e,#1,{\verafigura}
 \global\advance\numfig by 1
 \write15{\string\FU (#1){\equ(#1)}}
 \write16{ Fig. \equ(#1) ha simbolo  #1  }}

\newdimen\gwidth

\def\BOZZA{
 \def\alato(##1){
 {\vtop to \profonditastruttura{\baselineskip
 \profonditastruttura\vss
 \rlap{\kern-\hsize\kern-1.2truecm{$\scriptstyle##1$}}}}}
 \def\galato(##1){ \gwidth=\hsize \divide\gwidth by 2
 {\vtop to \profonditastruttura{\baselineskip
 \profonditastruttura\vss
 \rlap{\kern-\gwidth\kern-1.2truecm{$\scriptstyle##1$}}}}}
 \def\talato(##1){\rlap{\sixrm\kern -1.2truecm ##1}}
}

\def\alato(#1){}
\def\galato(#1){}
\def\talato(#1){}

\def\veroparagrafo{\ifnum\numsec<0 A\number-\numsec\else
   \number\numsec\fi}
\def\veraformula{\number\numfor}
\def\verotheo{\number\numtheo}
\def\verafigura{\number\numfig}

\def\Thm[#1]{\tetichetta(#1)}
\def\thf[#1]{\senondefinito{futh#1}#1\else
   \csname futh#1\endcsname\fi}
\def\thm[#1]{\senondefinito{theo#1}thf[#1]\else
   \csname theo#1\endcsname\fi}

\def\Eq(#1){\eqno{\etichetta(#1)\alato(#1)}}
\def\eq(#1){\etichetta(#1)\alato(#1)}
\def\eqv(#1){\senondefinito{fu#1}#1\else
   \csname fu#1\endcsname\fi}
\def\equ(#1){\senondefinito{e#1}eqv(#1)\else
   \csname e#1\endcsname\fi}

%%%%%%%%%%%%%%%%%% Numbering toward the future and preeceding
%%%%%%%      sections not inserted in the file to be compiled
\def\include#1{
\openin13=#1.aux \ifeof13 \relax \else
\input #1.aux \closein13 \fi}
\openin14=\jobname.aux \ifeof14 \relax \else
\input \jobname.aux \closein14 \fi
\openout15=\jobname.aux
%%%%%%%%%%%%%%%%%%%%%%%%%%%%

%\footline={\rlap{\hbox{\copy200}\ $\st[\number\pageno]$}\hss\tenrm
%\foglio\hss}%Comment this if you do not want date at the bottom
%---------------- fonti disponibili ---------------------------
%
\newcount\fnts
\fnts=0
\fnts=1 %-----comment if msam, msbm, eufm are not available

\def\tthsp{\kern .083333 em}

%------------------------ itemizing
%

\def\indbox#1{\hbox to \parindent{\hfil\break\ #1\hfil\break} }

\def\ref[#1]{[#1]}

\def\beginsubsection#1\par{
\vskip 0pt plus .05\vsize \penalty -250%
  \vskip0pt plus -.05\vsize \bigskip \vskip \parskip 
  \leftline {\bf #1}\nobreak%\bigskip\bigskip\bigskip
\noindent}

\def\beginsubsection#1\par{\bigskip\leftline{\it #1}\nobreak\smallskip
            \noindent}

%
%
%-------------------------------------------------------------------
%..................If there are not fonts
%
\newfam\msafam
\newfam\msbfam
\newfam\eufmfam
\ifnum\fnts=0

\def\mbox{
\vbox{ \hrule width 6pt
   \hbox to 6pt{\vrule\vphantom{k} \hfil\break\vrule}
   \hrule width 6pt}
}
\def\QED{\ifhmode\unskip\nobreak\fi\quad
  \ifmmode\mbox\else$\mbox$\fi}
\let\restriction=\lceil
%
%.................. if there are fonts
%
\else
\def\hexnumber#1{%
\ifcase#1 0\or 1\or 2\or 3\or 4\or 5\or 6\or 7\or 8\or
9\or A\or B\or C\or D\or E\or F\fi}
%--------------------------------------
\font\tenmsa=msam10
\font\sevenmsa=msam7
\font\fivemsa=msam5
\textfont\msafam=\tenmsa
\scriptfont\msafam=\sevenmsa
\scriptscriptfont\msafam=\fivemsa
\edef\msafamhexnumber{\hexnumber\msafam}%
\mathchardef\restriction"1\msafamhexnumber16
\mathchardef\square"0\msafamhexnumber03
%\mathchardef\supneq"3\msafamhexnumberXX
\def\QED{\ifhmode\unskip\nobreak\fi\quad
  \ifmmode\square\else$\square$\fi}
%--------------------------------------
%
\font\tenmsb=msbm10
\font\sevenmsb=msbm7
\font\fivemsb=msbm5
\textfont\msbfam=\tenmsb
\scriptfont\msbfam=\sevenmsb
\scriptscriptfont\msbfam=\fivemsb
\def\Bbb#1{\fam\msbfam\relax#1}
\edef\msbfamhexnumber{\hexnumber\msbfam}%
%\mathchardef\supneq"3\msbfamhexnumber22
%\mathchardef\supneqx"3\msbfamhexnumber2A
\mathchardef\supneq"3\msbfamhexnumber28
%--------------------------------------
%
\font\teneufm=eufm10
\font\seveneufm=eufm7
\font\fiveeufm=eufm5
\textfont\eufmfam=\teneufm
\scriptfont\eufmfam=\seveneufm
\scriptscriptfont\eufmfam=\fiveeufm

%--------------------------------------
%
%\font\teneusb=eusb10
%\font\seveneusb=eusb7
%\font\fiveeusb=eusb5
%\newfam\eusbfam
%\textfont\eusbfam=\teneusb
%\scriptfont\eusbfam=\seveneusb
%\scriptscriptfont\eusbfam=\fiveeusb
%\def\script#1{{\fam\eusbfam\relax#1}}
%%%%%%%%%%%%%%%%%%%%

\outer\def\nproclaim#1 [#2]#3. #4\par{\medbreak \noindent
   \talato(#2){\bf #1 \Thm[#2]#3.\enspace }%
   {\sl #4\par }\ifdim \lastskip <\medskipamount 
   \removelastskip \penalty 55\medskip \fi}

\def\thmm[#1]{#1}
\def\teo[#1]{#1}

%-------------------------------------------------------------------
%----------------- tilde
%
\def\sttilde#1{%
\dimen2=\fontdimen5\textfont0
\setbox0=\hbox{$\mathchar"7E$}
\setbox1=\hbox{$\scriptstyle #1$}
\dimen0=\wd0
\dimen1=\wd1
\advance\dimen1 by -\dimen0
\divide\dimen1 by 2
\vbox{\offinterlineskip%
   \moveright\dimen1 \box0 \kern - \dimen2\box1}
}
%

%
%-------------------------------------------------------------------

%%%%
%%%% Macros
%%%%

\def\NN{{\Bbb N}}
\def\QQ{{\Bbb Q}}
\def\ha{{1\over 2}}
\def\AAA{{\Bbb A}}
\def\pAAA{\widetilde{\AAA}}
\def\BBB{{\Bbb B}}
\def\pBBB{\widetilde{\BBB}}

\let\ZZ=\integer
\let \RR = \real
\let \CC = \complex
\def\Isom{\rm Isom}
\def\GL{\rm GL}
\def\eps{\epsilon}
\def\Si{{\cal D}^{{\rm sign}}}

\def\SC{{\cal D}_C^{\rm sign}}
\def\SCp{{\cal D}_C^{{\rm sign},+}}
\def\SCc{{\cal D}_C^{\rm sign,cyl}}

\def\SCcp{{\cal D}_C^{{\rm sign,cyl},+}}
\def\SCAPSp{{\cal D}_C^{{\rm sign,APS},+}}
\def\Di{D\kern-6pt/}

\def\cDi{{\cal D}}
\def \P{{\cal P}}

\def \Q{{\cal Q}}

\def \Bi{{\cal B}^\infty}
\def \Fl{{\cal V}}
\def \Fli{\Fl^\infty}

\def \whW {\widehat{W}}

\def\nuint{\raise10pt\hbox{$\nu$}\kern-6pt\int}

\def \ub{{}^b}

\def \Cl{ \mathop{\rm Cl}\nolimits}
\def \tM{M^\prime}

\let\pa=\partial
\def \eb2M{\tM^2_{{\rm e}b}}

\def \ch{ \mathop{\rm ch}\nolimits}

\def \TR{ \mathop{\rm TR}\nolimits}
\def \STR{{\rm STR}}
\def \STr{ \mathop{\rm STr}\nolimits}
\def \STR1{ {\rm STR}_{{\rm Cl}(1)}}
\def \STr1{ {\rm STr}_{{\rm Cl}(1)}}

\def \bSTR{ b\mathop{\rm -STR}\nolimits}

\def \Ind{ \mathop{\rm Ind}\nolimits}

\def \Hom{ \mathop{\rm Hom}\nolimits}

\def \Dom{ \mathop{\rm Dom}\nolimits}
\def \IM{ \mathop{\rm Im}\nolimits}

\let \id=\Iidentity

\def\wsquare{$\sqcap\kern-7pt\sqcup $}
\def\boxtimes{\sqcap\kern-7pt\sqcup$\vskip-13pt$\kern-.8pt\times}

\def \dw{D_{\whW}}
\def \tw{\tau_{\whW}}
\def \Dw{{\cal D}^{\rm sign}_{\whW}}
\def \dco{D_C^{\rm cone}}
\def \Ker{{\rm Ker}}

\def \SCb{{\cal D}_C^{{\rm sign}, b}}
\def \SCbp{{\cal D}_C^{{\rm sign}, b,+}}
\def \SCbs{{\cal D}_C^{{\rm sign},b,s}}
\def \interior{{\rm int}}
\def \HH{{\rm H}}
\def \vol{{\rm vol}}
\def \const{{\rm const. }}

\expandafter\ifx\csname sezioniseparate\endcsname%---Do not touch these
   \relax\input macro \fi                        %---2 lines
\font\ttlfnt=cmcsc10 scaled 1200 %small caps
 %bold italic text mode
%
%\BOZZA
%
\begingroup

\endgroup

\expandafter\ifx\csname sezioniseparate\endcsname%--- Do not touch
   \relax\input macro \fi                        %--- these 2 lines
%--------------------------------- #Beginning#

%\centerline{\ttlfnt On the Novikov conjecture on  manifolds}
%\centerline{\ttlfnt with boundary II}
\centerline{\ttlfnt On the homotopy invariance of higher signatures}
\centerline{\ttlfnt  for manifolds with boundary}
\vskip 0.5truecm
\centerline{ Eric LEICHTNAM
\footnote{$^{*}$}{Institut de Jussieu,
Tour 46-00, $3^{me}$ \'etage (Alg\`ebres d'op\'erateurs), 4
place Jussieu, 75005 Paris, France,
{\tt leicht@math.jussieu.fr}, research partially supported by a CNR-CNRS
cooperation project,}
$\;$ John LOTT
\footnote{$^{**}$}{Department of Mathematics,
University of Michigan, Ann Arbor, 
MI 48109-1109, USA,
\hfil\break
{\tt lott@math.lsa.umich.edu}, research partially supported by NSF
grant DMS-9704633} 
 and  $\;$ Paolo PIAZZA
\footnote{$^{***}$}{Universit\`a di 
Roma ``La Sapienza'', Istituto ``Guido Castelnuovo'',
P.le A. Moro 2, I-00185 Roma, Italy, {\tt piazza@mat.uniroma1.it},
research partially supported by a CNR-CNRS cooperation project and by
M.U.R.S.T.
}}
\vskip 0.5cm

%\address{\ninerm
%\'Ecole Normale Sup\'erieure, DMI, 45 rue d'Ulm, F-75230 Paris Cedex 05 
%\hfil\breakl\break
%and 
%\hfil\breakl\break
%Universit\`a di Roma ``La Sapienza'', Istituto Matematico ``G. Castelnuovo'',
%P.le A. Moro 2, I-00185 Roma\hfil\breakl\break 
%\smallno
%e-mail: Eric.Leichtnam@ens.fr \hfil\breakl\break
%e-mail: piazza@mat.uniroma1.it \hfil\breakl\break}
%---------------------------

\numsec=0
\numfor=1
\numtheo=1
\pgn=1
\noindent
{\bf Abstract :} If $M$ is a compact
oriented manifold-with-boundary whose
fundamental group is virtually nilpotent or Gromov-hyperbolic, we show
that the higher signatures of $M$ are oriented-homotopy invariants.

\beginsection 0. Introduction

The Novikov Conjecture hypothesizes that certain numerical invariants
of closed oriented manifolds, called higher signatures,
are oriented-homotopy invariants. It is
natural to ask if there is an extension of the Novikov Conjecture to
manifolds with boundary.  Such an extension was made in 
[L2],[L5]. In this paper we show that if the relevant
discrete group is virtually
nilpotent or Gromov-hyperbolic then the higher signatures defined in [L2],[L5]
are oriented-homotopy invariants.

Before giving our result, 
let us recall the statement of Novikov's conjecture.
Let $M$ be a closed 
oriented smooth manifold.  Let $L \in \HH^*(M;\QQ)$ be the
Hirzebruch $L$-class and let $*L \in \HH_*(M;\QQ)$ be its
Poincar\'e dual. If $\Gamma$ is a finitely-generated discrete group,
let $B\Gamma$ denote its classifying space. Recall that
$\HH^*(B\Gamma; \QQ) \cong \HH^*(\Gamma; \QQ)$, the rational
group cohomology of $\Gamma$. Let
$\nu : M \rightarrow B\Gamma$ be a continuous map, defined up to homotopy.  
The Novikov Conjecture 
hypothesizes that the higher signature $\nu_*(*L) \in \HH_*(B\Gamma;\QQ)$
is an oriented-homotopy invariant of the pair
$(M, \nu)$. Equivalently, if $\tau \in \HH^*(\Gamma;\QQ)$
then $\langle \nu_*(*L), \tau \rangle = 
\langle L \cup \nu^* \tau, [M] \rangle$
should be an oriented-homotopy invariant.
If $\Gamma$ is virtually nilpotent or Gromov-hyperbolic then the validity
of the Novikov Conjecture was proven by Connes and Moscovici 
[CM, Theorem 6.6].

Now let $M$ be a compact oriented manifold-with-boundary, equipped with a 
continuous
map $\nu : M \rightarrow B\Gamma$ which is defined up to homotopy.
The formal definition of the higher signature of $M$, from
[L2, (67)] and [L5, Definition 10], is
$$\sigma_M = \left(
\int_M L(TM) \wedge \omega \right) - \widetilde{\eta}_{\partial M} \in 
\overline{\HH}_*({\cal B}^\infty).
\Eq(highersig)$$
The terms of this equation will be defined later in the paper. 
Briefly, \hfil\break
\noindent
1. $L(TM) \in \Omega^*(M)$ is the $L$-form of
$M$ associated to a Riemannian metric which is a product near the
boundary,  \hfil\break
\noindent
2. $\Bi$ is a ``smooth'' subalgebra of the reduced group $C^*$-algebra 
$C^*_r(\Gamma)$, i.e. $\CC\Gamma \subset \Bi \subset C^*_r(\Gamma)$
and $\Bi$ is closed under the holomorphic
functional calculus in $C^*_r(\Gamma)$,
\hfil\break
\noindent
3. $\widetilde{\eta}_{\partial M}$, the higher eta-form [L2, Definition 11],
is an element of the space $\overline{\Omega}_*({\cal B}^\infty)$  
of noncommutative differential forms 
[K, Sections 1.3 and 4.1] and can be thought of
 as a boundary correction term,
\hfil\break
\noindent
4. $\omega$ is a certain closed 
biform in $\Omega^*(M) \, \widehat{\otimes} \, \overline{\Omega}_*(
\CC \Gamma)$ [L1, Section V]
and 
\hfil\break
\noindent
5. $\overline{\HH}_*({\cal B}^\infty)$ is the
noncommutative de Rham homology of ${\cal B}^\infty$ [Co, p. 185],[K, 
Section 4.1].

As in [L2, Section 4.7] and [L5, Assumption 2], 
in order to make sense of the higher eta-form $\widetilde{\eta}_{\partial M}$
we must make an assumption about $\partial M$. To be slightly more general,
let $F$ be a closed oriented manifold,
equipped with a continuous map $\nu_0 : F \rightarrow B\Gamma$ which is
defined up to homotopy. Associated to $\nu_0$ is a
normal $\Gamma$-cover $\pi : F^\prime \rightarrow F$ of $F$.
There is an associated flat $C^*_r (\Gamma)$-vector bundle 
${\cal V}_0 = C^*_r (\Gamma)
\times_\Gamma F^\prime$ on $F$.
Let $\HH^*(F; {\cal V}_0) = \Ker(d)/\IM(d)$ denote the usual
(unreduced) de Rham or simplicial
cohomology of $F$, computed using the local
system ${\cal V}_0$.
Let $\overline{\HH}^*(F; {\cal V}_0) = \Ker(d)/
\overline{\IM(d)}$ denote the 
reduced cohomology.
There is an obvious surjection
$s : \HH^*(F; {\cal V}_0) \rightarrow 
\overline{\HH}^*(F; {\cal V}_0)$. 
\medskip
\noindent 
{\bf Assumption 1 :} \hfil\break
a. The map $s : \HH^k(F; {\cal V}_0) \rightarrow 
\overline{\HH}^k(F; {\cal V}_0)$ is an isomorphism for
$k = \left[ {{\rm dim}(F) + 1 \over 2} \right]$. \hfil\break
b. If $\dim(F) = 2k$ then $\overline{\HH}^k(F; {\cal V}_0)$ admits a
(stable) Lagrangian subspace.

\medskip

Assumption 1 is a homotopy-invariant assumption on $F$. If $F$ is endowed with
a Riemannian metric then an equivalent
formulation of Assumption 1.a. is 
: \hfil\break
\noindent
1. If $\dim(F) = 2k$ then the differential form Laplacian on 
$\Omega^k(F^\prime)$ has a strictly positive spectrum on the orthogonal
complement of its kernel.
\hfil\break
\noindent
2. If $\dim(F) = 2k-1$ then the differential form Laplacian on 
$\Omega^{k-1}(F^\prime)/\Ker(d)$ has a strictly positive spectrum.

 Given Assumption 1.a, 
Assumption 1.b. is equivalent to saying that the index of the signature
operator of $F$, as an element of $K_0(C^*_r(\Gamma))$, vanishes. 
As examples,
\hfil\break
\noindent
(a) If $F$ has a cellular decomposition without any cells of
dimension $k = \left[ {{\rm dim}(F) + 1 \over 2} \right]$ then 
Assumption 1 is satisfied. \hfil\break
\noindent
(b) If $\Gamma$ is finite and the signature of $F$ vanishes then 
Assumption
1 is satisfied. \hfil\break
\noindent
(c) Let $F_1$ and $F_2$ be even-dimensional manifolds, with 
$F_1$ a connected closed hyperbolic manifold and $F_2$ a closed manifold with
vanishing signature. Put $\Gamma = \pi_1(F_1)$. If $F = F_1 \times F_2$
then Assumption 1 is satisfied. 
   \hfil\break
\noindent
(d) If $\dim(F) = 3$, $F$ is connected and $\Gamma = \pi_1(F)$ then, assuming
Thurston's geometrization conjecture, Assumption 1
is satisfied if and only if $F$ is a connected sum of
spherical space forms, $S^1 \times S^2$'s and twisted circle bundles
$S^1 \times_{\ZZ_2} S^2$ over
$\RR P^2$. 

Suppose that Assumption 1 is satisfied. 
If $\dim(F) = 2k$, choose a (stable) Lagrangian subspace $L$ of 
$\overline{\HH}^k(F; {\cal V}_0)$. 
Then the higher eta-form 
$\widetilde{\eta}_{F}$ is well-defined. 
In the case of a manifold-with-boundary $M$, let $i : \partial M \rightarrow
M$ be the boundary inclusion. We take $F = \partial M$ and
$\nu_0 = \nu \circ i$.
In this case, if Assumption 1.a. holds then
Assumption 1.b. holds.
 
The main result of this paper is the following :

\nproclaim Theorem [mainthm]. If $\partial M$ satisfies Assumption 1 then 
$\sigma_M$ is
an oriented-homotopy invariant of the pair $(M, \nu)$.

By oriented-homotopy 
invariance of $\sigma_M$, we mean the following.  Suppose that
$h : (M_2, \partial M_2) \rightarrow (M_1, \partial M_1)$ is an
degree $1$ homotopy equivalence of pairs.  In particular, 
$h(\partial M_2) \subset \partial M_1$, 
but $h \big|_{\partial M_2}$
is not assumed to be a homeomorphism from $\partial M_2$ to 
$\partial M_1$. Suppose that
there are continuous maps $\nu_i : M_i \rightarrow B\Gamma$ such that
$\nu_2$ is homotopic to $\nu_1 \circ h$.
If $\dim(M_i) = 2k+1$, we assume that the (stable) Lagrangian subspaces
for the boundaries are
related by $(\partial h)^*(L_1) = L_2$.
Then
$\sigma_{M_1}$, computed using $\nu_1$, equals $\sigma_{M_2}$, computed
using $\nu_2$.  
(If $\dim(M) = 2k+1$ then $\sigma_M$ generally depends on
the choice of $L$.)

In order to obtain numerical invariants from $\sigma_M$, we must make
an assumption about the smooth subalgebra ${\cal B}^\infty$.

\medskip
\noindent
{\bf Assumption 2 :} Each class $\tau \in \HH^*(\Gamma; \CC)$ has a
cocycle representative whose corresponding cyclic cocycle
$Z_\tau \in ZC^*(\CC \Gamma)$ 
has a continuous extension from $\CC \Gamma$ to ${\cal B}^\infty$.
\medskip

If $\Gamma$ is virtually nilpotent or Gromov-hyperbolic then it is known
that smooth subalgebras ${\cal B}^\infty$ of $C^*_r(\Gamma)$ satisfying
Assumption 2 exist 
[dlH, Section 2],[J, Theorem 4.1]. 
We write $\langle \sigma_M, \tau \rangle$ for the pairing
of $\sigma_M$ with $Z_\tau$. 

\nproclaim Corollary [hominv]. Under Assumptions 1 and 2, the
higher signatures $\langle \sigma_M, \tau \rangle$ are oriented-homotopy
invariants.

As special cases of Corollary \thm[hominv], if $\partial M = \emptyset$ 
then 
$
\langle \sigma_M, \tau \rangle= \const \, \langle L(TM) \cup\nu^*\tau,[M]
\rangle
$ [L1, Corollary 2] 
and so we recover the Connes-Moscovici result [CM, Theorem 6.6]. At the other
extreme, if 
$\partial M \neq \emptyset$, 
$\Gamma = \{e\}$, $\Bi = \CC$ and $\tau = 1 \in \HH^0(\{e\}; \CC)$ then 
$\langle \sigma_M, \tau \rangle$ is the Atiyah-Patodi-Singer formula for
the signature of $M$ [APS, Theorem 4.14], 
which is clearly an oriented-homotopy invariant
of $M$.

Let us comment on Assumptions 1 and 2.  Assumption 2 is a technical condition
on $\Gamma$.  Assumption 1 is more germane and
is necessary for both analytic and topological
reasons.  On the analytic side, something like Assumption 1 is necessary in 
order to make sense of the formal expression for
$\widetilde{\eta}_{F}$. On the topological side, Assumption 1
implies that the higher signature of $F$, with respect
to $\nu_0$, vanishes.
Of course, if $F = \partial M$ then its higher signature
vanishes simply because $\partial M$ is a boundary, but Assumption 1 gives
a reason for the vanishing which is intrinsic to $\partial M$.  

(For clarity, we note that if we just want to
{\it define} $\langle \sigma_M, \tau
\rangle$ then 
we can get by with something weaker than Assumption 2. Namely, for a
connected component $F$ of $\partial M$, put
$\Gamma_F = \IM(\pi_1(F) \rightarrow \pi_1(M) \rightarrow \Gamma)$. Let
$\Bi_F$ be a smooth subalgebra of $C^*_r(\Gamma_F)$.
Then it is enough to assume that for each $F$,
$\tau \big|_{\Gamma_F}$ extends to a
cyclic cocycle on $\Bi_F$. For example, if $\partial M = \emptyset$ then
there is no assumption on $\Gamma$ and we recover the Novikov higher
signatures $\langle \sigma_M, \tau \rangle$ in full generality.  
However, in order to {\it prove} the
homotopy-invariance of $\langle \sigma_M, \tau
\rangle$, we need Assumption 2.) 

From equation \equ(highersig) and
the smooth topological invariance of $\sigma_M$, we obtain a ``Novikov
additivity'' for higher signatures.
\medskip
\nproclaim Corollary [novadd]. Let $\Gamma$ satisfy Assumption 2.
Let $M$ be a closed oriented manifold and let
$F$ be a two-sided hypersurface which separates $M$ into pieces $A$ and $B$.
Let $\nu : M \rightarrow B\Gamma$ be a continuous map, defined up to
homotopy. Let $i : F \rightarrow M$ be the inclusion map and put 
$\nu_0 = \nu \circ i$.
Suppose that $F$ satisfies Assumption 1. If $\dim(M) = 2k+1$, choose
a (stable) Lagrangian subspace $L$ of $\overline{\HH}^k(F; {\cal V}_0)$ and
use $L$ to define $\sigma_A$, and $-L$ to define $\sigma_B$.
Then for any $\tau \in \HH^*(\Gamma; \CC)$, the
corresponding higher signature of $M$ satisfies
$$\const \,
\langle L(TM) \cup \nu^* \tau, [M] \rangle = \langle \sigma_A, \tau \rangle +
\langle \sigma_B, \tau \rangle.$$
\medskip

As a consequence of Corollary \thm[novadd], we obtain a sort of cut-and-paste
invariance of the higher signatures of closed manifolds.
\medskip
\nproclaim Corollary [cutandpaste]. Let $\Gamma$ satisfy Assumption 2. Let
$M_1$ and $M_2$ be closed oriented manifolds, equipped with continuous maps 
$\nu_{M_i} : M_i \rightarrow
B\Gamma$ which are defined up to homotopy. 
Suppose that there are splittings $M_1 = A \cup_F B$ and
$M_2 = A \cup_F B$ over separating two-sided hypersurfaces.
(That is, both $M_1$ and $M_2$ are constructed
by gluing $A$ to $B$, but the gluing diffeomorphisms 
$ \phi_i: \partial A \rightarrow \partial B$ can be different.) 
Suppose that $\nu_1 \big|_{A}$ is homotopic to 
$\nu_2 \big|_A$,
$\nu_1 \big|_{B}$ is homotopic to $\nu_2 \big|_B$ and 
$F$ satisfies Assumption 1.
If $\dim(M_i) = 2k+1$, we also assume that $(\phi_2 \circ \phi_1^{-1})^*$
preserves a (stable) 
Lagrangian subspace of $\overline{\HH}^k(F; {\cal V}_0)$.
Then the higher signatures of $M_1$ and $M_2$ coincide.
\medskip

Corollary \thm[cutandpaste] is relevant because the higher signatures
of closed manifolds are generally not cut-and-paste invariant 
(over $B \Gamma$). 
For example, it is not hard to see that this is the case
when $\Gamma = \ZZ$, using [KKNO, Theorem 1.2], 
and the case $\Gamma = \ZZ^k$ then follows from
[N, Lemma 8].
 This shows that some condition like Assumption 1 
is necessary if one wants to define higher 
signatures for manifolds-with-boundary so as to have Novikov additivity.
Such a situation does not arise for the usual ``lower'' signature.

In general, it seems to be an interesting question as to for which groups
$\Gamma$ and which cohomology classes $\tau \in \HH^*(\Gamma; \CC)$, the
corresponding Novikov higher signature (of closed manifolds) is a
cut-and-paste invariant (over $B\Gamma$);  
see [L5, Remark 4.1] and [R, Chapter 30] for further discussion.

We now give a brief description of the proof of Theorem \thm[mainthm]. 
In the case of closed manifolds, the analytic proofs of the Novikov Conjecture,
as in [CM, Theorem 6.6], consist of two steps. First, one shows that the
index of the signature operator, as an element of $K_*(C^*_r(\Gamma))$, is
an oriented-homotopy invariant. Second, one constructs a pairing of
$K_*(C^*_r(\Gamma))$ with $\HH^*(\Gamma; \CC)$ and one verifies
that the result is the Novikov higher signature. This last step amounts
to proving an index theorem.

In the case of closed manifolds, 
many of the proofs of the first step implicitly
use the cobordism invariance of the index. As even the
usual ``lower'' signatures of manifolds-with-boundary are not cobordism
invariant, this method of proof is ruled out for us.  
Instead, we give a direct
proof of the homotopy invariance which, in the closed case, 
was developed by Hilsum and Skandalis [HS].
To use their methods,
we need $C^*_r(\Gamma)$-Fredholm 
signature operators with $C^*_r(\Gamma)$-compact resolvents.
For this reason, in our case
we would like to cone off the boundary on $M$ to obtain a
conical manifold
$CM$ (deleting the vertex point). If ${\cal V}$ denotes the canonical
flat $C^*_r(\Gamma)$-bundle on $CM$, we would consider the signature operator 
acting on $\Omega^*(CM; {\cal V})$, with its index in $K_*(C^*_r(\Gamma))$. 
We would then extend homotopy equivalences between 
manifolds-with-boundary to
homotopy equivalences between conical manifolds, in order to compare their
indices.
However, as the boundary signature operator ${\cal D}_{\partial M}$ 
may well have
continuous spectrum which goes down to zero (see [L4] for examples), 
there are serious
technical problems in carrying out the
conical analysis. (The paper [LP3] looked at a special case
in which  $\Gamma$ is of the form $\Gamma' \times G$, with $G$ finite, and 
${\cal D}_{\partial M}$ can be made invertible by twisting with a
nontrivial representation of the  finite group $G$.
The corresponding index class was proven
to be an oriented-homotopy invariant using the results of [KM]. The
higher APS-index formula of [LP1] was then applied
in order to show that a twisted version
of \equ(highersig) was an oriented-homotopy invariant.   
However, we wish to deal with the
general case here.)

In order to get around the problem of low-lying spectrum of 
${\cal D}_{\partial M}$, we follow the method of proof sketched in
[L5, Appendix].
We basically add an algebraic complex to cancel out the
small spectrum.  More precisely, we consider a certain Hermitian complex
$\widehat{W}^*$ of finitely-generated projective $C^*_r(\Gamma)$-modules. We
form a new complex
$C^* = \Omega^*(CM; {\cal V}) \oplus \left( 
\Omega^*(0,2) \, \widehat{\otimes} \, \widehat{W}^* \right)$, where the 
algebraic
complex
$\Omega^*(0,2) \, \widehat{\otimes} \, \widehat{W}^*$ is endowed with a metric
which makes it ``conical'' at $0$ and $2$. Formally, the complex  
$\Omega^*(0,2) \, \widehat{\otimes} \, \widehat{W}^*$ has vanishing higher 
signature, and so by adding it we have not changed the putative
higher signature of $CM$.
Then we perturb the differential of $C^*$ in order to couple 
$\Omega^*(CM; {\cal V})$ and $\Omega^*(0,2) \, \widehat{\otimes} \, 
\widehat{W}^*$
near the endpoint $0$. That is, we do a mapping-cone-type construction along
the conical end, which is turned on by a function $\phi(x)$ with 
$\phi(x) = 1$ for $0 < x < 1/4$ and $\phi(x) = 0$ for $1/2 \le x \le 2$.
This mapping-cone-type construction is done in a way which preserves
Poincar\'e duality, and makes the new boundary operator invertible.  
The price to be paid is that the new ``differential'' $D_C$
no longer satisfies $(D_C)^2 = 0$, as $\phi$ is nonconstant.  However,
by increasing the length of the conical end, we can make $(D_C)^2$ 
arbitrarily small in norm.  Then we can apply the ``almost flat''
results of [HS, Theorem 4.2] to conclude that the signature index class 
$[{\cal D}_C^{{\rm conic}}] \in
K_*(C^*_r(\Gamma))$ is an oriented-homotopy invariant.  
The results of [HS, Theorem 4.2] were
designed to deal with the case of almost-flat vector bundles.  We do not
have such vector bundles in our case, but we can use the results of
[HS, Theorem 4.2] nevertheless.

As ${\cal B}^\infty$ is assumed to be a smooth subalgebra of $C^*_r \Gamma$,
there is an isomorphism 
$K_*(C^*_r(\Gamma)) \cong K_*({\cal B}^\infty)$. Hence there is a
Chern character $\ch([{\cal D}^{{\rm conic}}_C] ) \in 
\overline{\HH}_*({\cal B}^\infty)$.
The second main step
in the proof of Theorem \thm[mainthm] consists
of proving an index theorem, in order to show 
that $\ch([{\cal D}^{{\rm conic}}_C])$ 
is given by the right-hand-side of 
\equ(highersig). In principle one could do so within the framework of
analysis on cone manifolds, but this seems to be very difficult.
Instead, we introduce two new
$C^*_r(\Gamma)$-Fredholm signature operators, one being
an Atiyah-Patodi-Singer (APS)-type operator and the other being a 
Melrose $b$-type operator.
We show that both the conic index and 
the $b$-index equal
the APS-index:
$$[{\cal D}_C^{{\rm conic}} ]=[{\cal D}_C^{{\rm APS}}]=
[{\cal D}_C^b]\quad{\rm in}
\quad K_*(C^*_r(\Gamma)).$$
The advantage of this intermediate step
is that we can then compute the Chern character
of the $b$-index $[{\cal D}_C^b]$ by means of
an extension of the higher $b$-pseudodifferential
calculus developed in 
[LP1] and [LP3]. Thus
we (briefly) develop an enlarged
$b$-calculus which takes into account
the above mapping-cone construction, and show that the Chern character of
the $b$-index class is given by the right-hand-side of \equ(highersig). 
This completes the proof of Theorem \thm[mainthm].

The organization of the paper is as follows.  In Section 1 we 
establish our conventions for signature operators,
following [HS, Section 3.1].  
We also give the product decomposition of the signature
operator on a manifold-with-boundary near the boundary. In Section 2 
we review the definition of the higher eta-invariant of an 
odd-dimensional manifold. In Section 3 
we review the definition of the higher eta-invariant of an 
even-dimensional manifold.  In Section 4
we discuss the signature operator on a manifold-with-boundary, perturbed
by the afore-mentioned algebraic complex $\widehat{W}^*$.
In Section 5 we add a conic metric and show that we obtain a well-defined
conic index class in $K_0(C^*_r(\Gamma))$. 
In Section 6 we prove that the conic index class is an
oriented-homotopy invariant. In Section 7 we define the APS-index class
and prove that it equals the conic index class.  In Section 8 we
define the (perturbed) $b$-signature operator.  In Section 9 we
show that the $b$-signature operator has a well-defined index class.
In Section 10 we show that the APS-index class and the
$b$-index class coincide.  In Section 11 we prove an index theorem which
computes the index class of the $b$-signature operator.  In Section 12
we put the pieces together to prove Theorem \thm[mainthm] and Corollaries
\thm[hominv]-\thm[cutandpaste]. In the Appendix we sketch an argument which
relates the signature class considered in this paper to that defined in
[LP4], using symmetric spectral sections.

\medskip
\noindent {\bf Table of Contents:}
\medskip
\noindent {\bf 0. Introduction}

\noindent {\bf 1.  Signature operators }

\noindent {\bf 2.  The higher eta invariant of an odd-dimensional manifold}

\noindent {\bf 3.  The higher eta invariant of an even-dimensional manifold}

\noindent {\bf 4.  Manifolds with boundary: the perturbed signature operator}

\noindent {\bf 5. The conic index class}

\noindent {\bf 6. Homotopy invariance of the conic index class}

\noindent {\bf 7. Equality of the conic and APS-index classes}

\noindent {\bf 8. The enlarged $b$-calculus}

\noindent {\bf 9. The $b$-index class}

\noindent {\bf 10. Equality of the APS  and $b$-index classes}

\noindent {\bf 11. The higher index formula for the $b$-signature operator}

\noindent {\bf 12. Proofs of Theorem \thm[mainthm] and Corollaries
\thm[hominv]-\thm[cutandpaste] }

\noindent {\bf 13. Appendix}

%\fine
\expandafter\ifx\csname sezioniseparate\endcsname\relax%
   \input macro \fi

%--------------------------------------------------------
\numsec=1
\numfor=1
\numtheo=1
\pgn=1

\beginsection 1. Signature Operators

In this section we establish our conventions for signature operators,
following [HS, Section 3.1]. 
The only difference between our conventions and those of
[HS] is that we deal with left modules, whereas [HS] deals with right
modules. We also give the product decomposition of the signature
operator on a manifold-with-boundary near the boundary.

Let $\Lambda$ be a $C^*$-algebra with unit.  Let ${\cal B}^\infty$ be a
Fr\'echet locally $m$-convex
$*$-subalgebra of $\Lambda$ which is dense in $\Lambda$ and
closed under the holomorphic functional calculus in $\Lambda$ 
[Co, Section III.C].

\nproclaim Definition [Hcomplex].
A graded regular $n$-dimensional Hermitian complex consists of \hfil\break
1. A $\ZZ$-graded cochain complex $({\cal E}^*, D)$ of finitely-generated
projective left ${\cal B}^\infty$-modules, \hfil\break
2. A nondegenerate quadratic form 
$Q : {\cal E}^* \times {\cal E}^{n-*} \rightarrow {\cal B}^\infty$ and
\hfil\break
3. An operator $\tau \in Hom_{{\cal B}^\infty} \left( {\cal E}^*,
{\cal E}^{n-*} \right)$ \hfil\break
such that \hfil\break
1. $Q(bx, y) = b Q(x,y)$.\hfil\break
2. $Q(x,y)^* = Q(y,x)$.\hfil\break
3. $Q(Dx, y) + Q(x, Dy) = 0$.\hfil\break
4. $\tau^2 = I$.\hfil\break
5. $<x, y> \equiv Q(x, \tau y)$ defines a Hermitian metric on ${\cal E}$ 
([L3, Definition 7]). \hfil\break

Let $M$ be a closed oriented $n$-dimensional Riemannian
manifold. Let ${\cal V}^\infty$ be a flat
${\cal B}^\infty$-vector bundle on 
$M$, meaning in particular
that its fibers are finitely-generated projective left
${\cal B}^\infty$-modules and the transition functions are compatible with the
${\cal B}^\infty$-module structures. We assume that the fibers of
${\cal V}^\infty$ have 
${\cal B}^\infty$-valued Hermitian inner products which are compatible with
the flat structure. Put ${\cal V} = \Lambda \otimes_{{\cal B}^\infty}
{\cal V}^\infty$. It is a flat vector bundle of ${\Lambda}$-Hilbert modules.

Let $\Omega^*(M; {\cal V}^\infty)$ denote the vector space of smooth
differential forms with coefficients in ${\cal V}^\infty$. If $n=\dim(M) > 0$ then
$\Omega^*(M; {\cal V}^\infty)$ is not finitely-generated over 
${\cal B}^\infty$, but we wish to show that it still has all of
the formal properties of a graded regular $n$-dimensional Hermitian complex.
If $\alpha \in
\Omega^*(M; {\cal V}^\infty)$ is homogeneous, denote its degree by $|\alpha|$.
In what follows, $\alpha$ and $\beta$ will sometimes implicitly denote
homogeneous elements of $\Omega^*(M; {\cal V}^\infty)$.
Given $m\in M$ and $(\lambda_1 \otimes e_1), (\lambda_2 \otimes e_2)
\in \Lambda^*(T^*_mM) \otimes
{\cal V}^\infty_m$, we define 
$(\lambda_1 \otimes e_1) \wedge (\lambda_2 \otimes e_2)^*
\in \Lambda^*(T^*_mM) \otimes
{\cal B}^\infty$ by 
$$(\lambda_1 \otimes e_1) \wedge (\lambda_2 \otimes e_2)^* =
(\lambda_1 \wedge \overline{\lambda_2}) \otimes <e_1, e_2>.$$
Extending by linearity (and antilinearity), given
$\omega_1, \omega_2 \in \Lambda^*(T^*_mM) \otimes
{\cal V}^\infty_m$,
we can define 
$\omega_1 \wedge \omega_2^* \in 
\Lambda^*(T^*_mM) \otimes {\cal B}^\infty$.

Define a ${\cal B}^\infty$-valued
quadratic form $Q$ on $\Omega^*(M; {\cal V}^\infty)$ by
$$
Q(\alpha, \beta) = i^{- |\alpha| (n - |\alpha|)} \int_M \alpha(m) 
\wedge {\beta}(m)^*.
$$
It satisfies $Q(\beta, \alpha) = {Q(\alpha, \beta)}^*$.
Using the Hodge duality operator $*$, define
$\tau : \Omega^p(M; {\cal V}^\infty) \rightarrow 
\Omega^{n-p}(M; {\cal V}^\infty)$ by
$$
\tau (\alpha) = i^{- |\alpha| (n-|\alpha|)} * \alpha.
$$
Then $\tau^2 = 1$ and the inner product $<\cdot,\cdot>$ on 
$\Omega^*(M; {\cal V}^\infty)$ is 
given by
$< \alpha, \beta> = Q(\alpha, \tau \beta)$.
Define $D : \Omega^*(M; {\cal V}^\infty) \rightarrow 
\Omega^{*+1}(M; {\cal V}^\infty)$ by
$$
D \alpha = i^{|\alpha|} d \alpha.
\Eq(signed-diff)
$$
It satisfies $D^2 = 0$.
Its dual $D^\prime$ with respect to $Q$, i.e. the operator $D^\prime$
such that $Q(\alpha, D\beta) = Q(D^\prime \alpha, \beta)$, is given by
$D^\prime = -D$. The formal adjoint of $D$ with respect to $<\cdot, \cdot>$ is
$D^* = \tau D^\prime \tau = - \tau D \tau$.
\nproclaim Definition [sigop].
If $n$ is even, the signature operator is
$$
\Si = D + D^* = D - \tau D \tau.
\Eq(evensign)
$$
It is formally self-adjoint and 
anticommutes with the $\ZZ_2$-grading
operator $\tau$. If $n$ is odd, the signature operator is
$$
\Si = -i (D \tau + \tau D).
\Eq(oddsign)
$$
It is formally self-adjoint.

Let $\Omega^*_{(2)}(M; {\cal V})$ denote the completion of 
$\Omega^*(M; {\cal V})$ in the sense of $\Lambda$-Hilbert modules. 
If $n$ is even then the triple 
$(\Omega^*_{(2)}(M; {\cal V}), Q, D)$ defines an
element of ${\bf L}_{nb}(\Lambda)$ 
in the sense of [HS, D\'efinition 1.5].

Now suppose that $M$ is a compact oriented 
manifold-with-boundary of
dimension $n = 2m$. 
Let $\partial M$ denote
the boundary of $M$. 
We fix a non-negative boundary defining function
$x\in C^\infty(M)$ for $\pa M$ and a Riemannian metric
on $M$ which is isometrically a product in an (open) collar neighbourhood
${\cal{U}}\equiv (0,2)_x\times\pa M$ of $\partial M$.
The signature operator $\Si$ still makes sense
as a differential operator on $\Omega^*(\interior(M); {\cal V}^\infty)$. 
Let ${\cal V}^\infty_0$ 
denote the pullback of
${\cal V}^\infty$ from $M$ to $\partial M$; there is a natural
isomorphism 
$${\cal V}^\infty|_{{\cal{U}}}\cong (0,2) \times {\cal V}^\infty_0 \,.$$
Our orientation conventions are such that the
volume form on $(0,2) \times \partial M$ is $d\vol_M = dx \wedge
d\vol_{\partial M}$. Let $Q_{\partial M}$, $\tau_{\partial M}$, 
$D_{\partial M}$
and $\Si (\partial M)$ 
denote the expressions defined above on 
$\Omega^*(\partial M; {\cal V}^\infty_0)$. 
We wish to decompose $Q$, $\tau$, $D$ and $\Si$, when restricted to
compactly-supported forms on $(0,2) \times \partial M$, in terms of
$Q_{\partial M}$, $\tau_{\partial M}$, $D_{\partial M}$
and $\Si (\partial M)$. 

For notation, we let $\Omega^*_c(0,2)$ denote
compactly-supported forms on $(0,2)$. We let $\otimes$ denote a projective
tensor product and
we let $\, \widehat{\otimes} \,$
denote a graded projective tensor product.
We write a compactly-supported differential form on  
$(0,2) \times \partial M$ as $(1 \wedge \alpha(x)) + (dx \wedge
\beta(x))$, where for each $x \in (0,2)$, 
$\alpha(x)$ and $\beta(x)$ are in 
$\Omega^*(\partial M; {\cal V}^\infty_0)$. 
It is convenient to
introduce the notation
$$
\widehat{\alpha} = i^{|\alpha|} \alpha
$$
for $\alpha \in \Omega^*(\partial M; {\cal V}^\infty_0)$.
One finds
$$
Q (dx \wedge \alpha, 1 \wedge \beta) = \int_0^2
Q_{\partial M}(\alpha(x), \widehat{\beta(x)}) dx, 
$$
$$
Q (1 \wedge \alpha, dx \wedge \beta) = \int_0^2
Q_{\partial M}(\widehat{\alpha(x)}, \beta(x)) dx, 
$$
$$
\tau(1 \wedge \alpha) = dx \wedge \tau_{\partial M} \widehat{\alpha},
\Eq(Top)$$
$$
\tau(dx \wedge \alpha) = 1 \wedge
i^{-(2m-1)} \tau_{\partial M} \widehat{\alpha},
$$
$$
D(1 \wedge \alpha) = (1 \wedge D_{\partial M} \alpha) + 
(dx \wedge \partial_x \widehat{\alpha}),
$$
$$
D(dx \wedge \alpha) = dx \wedge -i D_{\partial M} \alpha.
$$
Then one can compute that $\Si$ takes the form
$$
\Si = \left(
\matrix{
D_{\partial M} - \tau_{\partial M} D_{\partial M} \tau_{\partial M} & 
- i^{- |\beta|} \partial_x \cr
i^{|\alpha|} \partial_x & -i(
D_{\partial M} + \tau_{\partial M} D_{\partial M} \tau_{\partial M})
\cr}
\right)
$$
when acting on
$\left( 
\matrix{
1 \wedge \alpha \cr
dx \wedge \beta \cr}
\right)$.
That is,
$$
\Si (1 \wedge \alpha) =
(1 \wedge (D_{\partial M} - \tau_{\partial M} D_{\partial M} 
\tau_{\partial M}) 
\alpha)
+ (dx \wedge i^{|\alpha|} \partial_x \alpha)
$$ 
and
$$
\Si (dx \wedge \beta) =
(1 \wedge - i^{-|\beta|} \partial_x \beta) 
+ (dx \wedge -i (
D_{\partial M} + \tau_{\partial M} D_{\partial M} \tau_{\partial M}) \beta).
$$
Let us define an operator 
$\Theta: \Omega^*_c((0,2) \times\pa M;{\cal{V}}^\infty_0)\rightarrow 
\Omega^*_c((0,2)\times\pa M;{\cal{V}}^\infty_0)$ by
$$
\Theta((1 \wedge \alpha) + (dx \wedge \beta)) =
(1 \wedge -i^{-\beta} \beta) + (dx \wedge i^{|\alpha|} \alpha).$$
Then $\Theta$
anticommutes with $\tau$ and we can write ${\cal D}^{\rm sign} =
\Theta (\partial_x + H)$, where $H$ commutes with $\tau$.
Acting on the $+1$-eigenvector
$$(dx \wedge \alpha) + \tau (dx \wedge \alpha) =
(dx \wedge \alpha) + (1 \wedge i^{-(2m-1-|\alpha|)} \tau_{\partial M} \alpha)$$
of $\tau$, one finds

$$H \left( (dx \wedge \alpha) + \tau (dx \wedge \alpha) \right) =
(dx \wedge -i(D_{\partial M} \tau_{\partial M} + \tau_{\partial M}
D_{\partial M})\alpha) + (1 \wedge -i^{-|\alpha|} (D_{\partial M} -
\tau_{\partial M} D_{\partial M} \tau_{\partial M}) \alpha).$$

Let $E^{\pm}$ be the $\pm 1$-eigenspaces of $\tau$ acting on
$\Omega^*_c((0,2) \times \partial M; {\cal V}^\infty_0)$.
We define  an  isomorphism $\Phi$
from 
$C^\infty_c(0,2) \otimes\Omega^*(\partial M; {\cal V}^\infty_0)$ to $E^+$, 
by setting 
$$
\Phi(\alpha) = (dx \wedge \alpha) + \tau (dx \wedge \alpha).
$$
We 
then obtain an isomorphism $$\Theta\circ \Phi: 
C^\infty_c(0,2)\otimes\Omega^*(\partial M; {\cal V}^\infty_0)\rightarrow 
E^-\;.$$
Denote as usual by $\Si_+$ the signature operator on $M$ going from
$E^+$ to $E^-$;
using the above isomorphisms we easily obtain
$$\Phi^{-1}\circ H|_{E^+} \circ \Phi=\Si(\pa M)$$
and
$$
\Si_+ =\Theta \circ \Phi ( \pa_x + \Si(\pa M) ) \Phi^{-1}.
\Eq(deco)
$$
This shows that $\Si(\pa M)$ is the boundary component of
$\Si$ in the sense of Atiyah-Patodi-Singer [APS, (3.1)].

\smallskip

Consider the $\ZZ_2$-graded vector space
$(\Omega^*_{(2)}(\partial M; {\cal V}_0)) 
\oplus (dx \wedge \Omega^*_{(2)}(\partial M; {\cal V}
_0))$, where the $\ZZ_2$-grading comes from the
operator $\tau$ of \equ(Top). 
The triple $\left( (\Omega^*_{(2)}(\partial M; {\cal V}_0)) 
\oplus (dx \wedge
\Omega^*_{(2)}(\partial M; {\cal V}_0)), Q, \Theta H \right)$ 
defines an element of
${\bf L}_{nb,odd}(\Lambda)$ in the sense of [HS, p. 81].

%\fine
\expandafter\ifx\csname sezioniseparate\endcsname\relax%
    \input macro \fi

%--------------------------------------------------------
\numsec=2
\numfor=1
\numtheo=1
\pgn=1

\beginsection 2. The higher eta invariant of an odd-dimensional manifold

In this section we review the definition of the
higher eta invariant 
([L2, Definition 11] and 
[L5, Section 3.2]). The material in this section comes from
these references, with minor variations.  The higher eta invariant 
is defined for closed oriented Riemannian manifolds
of either even or odd dimension. We first treat
the case of a closed oriented Riemannian manifold $F$ of
dimension $n=2m-1$.

Let us make a general remark about homotopy equivalences between
cochain complexes.  Suppose that $(C_1, d_1)$ and 
$(C_2, d_2)$ are cochain complexes,
with homotopy equivalences $f : C_1^* \rightarrow C_2^*$ and $g : C_2^*
\rightarrow C_1^*$. Then one implicitly understands that there are
maps $A : C_1^* \rightarrow C_1^{*-1}$ and $B : C_2^* \rightarrow
C_2^{*-1}$ so that
$I - gf = d_1 A + A d_1$ and $I - fg = d_2 B + B d_2$. It follows that
$gB-Ag : C_2^* \rightarrow C_1^{*-1}$ and 
$fA-Bf : C_1^* \rightarrow C_2^{*-1}$ are cochain maps.
We will say that such $f$ and $g$ form
a {\it double homotopy equivalence} if, in addition, there are maps
$\alpha : C_2^* \rightarrow C_1^{*-2}$ and 
$\beta : C_1^* \rightarrow C_2^{*-2}$ such that
$gB-Ag = d_1 \alpha - \alpha d_2$ and $fA-Bf = d_2 \beta - \beta d_1$.
One can check that the composition of two double 
homotopy equivalences is a double homotopy equivalence.
The notion of double homotopy equivalence is not strictly needed 
for this section
but will enter in the proof of Theorem \thf[conichominv].

Now let $\Gamma$ be a finitely-generated discrete group.
Let $\nu : F \rightarrow B\Gamma$ be a continuous map.  There is a
corresponding normal $\Gamma$-cover $F^\prime\rightarrow F$.
Let $C^*_r(\Gamma)$ be the reduced group $C^*$-algebra of $\Gamma$. 
Let ${\cal B}^\infty$ 
be a subalgebra of $C^*_r(\Gamma)$ as in Section 1.

We introduce two flat unitary
vector bundles of left modules on $F$ :
$$
\Fl=C^*_r(\Gamma) \times_\Gamma F^\prime,
%#New macros: $\Fl={\cal V}$ for ``flat''
\quad\quad
\Fli=\Bi\times_\Gamma F^\prime.
$$
%#New macros: $\Fli=\Fl^\infty$ for ``flat'' infinity
%and $\Bi={\cal B}^\infty$
Following [L2, Section 4.7], 
we make an assumption about the de Rham cohomology of $F$, with
value in the local system $\Fl$.
\medskip

\noindent
{\bf Assumption 1 :} The natural surjection
$\HH^{m}(F;\Fl) \rightarrow \overline{\HH}^{m}(F;\Fl)$ is an
isomorphism.
\medskip

\nproclaim Lemma [Riem].
If $F$ is equipped with a Riemannian metric then
Assumption 1 is equivalent to saying that the differential form Laplacian
on $\Omega^{m-1}(F^\prime)/\Ker(d)$ has a strictly positive spectrum.

\noindent
{\it Proof.} We give an outline of the proof.
Let $\Omega^*_{(2)}(F; \Fl)$ denote the completion of 
$\Omega^*(F; \Fl)$ as a $C^*_r(\Gamma)$-Hilbert module.
Assumption 1 is equivalent to saying that the differential
$\widetilde{D}_F : 
\Omega^{m-1}(F; \Fl) \rightarrow \Omega^{m}(F; \Fl)$ has a closed
image.
Using Hodge duality, this is equivalent to saying that
$\widetilde{D}_F^* : 
\Omega^{m}(F; \Fl) \rightarrow \Omega^{m-1}(F; \Fl)$ has a closed image.
Clearly $\widetilde{D}_F^*$ is adjointable, and we obtain an
orthogonal decomposition $\Omega^{m-1}(F; \Fl) =
\IM (\widetilde{D}_F^*) \oplus \Ker(\widetilde{D}_F)$.
From arguments as in [L3, Propositions 10 and 27], this is
equivalent to saying that $\widetilde{D}_F^* \widetilde{D}_F$ 
has a strictly positive spectrum as a densely-defined operator on
$\IM \left( \widetilde{D}_F^* : 
\Omega^{m}(F; \Fl) \rightarrow \Omega^{m-1}(F; \Fl) \right) \subset
\Omega^{m-1}_{(2)}(F; \Fl)/\Ker(\widetilde{D}_F)$; see
[We, Theorem 15.3.8] for the analogous result in the case of bounded
operators. 
Put ${\cal V}^{(2)} = l^2(\Gamma) \times_\Gamma F^\prime$.
Using the injective homomorphism $C^*_r(\Gamma) \rightarrow B(l^2(\Gamma))$ of
$C^*$-algebras, it follows as in [L3, Proposition 19] that the spectrum of
$\widetilde{D}_F^* \widetilde{D}_F$ acting on 
$\Omega^{m-1}_{(2)}(F; \Fl)/\Ker(\widetilde{D}_F)$ is the same as the
spectrum of  $d^* d$ acting on
$\Omega^{m-1}_{(2)}(F; {\cal V}^{(2)})/\Ker(d)$, where
$\Omega^{m-1}_{(2)}(F; {\cal V}^{(2)})$ denotes the Hilbert space of
square-integrable ${\cal V}^{(2)}$-valued $(m-1)$-forms on $F$. 
However, the latter is the same as the space 
$\Omega^{m-1}_{(2)}(F^\prime)$ of 
square-integrable $(m-1)$-forms on $F^\prime$. 
Since the Laplacian $d^* d + d d^*$ acts on
$\Omega^{m-1}_{(2)}(F^\prime)/\Ker(d)$ as $d^* d$, the lemma follows.  
\QED
\medskip

\nproclaim Lemma [examples]. \hfil\break
(a) If $F$ has a cellular decomposition without any cells of
dimension $m$ then Assumption 1 is satisfied. \hfil\break
\noindent
(b) If $\Gamma$ is finite then 
Assumption
1 is satisfied. \hfil\break
\noindent
(c) If $\dim(F) = 3$, $F$ is connected and $\Gamma = \pi_1(F)$ then, assuming
Thurston's geometrization conjecture, Assumption 1
is satisfied if and only if $F$ is a connected sum of
spherical space forms, $S^1 \times S^2$'s and twisted circle bundles
$S^1 \times_{\ZZ_2} S^2$ over
$\RR P^2$. 

\noindent
{\it Proof.} \hfil\break
(a) If $F$ has a cellular decomposition without any cells of
dimension $m$ then $\HH^{m}(F;\Fl)$ vanishes and Assumption 1 is
automatically satisfied. \hfil\break
(b) If $\Gamma$ is finite then $F^\prime$ is compact and from standard
elliptic theory, the result
of Lemma \thm[Riem] is satisfied. \hfil\break
(c) If $\Gamma = \pi_1(F)$ and $F$ is connected then in the notation of [LL],
the result of Lemma \thm[Riem] is equivalent to
saying that the $m$-th Novikov-Shubin invariant $\alpha_m(F)$ of $F$ is
$\infty^+$. In the present case, $m = 2$.
Let $F = F_1 \sharp F_2 \sharp \ldots \sharp F_N$ be the
connected sum decomposition of $F$ into prime $3$-manifolds.  
From [LL, Proposition 3.7.3], $\alpha_2(F) = \min_{i} \alpha_2(F_i)$.
Hence it suffices to characterize the prime closed $3$-manifolds $F$ with
$\alpha_2(F) = \infty^+$. If $F$ has finite fundamental group then
$\alpha_2(F) = \infty^+$ and the geometrization conjecture says that
$F$ is a spherical space form.  If $F$ has infinite fundamental group 
and $\alpha_2(F) = \infty^+$ then, assuming the geometrization conjecture,
[LL, Theorem 0.1.5] implies that
$F$ has an $\RR^3$, $S^2 \times \RR$ or $Sol$ structure.  
From [LL, Theorem 0.1.4], if $F$ has an
$\RR^3$-structure then $\alpha_2(F) = 3$, while if $F$ has an
$S^2 \times \RR$ structure then $\alpha_2(F) = \infty^+$. Finally,
a slight refinement of [L4, Corollary 5] shows that if $F$ has a $Sol$
structure then $\alpha_2(F) < \infty^+$. The claim follows.
\QED
\medskip

Hereafter we assume that Assumption 1 is satisfied.

\nproclaim Lemma [complex].
There
is a cochain complex $W^* = \bigoplus_{i=0}^{2m-1} W^i$
of finitely-generated projective $\Bi$-modules
such that \hfil\break
\noindent
1. $W^*$ is a graded regular $n$-dimensional Hermitian
complex.  \hfil\break
\noindent
2. The differential $D_W : W^{m-1} \rightarrow W^m$ vanishes. \hfil\break
\noindent
3. There is a double homotopy equivalence
$$f : \Omega^*(F;\Fli) \rightarrow W^* 
\Eq(pr2)$$
which, as an element of $(\Omega^*(F;\Fli))^* \otimes W^*$, is actually
smooth with respect to $F$.

\noindent
{\it Proof.}
The strategy of the proof is to first establish a double homotopy equivalence
between $\Omega^*(F;\Fli)$ and a simplicial cochain complex, and then 
to further homotope the simplicial cochain complex in order to end up with
a graded regular Hermitian complex.
We will implicitly use results from [L3, Proposition 10 and
Section 6.1] concerning spectral analysis involving $\Bi$.
Let $K$ be a triangulation of $F$. Let
$(C^*(K; \Fli), D_K)$ 
be the simplicial cochain complex, a complex of 
finitely-generated free $\Bi$-modules. We first construct a cochain embedding
from $C^*(K; \Fli)$ to $\Omega^*(F;\Fli)$, following the work of Whitney
[Wh, Chapter IV.27]. 
In order to have an embedding into smooth forms, we use the
modification of Whitney's formula given in [D, (3.4)]. If
${\cal V}^{(2)} = l^2(\Gamma) \times_\Gamma F^\prime$ then
the map $W$ of [D, (3.4)], which is defined in the $l^2$-setting,
gives a cochain embedding
$W : C^*(K; {\cal V}^{(2)}) \rightarrow
\Omega^*(F; {\cal V}^{(2)})$ which is a homotopy equivalence. Using the
same formula as in [D, (3.4)], but considering cochains and forms 
with values in the flat bundle $\Fli$, we obtain a
cochain embedding $w : C^*(K; \Fli) \rightarrow
\Omega^*(F;\Fli)$ which is a homotopy equivalence.
Give 
$C^*(K; \Fli)$ the induced ${\cal B}^\infty$-valued Hermitian inner product.

Using this embedding, let us decompose
$\Omega^*(F;\Fli)$ as $\Omega^*(F;\Fli) = C^*(K; \Fli) \oplus
C^\prime$ where $C^\prime$ is the orthogonal complement to 
the finitely-generated free submodule $C^*(K; \Fli)$ of
$\Omega^*(F;\Fli)$.
With respect to this decomposition, we can write $w = \left(
\matrix{
I \cr
0 \cr}
\right)$ and
$D_F = 
\left(
\matrix{
D_K & X \cr
0 & D_{C^\prime}
\cr}
\right)$ for some cochain map
$X \in \Hom_{{\cal B}^\infty} \left( (C^\prime)^*, C^{*+1}(K; \Fli) \right)$.
Then the complex $(C^\prime, D_{C^\prime})$ is acyclic.
Putting 
$\widetilde{w} = \id_{C^*_r(\Gamma)} \otimes_{{\cal B}^\infty} w :
C^*(K; \Fl) \rightarrow
\Omega^*(F;\Fl)$ 
and doing the analogous constructions, we see that
the complex $C^*_r(\Gamma) \otimes_{{\cal B}^\infty} C^\prime$ is also
acyclic.

As $C^\prime$ is acyclic, there is a operator $\delta_{C^\prime}$
of degree $-1$ such that
$\left( \delta_{C^\prime} \right)^2 = 0$ and
$D_{C^\prime} \delta_{C^\prime} + \delta_{C^\prime} D_{C^\prime} = I$.
We claim that we can take $\delta_{C^\prime}$ to be continuous. To see this,
put 
$E_F = 
\left(
\matrix{
0 & 0 \cr
0 & D_{C^\prime}
\cr}
\right)$. It is an element of the space 
$\Psi^1_{{\cal B}^\infty}(F; \Lambda^*(TF) \otimes \Fli,
\Lambda^*(TF) \otimes \Fli)$ of
${\cal B}^\infty$-pseudodifferential operators
of order $1$, as defined in [L3, Section 6.1].
Put ${\cal L} = E_F (E_F)^* + (E_F)^* E_F$, an element 
of $\Psi^2_{{\cal B}^\infty}(F; \Lambda^*(TF) \otimes \Fli,
\Lambda^*(TF) \otimes \Fli)$.

As $C^*_r(\Gamma) \otimes_{{\cal B}^\infty} C^\prime$ is acyclic, its
differentials have closed image. Hence 
$\id_{C^*_r(\Gamma)} \otimes_{{\cal B}^\infty}
E_F : \Omega^*(F; \Fl) \rightarrow 
\Omega^{*+1}(F; \Fl)$ 
also has closed image.
It follows, as in the proof of 
[L3, Propositions 10 and 27], that $0$ is isolated in the spectrum of
${\cal L}$, with $\Ker({\cal L}) = C^{*}(K; \Fli) \oplus 0$.
Let 
$G \in \Psi^{-2}_{{\cal B}^\infty}(F; \Lambda^*(TF) \otimes \Fli,
\Lambda^*(TF) \otimes \Fli)$ 
be the Green's operator for ${\cal L}$.
Then $(E_F)^* G$ is an element of
$\Psi^{-1}_{{\cal B}^\infty}(F; \Lambda^*(TF) \otimes \Fli,
\Lambda^*(TF) \otimes \Fli)$ 
and can be written in the form
$(E_F)^* G = 
\left(
\matrix{
0 & 0 \cr
0 & \delta_{C^\prime}
\cr}
\right)$. As in usual Hodge theory,
this operator $\delta_{C^\prime}$ satisfies 
$\left( \delta_{C^\prime} \right)^2 = 0$ and 
$D_{C^\prime} \delta_{C^\prime} + \delta_{C^\prime} D_{C^\prime} = I$.
It is everywhere-defined and continuous, as 
$(E_F)^* G$ has order $-1$ in the pseudodifferential operator
calculus.

Define $q : \Omega^*(F;\Fli) \rightarrow C^*(K; \Fli)$ by
$q = \left( 
\matrix{
I & -X\delta_{C^\prime} \cr}
\right)
$. Then one can
check that $I - q w = 0$ and $I - w q = D_F A + A D_F$, where
$$A =  
\left(
\matrix{
0 & 0 \cr
0 & \delta_{C^\prime}
\cr}
\right).$$
Furthermore, $qA = Aw = 0$. Hence $w$ and $q$ define a
double homotopy equivalence between $\Omega^*(F;\Fli)$ and $C^*(K; \Fli)$.

We now show that $C^*(K; \Fli)$
is double homotopy equivalent to an appropriate 
regular Hermitian complex $W^*$ of
finitely-generated projective $\Bi$-modules. In the even case, the
homotopy equivalence to a regular Hermitian complex was
proven in [KM, Proposition 2.4]. In order to extend the proof to the odd
case, we need Assumption 1. 

Let $\widetilde{D}_K$ denote the differential on $C^*(K; \Fl)$, an
adjointable operator.
Using the homotopy equivalence between $C^*(K; \Fl)$
and $\Omega^*(F; \Fl)$, along with the fact that all of the maps involved in
defining the homotopy equivalence
are continuous, it follows that
Assumption 1 is equivalent to saying that the natural surjection
$\HH^{m}(K; \Fl) \rightarrow \overline{\HH}^{m}(K; \Fl)$ 
is an isomorphism. Equivalently, $\widetilde{D}_K (C^{m-1}(K; \Fl))$ is
closed in $C^m(K; \Fl)$.
Let $C^m(K; \Fl) = \IM(\widetilde{D}_K) \oplus \Ker(\widetilde{D}_K^*)$ be
the corresponding orthogonal decomposition [We, Theorem 15.3.8].
Then the operator $\widetilde{D}_K \widetilde{D}_K^*$ is invertible
on $\IM(\widetilde{D}_K) \subset C^m(K; \Fl)$ [We, Theorem 15.3.8].
In particular, there is some $\epsilon > 0$ such that the intersection of
the spectrum of $\widetilde{D}_K \widetilde{D}_K^*$ (acting on
$C^m(K; \Fl)$) with the ball
$B_\epsilon(0) \subset \CC$ consists at most of the point $0$.
From [L3, Lemma 1], the same is true of the operator
$D_K D_K^*$, acting on $C^m(K; \Fli)$. 
Define a continuous operator $G$ on $C^m(K; \Fli)$ by
$$G = {1 \over 2 \pi i} \int_\gamma \quad {1 \over \lambda} \quad 
{d\lambda \over {D_K D_K^* - \lambda}},$$
where $\gamma$ is the circle of radius ${\epsilon \over 2}$ around $0 \in \CC$,
oriented counterclockwise. Then $G$ is the Green's operator for
$D_K D_K^*$. Put $\widetilde{G} = \id_{C^*_r(\Gamma)} \otimes_{{\cal B}^\infty}
G$, the Green's operator for $\widetilde{D}_K \widetilde{D}_K^*$.

We claim that ${D_K} (C^{m-1}(K; \Fli))$ is closed in $C^m(K; \Fli)$.
To see this, suppose that $\{z_i\}_{i=1}^\infty$ is a sequence in 
$C^{m-1}(K; \Fli)$ such that $\lim_{i \rightarrow \infty} D_K (z_i) = y$ for
some $y \in C^m(K; \Fli)$. 
Let $\widetilde{z}_i \in C^{m-1}(K; \Fl)$ and 
$\widetilde{y} \in C^m(K; \Fl)$ be the corresponding elements.
Then 
$$\widetilde{D}_K \widetilde{D}_K^* \widetilde{G} (\widetilde{y}) =
\lim_{i \rightarrow \infty} \widetilde{D}_K \widetilde{D}_K^* 
\widetilde{G} \widetilde{D}_K (\widetilde{z}_i) =
\lim_{i \rightarrow \infty} \widetilde{D}_K (\widetilde{z}_i) = 
\widetilde{y}.$$
It follows that
$D_K D_K^* G (y) = y$,
showing that $y \in \IM(D_K)$. 
Equivalently, the surjection
$\HH^{m}(K;\Fli) \rightarrow \overline{\HH}^{m}(K;\Fli)$ is an
isomorphism.
Similarly, using the fact that $D_K^* G D_K$ acts as the identity on
$\IM(D_K^*) \subset C^{m-1}(K; \Fli)$, one can show that
$\IM(D_K^*)$ is closed in $C^{m-1}(K; \Fli)$

We recall that 
$(C^*(K; \Fli), D_K)$ is a Hermitian complex. This means that it has a 
possibly-degenerate quadratic form
$Q_K$ which satisfies conditions $1$.- $3$. of Definition \thm[Hcomplex] and
for which the corresponding map 
$\Phi_K : C^*(K; \Fli) \rightarrow (C^{2m-1-*}(K; \Fli))^\prime$ is a 
homotopy equivalence.  (Here $\prime$ denotes the antidual space.)
{\it A priori}, $\Phi_K$ may not be an isomorphism. 
To construct the regular Hermitian complex $W^*$, we
need to homotope $C^*(K; \Fli)$ so that the map $\Phi_K$ 
becomes an isomorphism.
Using the construction of [Lu, Proposition 1.3], we can construct a Hermitian
complex $Z^*$ which is homotopy
equivalent to $C^*(K; \Fli)$ and whose map
$\Phi_Z : Z^* \rightarrow (Z^{2m-1-*})^\prime$ is an isomorphism in degrees
other than $m-1$ and $m$.
Looking at
the diagram in the proof of [Lu, Proposition 1.3], one sees that
$C^*(K; \Fli)$ is in fact
double homotopy equivalent to $Z^*$. We again have that
the surjection
$\HH^{m}(Z) \rightarrow \overline{\HH}^{m}(Z)$ is an
isomorphism, or equivalently, $D_Z (Z^{m-1})$ is
closed in $Z^m$. From the diagram in
the proof of [Lu, Proposition 1.3], there is an obvious 
${\cal B}^\infty$-valued Hermitian inner product on $Z^*$, and 
$\widetilde{D}_Z =
\id_{C^*_r(\Gamma)} \otimes_{{\cal B}^\infty} D_Z$ is adjointable.

Put 
$$
W^i =
\cases{
Z^i & if $i < m-1$, \cr
\Ker \left(D_Z : Z^{m-1} \longrightarrow Z^{m} \right) &
if $i = m-1$, \cr
Z^m/\IM \left(D_Z : Z^{m-1} \longrightarrow 
Z^{m} \right)
& if $i = m$, \cr
Z^{i} & if $i > m$. 
} \Eq(Vi)
$$
We give $W^*$ the differential induced from $Z^*$ 
in degrees other than $m-1$, and
the zero differential in degree $m - 1$. 

Using the fact that $D_Z (Z^{m-1})$ is
closed in $Z^m$, it follows as before that
there is some $\epsilon > 0$ such that the intersection of
the spectrum of $D_Z D_Z^*$ (acting on $Z^m$) with the ball
$B_\epsilon(0) \subset \CC$ consists at most of the point $0$.
Then with $\gamma$ as before, the projection operator
$
{1 \over 2 \pi i} \int_\gamma 
{d\lambda \over {\lambda - D_K D_K^*}}
$
gives a direct sum decomposition into closed ${\cal B}^\infty$-submodules :
$$Z^m = \IM \left(D_Z : Z^{m-1} \longrightarrow 
Z^{m} \right) \oplus \Ker
\left(D_Z^* : Z^{m} \longrightarrow 
Z^{m-1} \right).
$$
Using this decomposition, we can identify $W^m$ with $\Ker
\left(D_Z^* : Z^{m} \longrightarrow 
Z^{m-1} \right)$.
It also follows as before that
$D_Z^*(Z^m)$ is closed in $Z^{m-1}$,
and there is a
direct sum decomposition into closed ${\cal B}^\infty$-submodules :
$$
Z^{m-1} = \Ker \left(D_Z : Z^{m-1} 
\longrightarrow Z^{m} \right) \oplus 
\IM \left(D_Z^* : Z^{m} \longrightarrow 
Z^{m-1} \right).
$$
Let 
$p : Z^* \longrightarrow W^*$ be the corresponding projection operator and let
$i : W^* \longrightarrow Z^*$ be the inclusion operator.
Let $L : Z^* \rightarrow Z^{*-1}$ be the map which is an inverse to 
$$D_Z :  \IM \left(D_Z^* : Z^{m} \longrightarrow 
Z^{m-1} \right) \rightarrow
\IM \left(D_Z : Z^{m-1} \longrightarrow 
Z^{m} \right)$$ on $\IM \left(D_Z : Z^{m-1} \longrightarrow 
Z^{m} \right) \subset Z^m$, i.e. $L = D_Z^* 
\left( D_Z D_Z^* \big|_{\IM(D_Z)} \right)^{-1}$, and 
which vanishes on 
\hfil\break
$\Ker \left(D_Z^* : Z^{m} \longrightarrow 
Z^{m-1} \right) \subset Z^m$ and on the rest of $Z^*$.
Then one can check that $p$ and $i$ are cochain maps, that
$p \circ i = I$ and that $i \circ p = I - D_ZL - LD_Z$. Also, $Li= pL = 0$.
Thus $Z^*$ and $W^*$ are doubly homotopy
equivalent.  We give $W^*$ the structure of a Hermitian complex by saying
that $\Phi_W = i^\prime \circ \Phi_Z \circ i$. 
Equivalently, $Q_W$ is
the quadratic form induced from $Q_Z$ under $i : W^* \rightarrow Z^*$. 
Then $W^*$ and $Z^*$ are homotopy equivalent as Hermitian complexes.

We claim that $\Phi_W$ is an isomorphism. This is clear when $\Phi_W$ acts on
$W^*$, $* \notin \{m-1, m\}$, as $\Phi_Z$ is an isomorphism in those degrees.
Hence we must prove the following result :

\nproclaim Sublemma [sublemma].
Suppose that we have a homotopy equivalence $\Phi^* : W^* \rightarrow
(W^{2m-1-*})^\prime$ 
$$\matrix{
\ldots & \rightarrow & W^{m-2} & \rightarrow & 
W^{m-1} & \rightarrow & W^m & \rightarrow & W^{m+1} & \rightarrow
& \ldots \cr
&  &  \downarrow & 
&  \downarrow &  & \downarrow & &
\downarrow & &  \cr
\ldots & \rightarrow & 
(W^{m+1})^\prime & 
\rightarrow & (W^{m})^\prime & \rightarrow & (W^{m-1})^\prime & \rightarrow & 
(W^{m-2})^\prime & \rightarrow
& \ldots \cr}
$$
such that $\Phi^*$ is an isomorphism for $* \notin \{m-1, m\}$ and
$D_W^{m-1} : W^{m-1} \rightarrow W^m$ vanishes.  Then $\Phi^{m-1}$ and
$\Phi^m$ are isomorphisms.

\noindent
{\it Proof.}
We first show that $\Phi^{m-1}$ is injective.
Suppose that $x \in W^{m-1}$ and $\Phi^{m-1}(x) = 0$. 
As $\Phi^{m-1}$ is an isomorphism on 
cohomology, and $[\Phi^{m-1}(x)]$ vanishes in cohomology, there is a
$y \in W^{m-2}$ such that $x = D_W^{m-2} y$. Then $(D_W^m)^\prime 
(\Phi^{m-2}(y)) = \Phi^{m-1}(D_W^{m-2} y) = 0$. Hence $[\Phi^{m-2}(y)]$ 
represents a cohomology class and, as $\Phi^{m-2}$ is an isomorphism on
cohomology, there are some $z \in W^{m-2}$ and $u \in (W^{m+2})^\prime$
such that $D_W^{m-2} z = 0$ and
$\Phi^{m-2}(y) - \Phi^{m-2}(z) = (D_W^{m+1})^\prime(u)$. Put 
$v = (\Phi^{m-3})^{-1}(u)$. Then
$\Phi^{m-2}(y-z-D_W^{m-3} v) = 
(D_W^{m+1})^\prime(u) - \Phi^{m-2}(D_W^{m-3}v) =
(D_W^{m+1})^\prime(u) - (D_W^{m+1})^\prime(\Phi^{m-3}(v)) = 0$.
Thus $y-z-D_W^{m-3}v = 0$ and $x = D_W^{m-2}y = 
D_W^{m-2}(z + D_W^{m-3}v) = 0$,
which shows that $\Phi^{m-1}$ is injective.

We now show that $\Phi^m$ is injective.  Suppose that 
$x \in W^m$ and $\Phi^m(x) = 0$.
Then $\Phi^{m+1}(D_W^m x) = (D_W^{m-2})^\prime(\Phi^m(x)) = 0$, so
$D_W^m x = 0$. Thus $x$ represents a cohomology class.
As $\Phi^{m}$ is an isomorphism on 
cohomology, and $[\Phi^{m}(x)]$ vanishes in cohomology, 
it follows that $x \in \IM(D_W^{m-1}) = 0$.
This shows that
$\Phi^m$ is injective.

We now show that $\Phi^{m-1}$ is surjective.  Suppose that $x \in 
(W^m)^\prime$. As $(D_W^{m-1})^\prime(x) = 0$, there is a cohomology class
represented by $[x]$. As $\Phi^{m-1}$ is an isomorphism on cohomology,
there are some $y \in W^{m-1}$ and $z \in (W^{m+1})^\prime$ such that
$x = \Phi^{m-1}(y) + (D_W^m)^\prime(z)$. Put $w = (\Phi^{m-2})^{-1}(z)$.
Then $x = \Phi^{m-1}(y) + (D_W^m)^\prime(\Phi^{m-2}(w)) =
\Phi^{m-1}(y) + \Phi^{m-1}(D_W^{m-2} w) = \Phi^{m-1}(y + D_W^{m-2} w)$,
which shows that $\Phi^{m-1}$ is surjective.

We finally show that $\Phi^m$ is surjective.  Suppose that
$x \in (W^{m-1})^\prime$. Put $y = (\Phi^{m+1})^{-1}((D_W^{m-2})^\prime(x))$.
As $[\Phi^{m+1}(y)]$ vanishes in cohomology, and
$\Phi^{m+1}$ is an isomorphism on cohomology, there is some $z \in 
W^m$ such that $y = D_W^m z$. Then $(D_W^{m-2})^\prime (x - \Phi^m(z)) =
\Phi^{m+1}( D_W^m z) - (D_W^{m-2})^\prime (\Phi^m(z)) = 0$. Thus
$x - \Phi^m(z)$ represents a cohomology class.  As $\Phi^m$ is an isomorphism
on cohomology, there is some $w \in W^m$ such that $D_W^m w = 0$ and
$x - \Phi^m(z) = \Phi^m(w)$. Hence $x = \Phi^m(z + w)$, which proves the
sublemma.
\QED
\medskip

To finish the proof of Lemma \thm[complex],
as in [KM, Proposition 2.6], 
one can introduce a grading $\tau_W$ so that $(W^*, Q_W, \tau_W)$ satisfies
Definition \thm[Hcomplex].
Hence we have constructed the desired complex
$W^*$, along with a double homotopy equivalence 
$f : \Omega^*(F; \Fli) \rightarrow W^*$. 
From this,
$f$ is an element of $(\Omega^*(F;\Fli))^* \otimes_{{\cal B}^\infty} W^*$.
{\it A priori}, it could be distributional with respect to $F$.  
However, in the proof we constructed $f$ to actually
be smooth on $F$, 
i.e. $f\in C^{\infty}(F;
\Hom_{\Bi}(\Lambda ^* TF \otimes {\cal V}^{\infty}, W^*))$. 
The lemma follows.    \QED

\medskip

Following [L5, (3.23)], we define a new $n$-dimensional complex
$\whW^*$ %#New macro: $\whW$ for $\widehat{W}$
by 
$$\whW^i=
\cases{
W^{i+1}& if $-1\leq i \leq m-2$, \cr
0 & if $i = m-1$ or $m$, \cr
W^{i-1}& if $m+1 \leq i \leq 2m$. \cr} \Eq(newcomplex)$$
The differential $D_W$ induces
a differential $D_{\whW}$ in an obvious way.
We also obtain 
a Hermitian form
$Q_{\whW}(\cdot, \cdot)$ on $\whW$ by putting
$$Q_{\whW}(v^j,z^{(2m-1)-j})
= Q_W(v^j,z^{(2m-1)-j})$$
for $v^j\in\whW^j$, $z^{(2m-1)-j}\in\whW^{(2m-1)-j}$, and
a duality operator
$\tau_{\whW}:\whW^j\rightarrow \whW^{(2m-1)-j}$ by putting
$$\tau_{\whW}(v^j)=\tau_W(v^j).$$
The 
signature operator of $\whW^*$
is defined to be
$$\Dw=i (\dw\tw + \tw \dw).
\Eq(signw)$$
(The right-hand-side of \equ(signw) differs from the right-hand-side of
\equ(oddsign) by a sign; the reason for this will become apparent
in the formula for ${\cal D}_C^{\rm sign} (\epsilon)$
given below.)
%#New macros here: $\dw=d_{\whW}$, $\tw=\tau_{\whW}$
%\Dw={\rm D}_{\whW}

Let $g:W^*\rightarrow \Omega^*(F;\Fli)$
be the dual to $f$ with respect to the Hermitian forms, i.e.
$$Q_W(f(\alpha),z)=
Q_F( \alpha, g(z)).$$
Then we leave to the reader the proof of the following lemma.

\nproclaim Lemma [antil]. $g$ commutes with the differentials.
If $f^*$ denotes the adjoint
of $f$ with respect to the inner products $<\cdot, \cdot>_F$
and $<\cdot, \cdot>_W$ then
$f^* = \tau_F g \tau_W$.

Using the isomorphism between $W^*$ and $\widehat{W}^*$ in
\equ(newcomplex), 
let $\widehat{f} :  \Omega^*(F;\Fli)\rightarrow \widehat{W}^*$ and
$\widehat{g}: 
\widehat{W}^*\rightarrow \Omega^*(F;\Fli)$ be the obvious extensions
of $f$ and $g$.
Define a cochain complex $C^* = \bigoplus_{k=-1}^{2m} C^k$ by 
$C^k=\Omega^k(F;\Fli)\oplus \whW^k$. Given $\epsilon \in \RR$, 
define a differential 
$D_C$  on $C^*$ by
$$D_C= 
\left( \matrix{ D_F & \epsilon \widehat{g} \cr 0  & - \dw} \right)
\;{\rm if}\,*< m- {1 \over 2},\, \quad\quad 
D_C= 
\left( \matrix{ D_F & 0 \cr - \epsilon \widehat{f} & - \dw} \right)
\;{\rm if}\,*> m-{1 \over 2},
\Eq(diffp)$$
where $D_F$ has been defined in \equ(signed-diff).
Since $D_{\widehat{W}} \widehat{f} = \widehat{f} D_F$ and 
$D_F \widehat{g} = \widehat{g} D_{\widehat{W}}$,
we have
$(D_C)^2=0$.

If $\epsilon > 0$ then the complex $(C^*,D_C)$ has 
vanishing cohomology, as can be seen by Lemma \thm[complex] and the 
mapping-cone nature of the construction of $(C^*,D_C)$.
Define a duality operator $\tau_C$
on $C^*$ by
$$\tau_C = 
\left( \matrix{ \tau_F & 0 \cr 0 & \tw} \right).
\Eq(tauc)$$
There is also a Hermitian form $Q_C: C^k\times C^{(2m-1)-k}
\rightarrow {\cal B}^\infty$ given by
$$Q_C((\alpha,v),(\beta,z))
=
Q_F( \alpha , \beta ) + Q_{\whW}(v,z).$$ Note that $C$ has formal dimension
$2m-1$. 
We obtain a Hermitian inner product on $C^*$
by 
$$<\cdot, \cdot>_C = Q_C(\cdot, \tau_C \cdot).$$
The signature operator of $(C^*,D_C)$ is defined to be
${\cal D}_C^{\rm sign} (\epsilon)
=-i(\tau_C D_C+ D_C \tau_C)$ and is  given
on the degree-$j$ subspace by
$${\cal D}_C^{\rm sign} (\epsilon) =
(-i)\left( \matrix{ D_F\tau_F + \tau_F D_F & 0 \cr
                  0 & - (\dw\tw + \tw\dw) } \right)
+(-i)\cases{\left(\matrix{0 & \epsilon\tau_F \widehat{g} \cr
-\epsilon \widehat{f} \tau_F & 0 } \right)&if$\;j < m- {1 \over 2}$\cr 
\left(\matrix{0 & \epsilon \widehat{g} \tw \cr
-\epsilon \tw \widehat{f} & 0 } \right)&if$\;j > m - {1 \over 2}$\cr} 
\Eq(psign)$$
If $\epsilon >0$, it follows from the vanishing of the cohomology of $C^*$ 
that ${\cal D}_C^{\rm sign}(\epsilon)$ is an invertible self-adjoint 
$\Bi$-operator. 
Namely, $\left( {\cal D}_C^{\rm sign}(\epsilon) \right)^2$
is the Laplace operator $D_C D_C^* + D_C^* D_C$ on $C^*$. From the
method of proof of [L3, Propositions 10 and 27], the
vanishing of the cohomology of $C^*$ implies the
invertibility of $D_C D_C^* + D_C^* D_C$.

We are now in a position to recall the
definition of the higher eta invariant. Suppose that $F$ satisfies
Assumption 1.
Define the space $\overline{\Omega}_*(\Bi)$ of noncommutative differential
forms as in [L1, Section II].
Define a rescaling operator ${\cal R}$ on 
$\overline{\Omega}_*(\Bi)$ which acts on
$\overline{\Omega}_{2j}(\Bi)$
as multiplication by $(2 \pi i)^{-j}$
and acts on
$\overline{\Omega}_{2j-1}(\Bi)$
as multiplication by $(2 \pi i)^{-j}$.

Let $$\nabla^\Omega : \Omega^*(F; {\cal V}^\infty) \rightarrow
\Omega_1(\Bi) \otimes_{\Bi} \Omega^*(F; {\cal V}^\infty)$$
be the connection
constructed in [L1, Proposition 9], 
in terms of a function $h \in C^\infty_0(F^\prime)$
such that $\sum_{\gamma \in \Gamma} \gamma \cdot h = 1$.
(Recall that $F^\prime$ is a normal $\Gamma$-cover of $F$.)
As in [L5, (3.28)], let $$\nabla^W : W^* \rightarrow
\Omega_1(\Bi) \otimes_{\Bi} W^* \Eq(sdconnect)$$ 
be a connection on $W^*$ which is
invariant under $\tau_W$ and preserves $Q_W$. Let $\nabla^{\widehat{W}^*}$ 
be the obvious extension of $\nabla^W$ to $\widehat{W}^*$ and put
$\nabla^C = \nabla^\Omega \oplus \nabla^{\widehat{W}^*}$.

Let $\Cl(1)$ be the complex Clifford algebra of $\CC$
generated by 1 and $\sigma$, with $\sigma^2=1$, and let
$\STR1$ be the supertrace as in [LP1].
Let $\eps \in C^\infty(0, \infty)$ now be a nondecreasing function such that
%#New macro : \eps=\epsilon
$\eps(s)=0$ for $s\in (0,1]$ and $\eps(s)=1$
for $s\in [2,+\infty)$. 
Consider
$$\widetilde{\eta}_F(s)=
{1 \over \sqrt{\pi}} \, {\cal R} \, \STR1 \left( {d\over ds}
[\sigma s{\cal D}_C^{\rm sign}(\eps(s)) + \nabla^C] \right)
\exp[-(\sigma s{\cal D}_C^{\rm sign}(\eps(s)) + 
\nabla^C)^2]\in\overline{\Omega}_{\rm even}(\Bi).
\Eq(integrand)$$
%#New macro: \STR1={\rm STR}_{{\rm Cl}(1)}

The higher eta invariant of $F$ is, by definition,
$$\widetilde{\eta}_F=\int_0^\infty \widetilde{\eta}_F(s) ds \in
\overline{\Omega}_{\rm even}(\Bi)/d \overline{\Omega}_{\rm odd}(\Bi).
\Eq(heta)$$
It is shown in
[L5, Proposition 14] that $\widetilde{\eta}_F$ is independent of
the particular choices of the function 
$\epsilon$, the perturbing complex $W^*$ and the
self-dual connection $\nabla^W$.
Definition \equ(heta) can be seen as a way
of regularizing the {\it a priori} divergent integral
$${1 \over \sqrt{\pi}} \, {\cal R} \, \int_0^\infty \STr1 \left({d\over ds}
[\sigma s\Si + \nabla] \right)
\exp[-(\sigma s\Si + \nabla)^2]\, ds
\Eq(hetainv)$$
coming from the signature operator $\Si$ of $F$.
It is not clear that the integrand in \equ(hetainv) is integrable for
large $s$, as the spectrum of $\Si$ may include zero. To get around this
problem, we have first added the complex $\widehat{W}^*$, whose
higher eta-invariant formally vanishes by a duality argument.  Then we
have perturbed the direct sum differential so that for large $s$, we are 
dealing with the invertible signature operator ${\cal D}_C^{\rm sign}(1)$. 

The invertibility of ${\cal D}_C^{\rm sign}(1)$ ensures that
the integrand in \equ(heta) is integrable for
large $s$; see the proof of [L3, Proposition 28] in the analogous but
more difficult case of the analytic torsion form. From
[L2, Proposition 26], the integrand in \equ(heta) is integrable for
small $s$.

\smallskip
\noindent
{\bf Remark.}
A different regularization
of \equ(hetainv) has been proposed in [LP4]
using the notion of {\it symmetric spectral section}.
See the Appendix for an informal
argument showing the equality of the two regularizations.

\smallskip
The higher eta-invariant satisfies
$$d \widetilde{\eta}_F = \int_F L (R^F/2 \pi) \wedge \omega,
\Eq(funny)$$
where the closed biform $\omega \in \Omega^*(F) \, \widehat{\otimes} \,
\overline{\Omega}_*({\cal B}^\infty)$ is given in [L1, Section V]. 
In fact, $\omega$ is the image
of an element of $\Omega^*(F) \, \widehat{\otimes} \,
\overline{\Omega}_*(\CC \Gamma)$ under the map 
$\overline{\Omega}_*(\CC \Gamma)
\rightarrow \overline{\Omega}_*({\cal B}^\infty)$. 
(This follows from the
fact that the function $h$ used to define
$\nabla^\Omega$ and $\omega$ from [L1, (40)] 
has compact support on $F^\prime$.)
By abuse of notation,
we will also denote this element of $\Omega^*(F) \, \widehat{\otimes} \,
\overline{\Omega}_*(\CC \Gamma)$ by $\omega$. It satisfies
the property that if $Z_\tau$ is a cyclic cocycle which represents a
cohomology class $\tau \in \HH^*(\Gamma; \CC)$ then
$\langle \omega, Z_\tau \rangle \in \Omega^*(F)$ is an explicit closed
form on $F$ whose de Rham cohomology class is a nonzero 
constant (which only depends on the degree of $\tau$) times $\nu^* \tau$.

\medskip
\noindent
{\bf Conventions : } 
Let us take this occasion to 
establish our conventions for Chern characters.  If
$\nabla$ is a connection on a vector bundle then its Chern character
is 
$$\ch(\nabla) \, = \,  \TR \left( e^{-{\nabla^2 \over 2\pi i}} 
\right)
\, = \, {\cal R} \, \TR \left( e^{-\nabla^2} \right).$$
The de Rham cohomology class of $\ch(\nabla)$
is the representative of a {\it rational} cohomology class.
Similarly, the de Rham cohomology class of the $L$-form
$L(R^F/2 \pi)$ lies in the
image of the map $\HH^*(F; \QQ) \rightarrow \HH^*(F; \RR)$.
If $\AAA$ is a superconnection then its Chern character is
$\ch(\AAA) \, = \, {\cal R} \, {\rm STR} \left( e^{-\AAA^2} \right)$.
\medskip

%\fine
\expandafter\ifx\csname sezioniseparate\endcsname\relax%
    \input macro \fi

%--------------------------------------------------------
\numsec=3
\numfor=1
\numtheo=1
\pgn=1

\beginsection 3. The higher eta invariant of an even-dimensional manifold

Before dealing with the case of even-dimensional $F$, we introduce the
notion of a Lagrangian subspace of a ${\cal B}^\infty$-module.
Let ${\HH}$ be a finitely-generated
projective ${\cal B}^\infty$-module with
a nondegenerate quadratic
form $Q_{{\HH}} : {\HH} \times {\HH} \rightarrow
{\cal B}^\infty$ such that
$Q_{{\HH}}(bx, y) = b Q_{{\HH}}(x,y)$ and
$Q_{{\HH}}(x,y)^* = Q_{{\HH}}(y,x)$.
A {\it Lagrangian subspace} of ${\HH}$
is a finitely-generated projective
${\cal B}^\infty$-submodule $L$ on which
$Q_{{\HH}}$ vanishes, such that $L$ equals $L^{\perp}$, its
orthogonal space with respect to $Q_{{\HH}}$. Equivalently, let
$L$ be a finitely-generated projective
${\cal B}^\infty$-submodule of ${\HH}$. Let 
$L^\prime$ be the antidual to $L$, i.e. the set of $\RR$-linear maps
$\l^\prime : L \rightarrow {\cal B}^\infty$ such that
$l^\prime (bl) = l^\prime(l) b^*$ for all $b \in {\cal B}^\infty$ and
$l \in L$. Here $L^\prime$ is also a left ${\cal B}^\infty$-module, with the
multiplication given by $(a l^\prime)(l) = a l^\prime(l)$. 
Then for $L$ to be a Lagrangian subspace of ${\HH}$ amounts to
the existence of a short exact sequence
$$0 \longrightarrow L \longrightarrow {\HH} \longrightarrow
L^\prime \longrightarrow 0 \Eq(seq)$$
whose maps are an injection
$i : L \rightarrow {\HH}$ and its antidual (with
respect to $Q_{{\HH}}$)
$i^\prime : {\HH} \rightarrow L^\prime$.

If ${\cal L}$ is a finitely-generated projective 
${\cal B}^\infty$-module then there is a canonical quadratic form on 
${\cal L} \oplus {\cal L}^\prime$ given by
$$Q(l_1 + l_1^\prime, l_2 + l_2^\prime) = l_1^\prime(l_2) +
(l_2^\prime(l_1))^*.$$ It has a canonical Lagrangian subspace given by 
${\cal L}$. (In what follows, it will in fact
suffice to take ${\cal L}$ of the
form $({\cal B}^\infty)^N$.)
A {\it stable Lagrangian subspace} of ${\HH}$ is
a Lagrangian subspace $L$
of ${\cal H} = {\HH} \oplus  ({\cal L} \oplus {\cal L}^\prime)$ for some
${\cal L}$ as above. 
We say that two stable Lagrangian subspaces of $\HH$, 
$L_1 \subset {\cal H}_1$ 
and
$L_2 \subset {\cal H}_2$,
are equivalent if there are ${\cal L}_3$ and ${\cal L}_4$, and an isomorphism 
$j : {\cal H}_1 \oplus {\cal L}_3 \oplus
{\cal L}_3^\prime \rightarrow {\cal H}_2 \oplus {\cal L}_4 \oplus
{\cal L}_4^\prime$ of quadratic form spaces, such that
$j(L_1 \oplus {\cal L}_3) = L_2 \oplus {\cal L}_4$.

Now suppose that $F$ is a closed oriented manifold of dimension
$n = 2m$. Let $\nu : F \rightarrow
B\Gamma$ be a continuous map as before.
We make the following assumption.
\medskip
\noindent
{\bf Assumption 1.a :} The natural surjection
$\HH^{m}(F;\Fl) \rightarrow \overline{\HH}^{m}(F;\Fl)$ is an
isomorphism.
\medskip

\nproclaim Lemma [Riem2].
If $F$ is equipped with a Riemannian metric then
Assumption 1.a is equivalent to saying that the differential form Laplacian
on $\Omega^{m}(F^\prime)$ has a strictly positive spectrum on the 
orthogonal complement of its kernel.

\noindent
{\it Proof.} We give an outline of the proof.
Let $\Omega^*_{(2)}(F; \Fl)$ denote the completion of 
$\Omega^*(F; \Fl)$ as a $C^*_r(\Gamma)$-Hilbert module.
Assumption 1 is equivalent to saying that the differential
$D_F : \Omega^{m-1}(F; \Fl) \rightarrow \Omega^{m}(F; \Fl)$ has a closed
image.  From arguments as in [L3, Propositions 10 and 27], this is
equivalent to saying that $D_F D_F^*$ has a strictly positive spectrum on
$\IM( D_F : 
\Omega^{m-1}_{(2)}(F; \Fl) \rightarrow \Omega^{m}_{(2)}(F; \Fl))$.
Then by Hodge duality, $D_F^* D_F$ has a strictly positive spectrum on
$\Omega^{m}_{(2)}(F; \Fl)/\Ker( D_F)$.
Again as in the proof of
[L3, Propositions 10 and 27], under Assumption 1.a there is an
orthogonal direct sum decomposition of closed $C^*_r(\Gamma)$-Hilbert modules
$$
\Omega^{m}_{(2)}(F; \Fl) = \Ker(D_F D_F^* + D_F^* D_F) \oplus \IM(D_F) \oplus 
\Omega^{m}_{(2)}(F; \Fl)/\Ker( D_F).$$
Thus $D_F D_F^* + D_F^* D_F$ has a strictly positive spectrum on the
orthogonal complement of its kernel.
Put ${\cal V}^{(2)} = l^2(\Gamma) \times_\Gamma F^\prime$.
Using the injective homomorphism $C^*_r(\Gamma) \rightarrow B(l^2(\Gamma))$ of
$C^*$-algebras and [L3, Proposition 19], the lemma now 
follows as in the rest of
the proof of Lemma \thm[Riem].  
\QED
\medskip

Hereafter we assume that Assumption 1.a is satisfied.
The proof of the next lemma is similar to that of Lemma \thm[complex], but
easier, and will
be omitted.

\nproclaim Lemma [complexx].
There
is a cochain complex $W^* = \bigoplus_{i=0}^{2m} W^i$
of finitely-generated projective $\Bi$-modules
such that \hfil\break
\noindent
1. $W^*$ is a graded regular $n$-dimensional Hermitian
complex.  \hfil\break
\noindent
2. The differentials $D_W : W^{m-1} \rightarrow W^m$ and
$D_W : W^{m} \rightarrow W^{m+1}$ vanish. \hfil\break
\noindent
3. There is a double homotopy equivalence
$$f : \Omega^*(F;\Fli)\rightarrow W^* 
\Eq(pr22)$$
which, as an element of $(\Omega^*(F;\Fli))^* \otimes W^*$, is actually
smooth with respect to $F$. 

For brevity, let us denote
${\HH}^m(F; {\cal V}^\infty)$ by ${\HH}$.
Then ${\HH}$ is a finitely-generated projective
${\cal B}^\infty$-module which is isomorphic to the module $W^m$ of Lemma
\thm[complexx].
The quadratic form $Q_F$ restricts to a nondegenerate quadratic form
$Q_{{\HH}}$.
The grading operator $\tau_F$ induces
a grading operator $\tau_{{\HH}}$ on ${\HH}$.
Let
${\HH}^\pm$ be the $\pm 1$-eigenspace of
$\tau_{{\HH}}$. 
Put $\overline{\HH} = C^*_r(\Gamma) \otimes_{{\cal B}^\infty} \HH$, and
similarly for $\overline{\HH}^\pm$.

\medskip
\nproclaim Lemma [coho].
The index of the signature operator on $F$
equals $\left[ \overline{\HH}^+ \right] -
\left[ \overline{\HH}^- \right] \in K_0(C^*_r(\Gamma))$.

\noindent
{\it Proof.}
Put $\overline{W}^* = C^*_r(\Gamma) \otimes_{{\cal B}^\infty} W^*$.
The index of the signature operator of $F$ equals the index of the
signature operator of the complex $\overline{W}^*$ [KM, Theorem 4.1] 
and is independent of the
choice of grading operator $\tau$ [KM, Proposition 3.6]. Hence we may
work with the complex $\overline{W}^*$. Consider the regular Hermitian complex 
$\widetilde{W}^* = \oplus_{i \ne m} \overline{W}^i$. It is enough to show that
the index of the signature operator 
${\cal D}^{{\rm sign},+}_{\widetilde{W}}$ of 
$\widetilde{W}^*$ vanishes.
To show this, define an operator $\mu$ on $\widetilde{W}^*$ by
$$\mu(w) = 
\cases{
w & if $w \in \widetilde{W}^i, i < m$, \cr
- w & if $w \in \widetilde{W}^i, i > m$. \cr}$$
Then $\mu^2 = 1$, $\mu {\cal D}^{\rm sign}_{\widetilde{W}} = 
{\cal D}^{\rm sign}_{\widetilde{W}} \mu$ and
$\mu \tau_{\widetilde{W}} + \tau_{\widetilde{W}} \mu = 0$. 
Let  $\widetilde{W}^*_\pm$ be the $\pm 1$ 
eigenspaces of $\tau_{\widetilde{W}}$.
Then $\mu$ induces an isomorphism from $\widetilde{W}^*_+$ to
$\widetilde{W}^*_-$ and so
$\Ind({\cal D}^{{\rm sign},+}_{\widetilde{W}}) = 
[\widetilde{W}^*_+] - [\widetilde{W}^*_-] = 0$.
\QED 
\medskip

We make the following
further assumption.

\medskip
\noindent
{\bf Assumption 1.b :} ${\HH}^m(F; {\cal V}^\infty)$ admits a
(stable) Lagrangian subspace.
\medskip

\nproclaim Lemma [index].
Given Assumption 1.a, 
Assumption 1.b is equivalent to saying that the index of the
signature operator on $F$ vanishes in $K_0(C^*_r(\Gamma))$.

\noindent
{\it Proof.}
In general, if ${\cal L}$ is a finitely-generated projective 
${\cal B}^\infty$-module, put $\overline{\cal L} = C^*_r(\Gamma) 
\otimes_{{\cal B}^\infty} {\cal L}$. 
If $h$ is a ${\cal B}^\infty$-valued Hermitian
metric on ${\cal L}$, let $\overline{h}$ be its extension to
a $C^*_r(\Gamma)$-valued Hermitian metric on $\overline{\cal L}$. 
Define ${\cal I} : \overline{\cal L} \rightarrow \overline{\cal L}^\prime$ by
$$
\left({\cal I}(l_1) \right)
(l_2) = \overline{h}(l_1, l_2). \Eq(gr)
$$ 
Put
$\tau_{\overline{\cal L} \oplus \overline{\cal L}^\prime} = 
\left(
\matrix{
0 & {\cal I}^{-1} \cr
{\cal I} & 0
\cr}
\right)$.
If $\left( \overline{\cal L} \oplus \overline{\cal L}^\prime \right)^\pm$ 
denotes the $\pm 1$ eigenspaces of 
$\tau_{\overline{\cal L} \oplus \overline{\cal L}^\prime}$ then
$( \overline{\cal L} \oplus \overline{\cal L}^\prime )^+$ is isomorphic to
$( \overline{\cal L} \oplus \overline{\cal L}^\prime )^-$, under
$x + {\cal I}(x) \rightarrow x - {\cal I}(x)$.

Suppose that Assumption 1.b is satisfied. Then there is some 
finitely-generated projective ${\cal B}^\infty$-module
${\cal L}$ such that
$\HH \oplus {\cal L} \oplus {\cal L}^\prime$ has a Lagrangian subspace $L$.
Give $\overline{\HH} \oplus \overline{\cal L} \oplus 
\overline{\cal L}^\prime$ the grading operator $\tau_{\overline{\HH}}
\oplus \tau_{\overline{\cal L} \oplus \overline{\cal L}^\prime}$. Then
from what has been said,
$$[\overline{\HH}^+] - [\overline{\HH}^-] =
[(\overline{\HH} \oplus \overline{\cal L} \oplus 
\overline{\cal L}^\prime)^+] -
[(\overline{\HH} \oplus \overline{\cal L} \oplus 
\overline{\cal L}^\prime)^-] \Eq(K0) 
$$
in $K_0(C^*_r(\Gamma))$. However, as the element of $K_0(C^*_r(\Gamma))$
coming from $\overline{\HH} \oplus \overline{\cal L} \oplus 
\overline{\cal L}^\prime$ is independent of the choice of the grading
operator $\tau$, we can use $\overline{L}$ to define a grading operator on
$\overline{\HH} \oplus \overline{\cal L} \oplus 
\overline{\cal L}^\prime$, as 
in \equ(gr), to see that 
$[(\overline{\HH} \oplus \overline{\cal L} \oplus 
\overline{\cal L}^\prime)^+] -
[(\overline{\HH} \oplus \overline{\cal L} \oplus 
\overline{\cal L}^\prime)^-] = 0
$ in $K_0(C^*_r(\Gamma))$. From \equ(K0) and
Lemma \thm[coho], this implies that the index of the
signature operator on $F$ vanishes in $K_0(C^*_r(\Gamma))$.

Now suppose that the index of the
signature operator on $F$ vanishes in $K_0(C^*_r(\Gamma))$. 
Then there is some $N \ge 0$ such
that  $\overline{\HH}^+ \oplus C^*_r(\Gamma)^N$ is
isomorphic to $\overline{\HH}^- \oplus C^*_r(\Gamma)^N$.
We can take the isomorphism 
$\overline{j} :  \overline{\HH}^+ \oplus C^*_r(\Gamma)^N \rightarrow
\overline{\HH}^- \oplus C^*_r(\Gamma)^N$ to be an isometry.
Using arguments as in [LP1, Appendix A], we can assume that
$\overline{j} = j \otimes_{{\cal B}^\infty} I$ for some isometric isomorphism
${j} :  {\HH}^+ \oplus ({\cal B}^\infty)^N \rightarrow
{\HH}^- \oplus ({\cal B}^\infty)^N$. Then 
graph(${j}$) is a Lagrangian
subspace of ${\HH} \oplus ({\cal B}^\infty)^N \oplus ({\cal B}^\infty)^N$.
Thus Assumption 1.b is satisfied.
\QED
\medskip

\nproclaim Corollary [existence].
Given Assumption 1.a,
if $F$ is the boundary of a compact oriented manifold $M$ and
$\nu$ extends over $M$  
then Assumption 1.b is satisfied.

\noindent
{\it Proof.}  This follows from the cobordism invariance of the index,
along with Lemma \thm[index].
\QED
\medskip

\nproclaim Lemma [examples2]. \hfil\break
(a) If $F$ has a cellular decomposition without any cells of
dimension $m$ then Assumptions  1.a and 1.b are satisfied. \hfil\break
\noindent
(b) If $\Gamma$ is finite and the signature of $F$ vanishes then
Assumptions
1.a and 1.b are satisfied. \hfil\break
\noindent
(c) Let $F_1$ and $F_2$ be even-dimensional manifolds, with 
$F_1$ a connected closed hyperbolic manifold and $F_2$ a closed manifold with
vanishing signature. Put $\Gamma = \pi_1(F_1)$. If $F = F_1 \times F_2$
then Assumptions 1.a and 1.b are satisfied. 

\noindent
{\it Proof.} \hfil\break
(a) If $F$ has a cellular decomposition without any cells of
dimension $m$ then $\HH^{m}(F;\Fl)$ vanishes and Assumptions 1.a and 1.b are
automatically satisfied \hfil\break
(b) If $\Gamma$ is finite then $F^\prime$ is compact and from standard
Hodge theory, the result
of Lemma \thm[Riem2] is satisfied. From Lemma \thm[index],
it remains to show that $F$ has vanishing index in $K_0(C^*_r(\Gamma))$.
Now $K_0(C^*_r(\Gamma))$ is isomorphic to the ring of complex virtual
representations of $\Gamma$. Given a representation 
$\rho : \Gamma \rightarrow U(N)$,
the corresponding component of the $K_0(C^*_r(\Gamma))$-index is the 
usual index of the signature operator acting on
$\Omega^*(F^\prime) \otimes_\rho \CC^N$. By the Atiyah-Singer index theorem
this equals $N$ times the signature of $F$, and hence vanishes. 
\hfil\break
(c) From [Do], the spectrum of the differential form Laplacian on the
hyperbolic space $\widetilde{F_1}$ is strictly positive on the orthogonal
complement of its kernel (which is concentrated in the middle degree).
Then by separation of variables and using the fact that the
universal cover $\widetilde{F_2}$ is compact, it follows that 
the result of Lemma \thm[Riem2] is satisfied.  From Lemma
\thm[index] and the multiplicativity of the index, along with
the vanishing of the signature of $F_2$, we obtain that Assumption 1.b is
satisfied.
\QED
\medskip

Hereafter we assume that $\HH$ admits a stable Lagrangian subspace.
Let $L \subset
{\HH} \oplus {\cal L} \oplus {\cal L}^\prime$ be one such.
We define a new complex
$\whW^*$ 
by 
$$\whW^i=
\cases{
W^{i+1}& if $-1\leq i \leq m-2$, \cr
L & if $i = m-1$ \cr
0 & if $i = m$ \cr
L^\prime & if $i = m+1$ \cr
W^{i-1}& if $m+2 \leq i \leq 2m+1$. \cr} \Eq(newcomplexx)$$
There is an obvious extension of $D_W$ to a differential
$D_{\whW}$ and obvious extensions of $Q_W$ and $\tau_W$ to 
$\whW$, at least on the part of $\whW$ that does not involve $L$ or $L^\prime$.
Define
$Q_{\whW} : L \times L^\prime \rightarrow {\cal B}^\infty$ by
$Q_{\whW}(l, l^\prime) = (l^\prime(l))^*$. 
Let $h$ be the ${\cal B}^\infty$-valued Hermitian metric on $L$
induced from ${\HH} \oplus {\cal L} \oplus {\cal L}^\prime$. Define
$\tau_{\whW} : L \rightarrow L^\prime$ by
$\left( \tau_{\whW} (l_1) \right)(l_2) = h(l_1, l_2)$. Let
$\tau_{\whW} : L^\prime \rightarrow L$ be the inverse.
Then we obtain a well-defined
triple $(D_{\whW}, Q_{\whW}, \tau_{\whW})$ on $\whW$. 

Let $\Omega^*(F; \Fli) \oplus {\cal L} \oplus {\cal L}^\prime$ be the
direct sum cochain complex, with ${\cal L} \oplus {\cal L}^\prime$ 
concentrated in
degree $m$. Recall the notation $i$ and $i^\prime$ for the maps in \equ(seq),
where $\HH$ is again $\HH^m(F; \Fli) = \overline{\HH}^m(F; \Fli)$. 
Let ${\cal I} : \HH^m(F; \Fli) \rightarrow \Omega^m(F;\Fli)$ be the
inclusion coming from Hodge theory and let  
${\cal I}^* : \Omega^m(F;\Fli) \rightarrow \HH^m(F; \Fli)$ be orthogonal
projection.
Define $\widehat{f} : \Omega^*(F;\Fli) \oplus {\cal L} \oplus
{\cal L}^\prime \rightarrow \widehat{W}^*$ to be the obvious extension of
$f$ outside of degree $m$, and to be given in degree $m$ by
$\widehat{f}(\omega + l +l^\prime) = (i^\prime \circ {\cal I}^*)(\omega)
\in \widehat{W}^{m+1}$. 
Define
$\widehat{g}: 
\widehat{W}^*\rightarrow \Omega^*(F;\Fli) \oplus {\cal L} \oplus
{\cal L}^\prime$ to be the obvious extension
of $g$ outside of degrees $m-1$, $m$ and $m+1$, 
to be given in degree $m-1$ by
$\widehat{g}(l) = ({\cal I} \circ i)(l) \in \Omega^m(F; {\cal V}^\infty)$ 
and to vanish in degrees $m$ and
$m+1$. 
Define a cochain complex $C = \bigoplus_{k=-1}^{2m+1} C^k$ by 
$C^*=\Omega^*(F;\Fli) \oplus {\cal L} \oplus {\cal L}^\prime
\oplus \whW^*$. Given $\epsilon \in \RR$, 
define a differential 
$D_C$  on $C$ by
$$D_C= 
\left( \matrix{ D_F & \epsilon \widehat{g} \cr 0  & - \dw} \right)
\;{\rm if}\,* \leq m- 1,\, \quad\quad 
D_C= 
\left( \matrix{ D_F & 0 \cr - \epsilon \widehat{f} & - \dw} \right)
\;{\rm if}\,* \geq m.
\Eq(diffpp)$$
We can then define $\tau_C$, $Q_C$ and ${\cal D}_C^{\rm sign}(\epsilon)$ in
analogy to what we did in the odd-dimensional case.

We put
$$\widetilde{\eta}_F(s)=  \, 
{\cal R} \, {\rm STR} \left({d\over ds}
[s{\cal D}_C^{\rm sign}(\eps(s)) + \nabla^C] \right)
\exp[-(s{\cal D}_C^{\rm sign}(\eps(s)) + 
\nabla^C)^2]\in\overline{\Omega}_{\rm odd}(\Bi) \Eq(oddeta)
$$
where ${\rm STR}$ is the supertrace and $\nabla^C$ is a self-dual connection
as before.
The function $\epsilon(s)$ is the same as in the odd-dimensional case.

The higher eta invariant of $F$ is, by definition,
$$\widetilde{\eta}_F=\int_0^\infty \widetilde{\eta}_F(s) ds \in
\overline{\Omega}_{\rm odd}(\Bi)/d \overline{\Omega}_{\rm even}(\Bi).
\Eq(hetaa)$$
As in [L5, Proposition 14], $\widetilde{\eta}_F$ is independent of
the particular choices of $\epsilon$, the perturbing complex $W^*$ and
the self-dual connection $\nabla^W$.
It satisfies \equ(funny).

Let us consider how $\widetilde{\eta}_F$ depends on the
choice of (stable) Lagrangian subspace $L$. For the moment, let us
denote the dependence by $\widetilde{\eta}_F(L)$. From equation
\equ(funny), if $L_1$ and $L_2$ are two (stable) Lagrangian subspaces then
$d (\widetilde{\eta}_F(L_1) -  \widetilde{\eta}_F(L_2)) = 0$.
Thus $\widetilde{\eta}_F(L_1) -  \widetilde{\eta}_F(L_2)$ represents an
element of $\overline{\HH}_{\rm odd}({\cal B}^\infty)$. To describe it,
we construct a characteristic class coming from two (stable) Lagrangian
subspaces.

Let $\HH$ be a finitely-generated
projective ${\cal B}^\infty$-module as above, equipped with a quadratic form
$Q_{\HH}$. For simplicity, we will only deal with 
honest Lagrangian
subspaces of $\HH$; the case of stable Lagrangian subspaces can be dealt
with by replacing $\HH$ by
${\HH} \oplus ({\cal B}^\infty)^N \oplus ({\cal B}^\infty)^N$.
 
As in the proof of Lemma \thm[index], after choosing a
grading $\tau_H$, the set of Lagrangian subspaces of $\HH$ can be identified
with $\Isom_{{\cal B}^\infty}(\HH^+, \HH^-)$, the set 
of isometric isomorphisms from $\HH^+$ to $\HH^-$.
If $j_1, j_2 \in \Isom_{{\cal B}^\infty}(\HH^+, \HH^-)$ then
$j_1 \circ j_2^{-1} \in \Isom_{{\cal B}^\infty}(\HH^-, \HH^-)$.
Now $\Isom_{{\cal B}^\infty}(\HH^-, \HH^-)$ is homotopy-equivalent to
$\GL_{{\cal B}^\infty}(\HH^-)$. Hence given two Lagrangian subspaces
$L_1$ and $L_2$ of $\HH$, we obtain an element of
$\pi_0 \left( \GL_{{\cal B}^\infty}(\HH^-) \right)$
represented by $j_1 \circ j_2^{-1}$. Let $[L_1 - L_2]$ denote
its image in $K_1({\cal B}^\infty) \cong \pi_0
\left(\GL_{{\cal B}^\infty}(\infty) \right)$. 

\nproclaim Proposition [etadiff].
$$\widetilde{\eta}_F(L_1) - \widetilde{\eta}_F(L_2) =
\, \ch([L_1 - L_2]) \quad \quad {\rm in } \, \, 
\overline{\HH}_{\rm odd}({\cal B}^\infty).$$

\noindent
{\it Proof.} Fix, for the moment,
a Lagrangian subspace $L$ of $\HH$. Writing $L = {\rm graph}(j)$
with $j \in \Isom_{{\cal B}^\infty}(\HH^+, \HH^-)$,
we can identify 
$L$, and hence $L^\prime$, with $\HH^+$. Under these identifications, 
the short exact sequence \equ(seq) becomes
$$0 \longrightarrow \HH^+ \longrightarrow \HH^+ \oplus \HH^- 
\longrightarrow \HH^+ \longrightarrow 0.$$ 
To describe the maps involved explicitly, let us consider this to be a 
graded regular $2m$-dimensional Hermitian complex ${\cal E}^*$ concentrated
in degrees $m-1$, $m$ and $m+1$. Then the maps are given by saying that
if $h_+ \in {\cal E}^{m-1}$ then
$D_{\cal E}(h_+) = \, {1 \over \sqrt{2}} \, (h_+, j(h_+))$, while if
$(h_+, h_-) \in {\cal E}^m$ then
$D_{\cal E}(h_+, h_-) =  \, {1 \over \sqrt{2}} \,(- h_+ + j^{-1}(h_-))$.
If $(h_+, h_-), (k_+, k_-) \in 
{\cal E}^m$ then 
$Q_{\cal E}((h_+, h_-), (k_+, k_-)) =
\langle h_+, k_+ \rangle - \langle h_-, k_- \rangle$, while if
$h_+ \in {\cal E}^{m-1}$ and $k_+ \in {\cal E}^{m+1}$ then
$Q_{\cal E}(h_+, k_+) = \langle h_+, k_+ \rangle$.
If $(h_+, h_-) \in {\cal E}^m$ then
$\tau_{\cal E}(h_+, h_-) = (h_+, - h_-)$, while
$\tau_{\cal E} : {\cal E}^{m \pm 1} \rightarrow {\cal E}^{m \mp 1}$ is the
identity map on $\HH^+$. 

The connection $\nabla^{\HH}$, induced from $\nabla^{\Omega}$,
breaks up as a direct sum $\nabla^{\HH^+} \oplus \nabla^{\HH^-}$. We choose
to put
the connection $\nabla^{\HH^+}$ on both 
${\cal E}^{m-1}$ and ${\cal E}^{m+1}$. We obtain a self-dual connection
$\nabla^{\cal E}$ on ${\cal E}$.

It is convenient to perform a change of basis
by means of the isomorphism
${\cal K} : {\cal E}^{m-1} \oplus {\cal E}^{m+1} \rightarrow
\HH^+ \oplus \HH^-$ given by
${\cal K} (o_1, o_2) = \, {1 \over \sqrt{2}} \, (o_1 - o_2, j(o_1 + o_2)).$
One can compute that
the signature operator ${\cal D}_{\cal E} = D_{\cal E} - 
\tau_{\cal E} D_{\cal E} \tau_{\cal E}$ acts as ${\cal K}$ on 
${\cal E}^{m-1} \oplus {\cal E}^{m+1}$, and as 
${\cal K}^{-1}$ on 
${\cal E}^m = \HH^+ \oplus \HH^-$. Thus using the isomorphism ${\cal K}$ to
identify ${\cal E}^{m-1} \oplus {\cal E}^{m+1}$ with $\HH^+ \oplus \HH^-$,
the signature operator ${\cal D}_{\cal E}$ acts on the total space
${\cal E}^m \oplus ({\cal E}^{m-1} \oplus {\cal E}^{m+1}) \cong
(\HH^+ \oplus \HH^-) \oplus (\HH^+ \oplus \HH^-)$ as
$\left(
\matrix{0 & I \oplus I \cr
I \oplus I & 0 } 
\right)$.
We note that this is indeed an odd operator with respect to the
$\ZZ_2$-grading, as the induced duality operator on the second
$\HH^+ \oplus \HH^-$ factor is
${\cal K} \tau_{\cal E} {\cal K}^{-1} = (-I, I)$. 
One can compute that in the new basis, the connection $\nabla^{\cal E}$
becomes
${\cal K} \nabla^{\cal E} {\cal K}^{-1} =
(\nabla^{\HH^+} \oplus \nabla^{\HH^-}) \oplus 
(\nabla^{\HH^+} \oplus j \nabla^{\HH^+} j^{-1})$.

We now consider the complex
$C^*=\Omega^*(F;\Fli) \oplus {\cal L} \oplus {\cal L}^\prime
\oplus \whW^*$ used to define $\widetilde{\eta}_F$. 
Then ${\cal E}^*$ is a subcomplex of $C^*$ and there is a direct sum
decomposition
$C^* = {\cal E}^* \oplus ({\cal E}^*)^\perp$. 
As a
$\ZZ_2$-graded vector space, we have shown that $C^*$ is isomorphic to
$(\HH^+ \oplus \HH^-) \oplus (\HH^+ \oplus \HH^-)
\oplus ({\cal E}^*)^\perp$, regardless of $L$. For $s > 0$,
put $\AAA_s = s{\cal D}_C^{\rm sign}(\eps(s)) + \nabla^C$, thought of as a
superconnection on this $\ZZ_2$-graded vector space.
With our identifications, the $0$-th order part of $\AAA_s$, namely
$s{\cal D}_C^{\rm sign}(\eps(s))$ is independent of $L$. Furthermore,
for large
$s$, the operator ${\cal D}_C^{\rm sign}(\eps(s))$ is invertible.
However, the connection part of $\AAA_s$,
$$\nabla^C = (\nabla^{\HH^+} \oplus \nabla^{\HH^-}) \oplus 
(\nabla^{\HH^+} \oplus j \nabla^{\HH^+} j^{-1}) \oplus 
\nabla^{{\cal E}^\perp}, \Eq(connpart)$$ 
does depend on $L$ through the map $j : \HH^+ \rightarrow \HH^-$.

We are reduced to studying how $\widetilde{\eta}_F$ depends on the
connection part of the superconnection.  
In general, if $\{\AAA(u)\}_{u \in
[0,1]}$ is a smooth $1$-parameter family of superconnections of the form
$$\AAA(u) = \sum_{j=0}^\infty \AAA_{[j]}(u)$$ and we put
$$\AAA_s(u) = \sum_{j=0}^\infty s^{1-j} \AAA_{[j]}(u)$$
then from [L2, (49)], modulo exact forms,
$${d \over du} \widetilde{\eta}(s) =  \,
{d \over ds} \, {\cal R} \, {\rm STR} \, {dA_s \over du} \, e^{- A_s^2}.
$$
Hence, when it can be justified,
$${d \over du} \widetilde{\eta} = 
\, {\cal R} \, {\rm STR} \, {dA_s \over du} \, 
e^{- A_s^2} \Big|_{s = \infty} - 
\, {\cal R} \, {\rm STR} \, {dA_s \over du} \,
e^{- A_s^2} \Big|_{s = 0}.$$

Now consider two Lagrangian subspaces $L_1$ and $L_2$ of $\HH$. 
Choose a $1$-parameter family 
$\{
\nabla^{\widehat{W}^*(u)} 
\}_{u \in [0,1]}$ of 
self-dual connections on 
$\widehat{W}^*$ such that $\nabla^{\widehat{W}^*}(1)$ 
is the connection coming from
$L_1$ and $\nabla^{\widehat{W}^*}(0)$ is the connection coming from
$L_2$.
In our case, the invertibility of ${\cal D}_C^{\rm sign}(\eps(\infty))$ implies
that ${\cal R} \, {\rm STR} \, {dA_s \over du} \, 
e^{- A_s^2} \Big|_{s = \infty} = 0$. For small $s$, the complexes
$\Omega^*(F; {\cal V}^\infty)$ and $\widehat{W}^*$ decouple.
After making a change of basis as above, 
the only $u$-dependence of $\AAA_s(u)$ arises from the
$u$-dependence of $\nabla^{\widehat{W}^*}$. Hence
$${\cal R} \, {\rm STR} \, {dA_s \over du} \,
e^{- A_s^2} \Big|_{s = 0} =
\, {\cal R} \, {\rm STR} \, {d\nabla^{\widehat{W}^*} \over du} \, 
e^{- (\nabla^{\widehat{W}^*})^2}$$
and so
$$\widetilde{\eta}_F(L_1) - \widetilde{\eta}_F(L_2) =
-  \int_0^1 \, {\cal R} \, {\rm STR} \, 
{d\nabla^{\widehat{W}^*} \over du} \, 
e^{- (\nabla^{\widehat{W}^*})^2} du.
$$
Let $j_1, j_2 \in \Isom_{{\cal B}^\infty}(\HH^+, \HH^-)$ be the maps
corresponding to $L_1$ and $L_2$.
Recall that $[L_1 - L_2]$ denote
the element of $K_1({\cal B}^\infty) \cong \pi_0
\left(\GL_{{\cal B}^\infty}(\infty) \right)$ corresponding to
$j_1 \circ j_2^{-1}$.
As $- \int_0^1 \, {\cal R} \, {\rm STR} \, 
{d\nabla^{\widehat{W}^*} \over du} \, 
e^{- (\nabla^{\widehat{W}^*})^2} du$ is the Chern character of
$j_1 \circ j_2^{-1}$ [Ge, Definition 1.1], the proposition follows. \QED

\medskip
\noindent
 {\bf Remark :} In Assumption 1.b of the introduction, 
instead of assuming that
the index  $\Ind_F$ of the signature operator
vanishes in $K_0(C^*_r(\Gamma))$, for our purposes it suffices
to assume that it vanishes in $K_0(C^*_r(\Gamma)) \otimes_{\ZZ} \QQ$. If
this is the case then there is an integer $N > 0$ such that $N \Ind_F$
vanishes in $K_0(C^*_r(\Gamma))$. We can then take $N$ disjoint copies of
$F$, choose a (stable) Lagrangian subspace of $\CC^N \otimes \HH$,
go through the previous construction of $\widetilde{\eta}_F$ and divide
by $N$.  

%\fine
\expandafter\ifx\csname sezioniseparate\endcsname\relax%
   \input macro \fi

%--------------------------------------------------------
\numsec=4
\numfor=1
\numtheo=1
\pgn=1

\beginsection 4. Manifolds with boundary: the perturbed
signature operator

Let $M$ be a compact oriented manifold-with-boundary of dimension $2m$. 
We fix
a non-negative boundary defining function
$x\in C^\infty(M)$ for $\pa M$ and 
a Riemannian metric on $M$ which is isometrically a product
in an (open) collar neighbourhood ${\cal U}\equiv (0,2)_x \times\pa M$
of  the boundary. We let $\Omega^*_c(M)$ denote the compactly-supported
differential forms on the interior of $M$.

Let $\Gamma$ be a finitely-generated discrete group. Consider a continuous map
$M \rightarrow B\Gamma$, with corresponding normal $\Gamma$-cover
$M^\prime \rightarrow M$. Let ${\cal B}^\infty$ be a subalgebra of
$C^*_r\Gamma$ as in Section 1.

We consider the following bundles of left modules over $M$:
$$
\Fl=C^*_r(\Gamma)\times_\Gamma \tM %#New macro \t* for \widetilde{*}
\quad\quad
\Fli=\Bi\times_\Gamma \tM,
$$
and denote their restrictions to the boundary $\partial M$ by
$\Fl_0$ and 
$\Fli_0$.
We suppose that Assumption 1 of the introduction is satisfied, with
$F = \partial M$.
Under this assumption, we shall define,
following [L5, Appendix A],
a perturbation of the differential complex
on $M$.
We shall also give the product structure, near the boundary,
of the
associated signature operator.
\medskip

Using Assumption 1 and following the construction 
of the previous
section with $F=\pa M$,
we  obtain a  perturbed differential complex
on the boundary $\pa M$; this complex
is constructed in terms
of a graded regular Hermitian complex
$W^*$ which is homotopy equivalent to $\Omega^*(\pa M;\Fli_0)$.

\smallskip
\noindent
{\bf Notation.}  
We shall denote the complex on the boundary by
$C^*_0$; thus $C^*_0=\Omega^*(\pa M;\Fli_0)\oplus \widehat{W}^*$.
In general, we let the subscript $0$ 
denote something living
on $\partial M$.

\smallskip

Equation \equ(psign), with $F=\pa M$, defines
 an {\it invertible}
boundary-signature operator 
$${\cal D}_C^{\rm sign} (1)_0:\Omega^*(\pa M;\Fli_0)\oplus \widehat{W}^*
\rightarrow \Omega^*(\pa M;\Fli_0)\oplus \widehat{W}^*.
$$
We wish to  realize ${\cal D}_C^{\rm sign} 
(1)_0$ as the boundary component of 
the signature operator ${\cal D}_C^{\rm sign}(1)$ associated to a {\it
perturbed complex} $(C^*,D_C)$ on $M$. 

To this end, consider the Hermitian $\Bi$-cochain complex 
$\Omega^*_c(0,2) \, \widehat{\otimes} \, \widehat{W}$. We imitate the results of
\equ(Top), thinking of $\widehat{W}$ as algebraically similar to
$\Omega^*(\partial{M}; \Fli_0)$. Thus we have objects
$Q_{\rm alg}$, $\tau_{\rm alg}$ and $D_{\rm alg}$ defined on 
$\Omega^*_c(0,2) \, \widehat{\otimes} \, \widehat{W}$ by
the formulas in \equ(Top), replacing the ``$\partial M$'' on the
right-hand-side of \equ(Top) by ``$\widehat{W}$'' and changing the sign
of $D_{\widehat{W}}$. Recalling that $\widehat{\alpha}= i^{|\alpha|}\alpha$
we thus have
$$
Q_{\rm alg} (dx \wedge \alpha, 1 \wedge \beta) = \int_0^2
Q_{\widehat{W}}(\alpha(x), \widehat{\beta(x)}) dx, 
$$
$$
Q_{\rm alg} (1 \wedge \alpha, dx \wedge \beta) = \int_0^2
Q_{\widehat{W}}(\widehat{\alpha(x)}, \beta(x)) dx, 
$$
$$
\tau_{\rm alg}(1 \wedge \alpha) = dx \wedge \tau_{\widehat{W}} 
\widehat{\alpha},
\Eq(Top2)$$
$$
\tau_{\rm alg}(dx \wedge \alpha) = 1 \wedge i^{-(2m-1)} 
\tau_{\widehat{W}} \widehat{\alpha},
$$
$$
D_{\rm alg}(1 \wedge \alpha) = (1 \wedge - D_{\widehat{W}} \alpha) + 
(dx \wedge \partial_x \widehat{\alpha}),
$$
$$
D_{\rm alg}(dx \wedge \alpha) = dx \wedge i D_{\widehat{W}} \alpha.
$$
One easily checks that the dual to $D_{\rm alg}$, with respect
to $Q_{\rm alg}$, is $D_{\rm alg}^\prime = - D_{\rm alg}$.

Define a new $\Bi$-cochain
complex $C^*$ by
$$C^*=\Omega^*_c(M;\Fli)\oplus (\Omega^*_c(0,2) \, \widehat{\otimes} \,
      \widehat{W}^*).
$$
It inherits objects $Q_C$, $\tau_C$ and $D_C$ from the direct sum
decomposition. Consider the open collar ${\cal U}$ of 
$\partial M$, with ${\cal U}\cong
(0,2)\times\pa M$. 
The bundle $\Fli|_{\cal U}$ is isomorphic to
$(0,2)\times\Fli_0$. Using this isomorphism, 
we can identify the elements of 
$\Omega^*_c(M;\Fli)$ with support in ${\cal U}$, with
$\Omega^*_c(0,2) \, \widehat{\otimes} \, \Omega^*(\pa M;\Fli_0)$.

Now we construct a perturbation of the differential $D_C$ to 
an ``almost'' differential on the complex 
$C^*$. Let $\phi \in C^\infty(0,2)$ be a nonincreasing function satisfying
$\phi(x) = 1$ for $0 < x \leq {1 \over 4}$ and  
$\phi(x) = 0$ for ${1  \over 2} \leq x < {2 }$.
Using Assumption 1, we construct
a homotopy equivalence  $f : \Omega^*(\partial M; \Fli_0) \rightarrow W^*$ 
with adjoint
$g : W^* \rightarrow \Omega^*(\partial M; \Fli_0)$,
exactly as in \equ(pr2) of Section 2 (but with $F=\pa M$).
We define $\widehat{f} : \Omega^*(\partial M; \Fli_0) \rightarrow 
\widehat{W}^*$ and
$\widehat{g} : \widehat{W}^* \rightarrow \Omega^*(\partial M; \Fli_0)$ as
in Section 2. We extend $\widehat{f}$ and $\widehat{g}$ to act on
$\Omega^*_c(0,2) \, \widehat{\otimes} \, \Omega^*(\partial M; \Fli_0)$ and
$\Omega^*_c(0,2) \, \widehat{\otimes} \, \widehat{W}^*$, respectively, by
$$\widehat{f} (\omega_0 + dx \wedge \omega_1) \, = \,
\widehat{f} (\omega_0) - \, i \,  dx \wedge \widehat{f}(\omega_1)$$
and
$$\widehat{g} (w_0 + dx \wedge w_1) \, = \,
\widehat{g} (w_0) - \, i \,  dx \wedge \widehat{g}(w_1).$$
Using the cutoff function $\phi$ and the product structure
on ${\cal U}$, it makes sense to define an operator on $C^*$ by
$$D_C = 
\cases{
\left(
\matrix{D_M & \phi \widehat{g}  \cr
0 & D_{\rm alg} } 
\right) &if$\;* \le m - 1$, \cr
& \cr
\left(
\matrix{D_M & 0 \cr
0 & D_{\rm alg} } 
\right) &if$\;* = m$, \cr
& \cr
\left(
\matrix{D_M & 0 \cr
- \phi \widehat{f} & D_{\rm alg} } 
\right) & if$\;* \ge m + 1$
.\cr
}
\Eq(diffdef)$$
Note that $(D_C)^2 \ne 0$, as $\phi$ is nonconstant.
With our conventions, we have
$$D_C +(D_C)^* = D_C -\tau_C D_C \tau_C\,.
$$
The next lemma follows from the same calculations as at the end
of  Section 1.

\nproclaim Lemma [ident.2].
Define an operator $\Theta$ on 
$\Omega^*_c(0,1/4) \, \widehat{\otimes} \, C^*_0$ by
$$
\Theta((1 \wedge \alpha) + (dx \wedge \beta)) =
(1 \wedge - i^{-|\beta|} \beta) + (dx \wedge i^{|\alpha|} \alpha).
$$
Then when restricted to ${\cal U}$, we can write $D_C +(D_C)^*$ in the form
$D_C +(D_C)^* = \Theta (\partial_x + H)$.
Define an isomorphism $\Phi$ from $C^\infty_c(0,1/4)
\otimes C^*_0$ to the 
+1-eigenspace
$E^+ \subset \Omega^*_c(0,1/4) \otimes C^*_0$ of $\tau_C$ by
$$\Phi(\alpha) = (dx \wedge \alpha) + \tau_C(dx \wedge \alpha).$$
Then
$$
\Phi^{-1} H \big|_{E^+} \Phi = {\cal D}_C^{\rm sign} (1)_0
$$
and 
$$
\left( D_C + (D_C)^* \right)^+ = \Theta \cdot \Phi (\pa_x + 
{\cal D}_C^{\rm sign} (1)_0)
\Phi^{-1}\;.
$$

\smallskip
\noindent
{\bf Notation.} We shall denote the signature operator
$D_C + (D_C)^*$ by ${\cal D}_C^{\rm sign}$.
%# New macro: \SC={\cal D}_C^{\rm sign}
The content of the  above lemma is that the boundary operator corresponding to
$\SC$ 
is precisely the odd perturbed signature 
operator \equ(psign) introduced in Section 2 for closed manifolds,
with $\epsilon = 1$. 

\smallskip

%\fine
\expandafter\ifx\csname sezioniseparate\endcsname\relax%
   \input macro \fi

%--------------------------------------------------------
\numsec=5
\numfor=1
\numtheo=1
\pgn=1

\beginsection 5. The conic index class

We wish to apply the formalism of Hilsum-Skandalis [HS] to show that
the higher signatures of manifolds-with-boundary are homotopy-invariant.
The approach of [HS] is to show that the index of an appropriate
Fredholm operator is homotopy-invariant.  In particular, to apply the
results of [HS, Sections 1 and 2], 
we need to have an operator with a $C^*_r(\Gamma)$-compact
resolvent. For this reason, we will replace the product metric on
$(0,2) \times \partial M$ with a conic metric. 
Recall that $M$ has dimension $2m$.

We keep the notation of Section 4
and assume that $\pa M$ satisfies {\bf Assumption 1}. We take an (open) collar
neighborhood of $\pa M$ which is diffeomorphic to $(0,2) \times \pa M$.
Let $\varphi \in C^\infty(0,2)$ be a nondecreasing function such that
$\varphi(x) = x$ if $x \le 1/2$ and $\varphi(x) = 1$ if $x \ge 3/2$. Given
$t > 0$, consider a Riemannian metric on $\interior(M)$ whose restriction to
$(0,2) \times \partial M$ is
$$g_M  =  t^{-2}  dx^2  + \varphi^2(x)  g_{\partial M}. \Eq(metric)$$ 

We have a triple $(Q, D, \tau)$ for $\Omega^*_c(M; {\cal V}^\infty)$ as in
\equ(Top), with the difference that $\tau$ is now given
on $\Omega^*_c(0,2) \, \widehat{\otimes} \, \Omega^*(\partial M)$ by
$$
\tau(1 \wedge \alpha) = dx \wedge t^{-1} \varphi(x)^{2m-1-2|\alpha|}
\tau_{\partial M} \widehat{\alpha},
$$
$$
\tau(dx \wedge \alpha) = 1 \wedge i^{-(2m-1)} t \varphi(x) ^{2m-1-2|\alpha|}
\tau_{\partial M} \widehat{\alpha}.
$$
We define $Q_{\rm alg}$ and $D_{\rm alg}$ as in \equ(Top2).
We modify $\tau_{\rm alg}$ of \equ(Top2) to act on
$\Omega^*_c(0,2) \, \widehat{\otimes} \, \whW^*$ by
$$
\tau_{\rm alg}(1 \wedge \alpha) = dx \wedge t^{-1} 
(\varphi (x)\varphi(2-x))^{2m-1-2|\alpha|}
\tau_{\widehat{W}} 
\widehat{\alpha},
$$
$$
\tau_{\rm alg}(dx \wedge \alpha) = 1 \wedge i^{-(2m-1)} 
t (\varphi(x) \varphi(2-x))^{2m-1-2|\alpha|}
\tau_{\widehat{W}} \widehat{\alpha}.
$$
That is, we metrically cone off the algebraic complex at both
$0$ and $2$.
Then we obtain a direct sum duality operator
$\tau_C$ on $C^*=\Omega^*_c(M ; \Fl^\infty) \oplus \left(
\Omega^*_c(0,2) \, \widehat{\otimes} \, \whW^* \right)$, and a corresponding
``conic'' inner product on $C^*$.

By the definition of the ${\cal B}^\infty$-module 
associated to $\Omega^*(M;\Fli)$  (endowed 
with the conic metric above), the following maps are
isometries :
$$J'_p: C^\infty_c(0,1/2) \otimes \left( \Omega^{p-1}(\pa M;\Fli_0)\oplus
\Omega^p(\pa M;\Fli_0) \right) \rightarrow
\Omega^p_c((0,1/2) \times \partial M ; \Fli)  
$$
$$J'_p(\phi_{p-1},\phi_p) = 
\left( dx\wedge t^{-1/2} x^{p-1-(2m-1)/2} \phi_{p-1} \right) + \left( 
1 \wedge t^{1/2} x^{p-(2m-1)/2}\phi_p \right) \,.$$
Similarly, by definition, the following maps are isometries
$$\widehat{J}_p: C^\infty_c(0, 1/2) \otimes 
\left( \widehat{W}^{p-1} \oplus \widehat{W}^p \right) \rightarrow
\bigoplus_{q+r=p} \left(
\Omega^q_c(0, 1/2) \, \widehat{\otimes} \, \whW^r \right)
$$
$$\widehat{J}_p(w^{p-1},w^p) =
\left( dx\wedge t^{-1/2} x^{p-1-(2m-1)/2} w^{p-1} \right) + \left( 1 \wedge
t^{1/2} x^{p-(2m-1)/2} w^p \right) \,.$$
Put
$$J=\left( \matrix{ J' & 0 \cr 0 & \widehat{J} } \right)\,.$$
We define an ``almost'' differential $\dco$ on the conic
complex $C^*=\Omega^*_c(M ; \Fli) \oplus \left(
\Omega^*_c(0,2) \, \widehat{\otimes} \, \whW^* \right)$
by the same formula as in \equ(diffdef).
Let $C^*_{(2)}$ denote the completion of $C^*$ in the sense of
$C^*_r(\Gamma)$-Hilbert modules.

\nproclaim Lemma [regular].
$\dco$ is a regular operator in the sense of [BaJ, Definition 1.1] when
acting on $C^*_{(2)}$.

\noindent
{\it Proof.} We give a sketch of the proof and omit some
computational details.
We will use throughout the general fact that if $T$ is a
regular operator and $a$ is an adjointable 
bounded operator then $T + a$ is a regular
operator.  This is proven in [RW, Lemma 1.9] when $T$ and $a$ are self-adjoint,
but one can check that the proof goes through without this additional
assumption. In addition, it follows from [H, Lemme 2.1] that a 
compactly-supported change in
the Riemannian metric does not affect the regularity question. 
Hence, for simplicity,
we will only specify our Riemannian metrics up to a compactly-supported
perturbation.

We define three new complexes. Put 
$$C_1^* = \left( \Omega^*_c(-\infty, \infty) \,
\widehat{\otimes} \, \Omega^*(\partial M; \Fl_0) \right) \oplus
\left( \Omega^*_c(-\infty, \infty) \,
\widehat{\otimes} \, \widehat{W}^* \right),$$
$$C_2^* = \left( \Omega^*_c(0, \infty) \,
\widehat{\otimes} \, \Omega^*(\partial M; \Fl_0) \right) \oplus
\left( \Omega^*_c(0, \infty) \,
\widehat{\otimes} \, \widehat{W}^* \right),$$
$$C_3^* =  \Omega^*_c( ((- \infty, 0] \times \partial M) \cup_{\partial M}
M; \Fl) \oplus
\left( \Omega^*_c(-\infty, 2) \,
\widehat{\otimes} \, \widehat{W}^* \right).$$
For $1\leq i \leq 3$, the differentials $D_{C_i}$
will roughly be of the form \equ(diffdef),
but the ``coupling'' between the geometric and algebraic subcomplexes will
depend on $i$. Namely, $D_{C_1}$ will be uncoupled on $(-\infty, 0)$
and fully coupled on $(1/4, \infty)$. The differential $D_{C_2}$ will always
be fully coupled.  The differential $D_{C_3}$ will always be completely 
uncoupled.
The metric on $C^*_1$ will be product-like on $(-\infty, 0)$ and conic on
$(1/2, \infty)$. The metric on $C^*_2$ will be fully conic on
$(0, \infty)$. The metric on $C^*_3$ will be product-like on 
$(- \infty, 0)$, and conic on its algebraic part for $(3/2, 2)$.
Let $C^*_{i, (2)}$ be the completion of $C^*_i$ in the sense of
$C^*_r(\Gamma)$-Hilbert modules. By abuse of notation, we will also
let $D_{C_i}$ denote the densely-defined differential on
$C^*_{i, (2)}$.

We claim that $D_{C_2}$ is regular, when acting on $C^*_{2, (2)}$. 
It is not hard to show that $D_{C_2}$ is closable and that
$D_{C_2}^*$ is densely defined. It remains to show that
$I + D^*_{C_2} D_{C_2}$ is surjective. To see
this, we can use separation of variables and adapt the functional calculus
of [Ch, Section 3] to our setting. That is, it is possible to write down an
explicit inverse to $I + D^*_{C_2} D_{C_2}$. 
Namely, as in [Ch, p. 586], we can span
$\Omega^*(0, \infty) \,
\widehat{\otimes} \, \left( \Omega^*(\partial M; \Fl_0) \oplus
\widehat{W}^* \right)$
by forms of type 1-4, E and O.  As in [Ch, p. 587], the 
operator $I + D^*_{C_2} D_{C_2}$ acts as the identity on forms of type 
2, 4 and O and as $I + D^*_{C_2} D_{C_2} + D_{C_2} D_{C_2}^*$ on forms of
type 1, 3 and E. Then using the equivalent of [Ch, (3.37) and (3.40)], 
one can write down an
explicit inverse to $I + D^*_{C_2} D_{C_2}+ D_{C_2} D_{C_2}^*$ when acting
on forms of type 1, 3 and E.
(The spectrum of the transverse Laplacian is discrete in
[Ch], but in our case the spectrum of the Laplacian on 
$\Omega^*(\partial M; \Fl_0) \oplus
\widehat{W}^*$ is generally not discrete.  Thus one must make the
notational change of replacing the
eigenvalue sums in [Ch] by a functional calculus.)
This proves the claim.

We claim that $D_{C_3}$ is also regular, when acting on
$C^*_{3, (2)}$. To see this, let
$d_{C_3}$ denote the differential on 
$\Omega^*_{(2)}( ((- \infty, 0] \times \partial M) \cup_{\partial M} M)$.
It follows from the analysis in [APS] that $d_{C_3}$ is regular. 
Hence for any $N \in \NN$, $d_{C_3} \otimes \id_N$ is regular when
acting on 
$\Omega^*_{(2)}( ((- \infty, 0] \times \partial M) \cup_{\partial M} M)
\otimes C^*_r(\Gamma)^N$.
Now we can find some
$N > 0$ and a projection $p \in C^\infty(M; M_N(C^*_r(\Gamma)))$ so that
$\Fl = \IM(p)$. Taking $p$ to be a product near $\partial M$, we can
extend it to a projection 
$p \in C^\infty(((- \infty, 0] \times \partial M) \cup_{\partial M} M; 
M_N(C^*_r(\Gamma)))$. As $p (d_{C_3} \otimes \id_N) p + (1-p)
(d_{C_3} \otimes \id_N) (1-p)$ differs from $d_{C_3} \otimes \id_N$ by
an adjointable 
bounded operator, it follows that $p (d_{C_3} \otimes \id_N) p + (1-p)
(d_{C_3} \otimes \id_N) (1-p)$ is regular.  Hence 
$p (d_{C_3} \otimes \id_N) p$ is regular. Now 
$p (d_{C_3} \otimes \id_N) p$ differs from the differential on 
$\Omega^*_{(2)}( ((- \infty, 0] \times \partial M) \cup_{\partial M}
M; \Fl)$ by an adjointable
bounded operator.  Finally, one can show by hand that
the differential on $\left( \Omega^*_{c}(-\infty, 2) \,
\widehat{\otimes} \, \widehat{W}^* \right)_{(2)}$ is regular. Thus
$D_{C_3}$ is regular.

We now define a certain unitary operator $U$ from
$C_{2, (2)}^* \oplus C_{3, (2)}^*$ to $C_{(2)}^* \oplus C_{1, (2)}^*$. 
The construction of $U$ is as in [Bu, Section 3.2], with some
obvious changes in notation. We refer to [Bu, Section 3.2] for the details.
Clearly $U (D_{C_2} \oplus D_{C_3}) U^{-1}$ is regular.
From the method of construction of [Bu, Section 3.2], one sees that
$U (D_{C_2} \oplus D_{C_3}) U^{-1}$ differs from
$\dco \oplus D_{C_1}$ by an adjointable bounded operator. Hence 
$\dco \oplus D_{C_1}$ is regular when acting on 
$C_{(2)}^* \oplus C_{1, (2)}^*$. In particular, $\dco$ is regular when
acting on $C_{(2)}^*$. 
\QED
\medskip
The perturbed conic signature operator
${\cal D}_C^{{\rm sign,cone}} = \dco + (\dco)^*$ satisfies
$${\cal D}_C^{{\rm sign,cone}} = \dco - \tau \dco \tau.$$
A straightforward calculation shows that
on the part of $C^*$ corresponding to $x \in (0, 1/4)$, we have
$$\left( J^{-1} \, {\cal D}_{C}^{{\rm sign,cone}} \, J \right)^+ =
\Theta \cdot \Phi
\left( 
t \left( \partial_x + { m-{1 \over 2}- {\rm degree} \over
x} \right) + {{\cal D}_C^{\rm sign}(1)_0 \over x}
\right)
\Phi^{-1},
\Eq(norform)$$
where $\Phi$,
$\Theta$ and 
${\cal D}_C^{\rm sign}(1)_0$ are as in Lemma \thm[ident.2], and
``degree'' is the $\ZZ$-grading operator. 
As ${\cal D}_C^{\rm sign}(1)_0$ is invertible, we can evidently
choose a $t> 0$ small enough
so that for any $s\in [0,1]$,
$${\rm Spec}\left(
s ( m-{1 \over 2}- {\rm degree})
+ t^{-1} {\cal D}_C^{\rm sign}(1)_0
\right) \,\cap\,(-1,1) = \emptyset\,. \Eq(speccon)
$$
In the rest of this section, we will fix such a number $t$.

\nproclaim Proposition [fred]. For $t>0$ small enough the triple $\left(
C^*_{(2)}, Q_C, D_C^{\rm cone} \right)$ defines an element of
${\bf L}_{nb}(C^*_r(\Gamma))$ in the sense of
[HS, D\'efinition 1.5].

  \noindent
{\it Proof.} We have shown in Lemma \thm[regular] 
that $D_C^{\rm cone}$ is regular.
We must show in addition that \hfil\break 
1. $(D_{C}^{\rm cone})^\prime + D_C^{\rm cone}$ is $C^*_r(\Gamma)$-bounded.
\hfil\break
2. $(D_C^{\rm cone})^2$ is $C^*_r(\Gamma)$-bounded.
\hfil\break
3. There are $C^*_r(\Gamma)$-compact operators $S$ and $T$ such that
$S D_C^{\rm cone}$ is $C^*_r(\Gamma)$-bounded, $\IM(T) \subset
\Dom(D_C^{\rm cone})$, $D_C^{\rm cone} T$ is $C^*_r(\Gamma)$-bounded and
$S D_C^{\rm cone} + D_C^{\rm cone} T - I$ is $C^*_r(\Gamma)$-compact.

For 1., we have
$(D_{C}^{\rm cone})^\prime + D_C^{\rm cone} = 0$. For 2., we have
$$(D_C^{\rm cone})^2 = 
\cases{
\left(
\matrix{0 & \# \, dx \, (\partial_x\phi) \, 
\widehat{g}  \cr
0 & 0 } 
\right) &if$\;* < m - \ha$, \cr
&\cr
&\cr
\left(
\matrix{0 & 0 \cr
\# \, dx \, (\partial_x\phi) \, \widehat{f} & 0 } 
\right) & if$\;* > m - \ha$
,\cr
}
\Eq(square)$$
where $\#$ denotes a power of $i$.
This is clearly a bounded operator.

Following [BS, Section 2], for $f\in L^2(0, \infty)$, put
$$
P_1(s)[f](x) = \int_0^x\, (y/x)^s f(y) dy,\quad s > -1/2
$$
$$
P_2(s)[f](x) = \int_1^x\, (y/x)^s f(y) dy,\quad s < 1/2
$$ Let $\phi_1(s),\, \phi_2(s) \in C^\infty(\RR)$
be such that $\phi_1(s) = 1$ for $s \geq 1$,
 $\phi_1(s) = 0$ for $s \leq 1/2$,  $\phi_2(s) = 1$ for $s \leq -1$ and
 $\phi_1(s) = 0$ for $s \geq -1/2$. Moreover, put
$$
X= {1\over t}{\cal D}_C^{\rm sign}(1)_0 + (m - {1 \over 2} - {\rm degree}).
$$
There is a standard interior parametrix for ${\cal D}_{C}^{{\rm sign,cone}}$.
Furthermore, as in [BS, Theorem 2.1],
$$
t^{-1} \left( P_1( \phi_1 (X ) )\, + \, P_2( \phi_2 (X ) ) \right)
$$ 
is a parametrix for $t (\partial_x + X)$ on 
$(0, 1/4) \times \partial M$.
Finally, if $z = 2 - x$ then with an evident notation, when acting on 
$\Omega^*_c(7/4, 2) \, \widehat{\otimes} \, \widehat{W}^*$,
we can write
$$\left( J^{-1} \, {\cal D}_{C}^{{\rm sign,cone}} \, J \right)^+ =
\Theta \cdot \Phi
\left( 
t \left( \partial_z + { m-{1 \over 2}- {\rm degree} \over
z} \right) + {{\cal D}_{\widehat{W}}^{\rm sign} \over z}
\right)
\Phi^{-1}.
\Eq(norform2)$$
We remark that, as in [BS, Lemma 5.4],
the middle-dimensional vanishing in \equ(newcomplex)
ensures that the conic operator \equ(norform2) exists without a further
choice of boundary condition. 
Put
$$
X^\prime = {1\over t}{\cal D}_{\widehat{W}}^{\rm sign} + 
(m - {1 \over 2} - {\rm degree}).
$$ 
Then
a parametrix for $t (\partial_x + X^\prime)$
on $\Omega^*_c(7/4, 2) \, \widehat{\otimes} \, \widehat{W}^*$ is given by
$$
t^{-1} \left( P_1( \phi_1 (X^\prime) )\, + \, P_2( \phi_2 (X^\prime) ) \right).
$$ 
One constructs an 
(adjointable)
parametrix $G$ for
${\cal D}_{C}^{{\rm sign,cone}}$ by patching 
these three parametrices together, using
\equ(norform) and \equ(norform2).
Put $S = T = G (D_C^{\rm cone})^* G$.
The proposition follows.
\QED

\medskip

From [HS, Proposition 1.6], 
the conic signature
operator defines a higher index class 
$\Ind({\cal D}_{C}^{{\rm sign,cone},+})\in
 K_0(C^*_r(\Gamma))$
which depends neither on the choice
of Riemannian metric on $M$  nor on 
$t$ (provided that $t$ is sufficiently small
for the constructions to make sense).  
As in [L5, Proposition 14], it is also independent
of the choices of $\phi$, 
$W$ and  the homotopy equivalence $f$.

%\fine
\expandafter\ifx\csname sezioniseparate\endcsname\relax%
   \input macro \fi

%--------------------------------------------------------
\numsec=6
\numfor=1
\numtheo=1
\pgn=1

\beginsection 6. Homotopy invariance of the conic index class

We keep the notation of Section 4. In this section alone, 
we put ${\cal B}^\infty = C^*_r(\Gamma)$.
Let $M_1$ and $M_2$ be compact oriented manifolds-with-boundary.  Suppose that
we have oriented-homotopy equivalences $h_1 : (M_1, \partial M_1) \rightarrow 
(M_2, \partial M_2)$ and
$h_2 : (M_2, \partial M_2) \rightarrow (M_1, \partial M_1)$ 
which are homotopy inverses to each other.
We can homotop $h_1$ and $h_2$ to assume that they are product-like
near the boundaries.  That is, 
for $i \in \{ 1,2 \}$, put $\partial h_i = h_i \big|_{\partial M_i}$.
Then when restricted to the collar neighborhood
${\cal U}_i =  (0,2) \times \partial M_i$, 
we assume that $h_i(x, b_i) = (x, \partial h_i(b_i))$ for $x \in (0,2)$ and
$b_i \in \partial M_i$.

We assume that $\partial M_1$ and $\partial M_2$ satisfy {\bf Assumption 1}.
Let $\widehat{W}^*_i$ be cochain complexes as in \equ(newcomplex), with
corresponding maps
$$\widehat{f}_i : \Omega^*(\partial M_i; ({\cal V}_i)_0) \rightarrow
\widehat{W}^*_i$$
and
$$\widehat{g}_i : \widehat{W}^*_i \rightarrow
\Omega^*(\partial M_i; ({\cal V}_i)_0).
$$

We would like to compare $\Omega^*_{(2)}(M_2; {\cal V}_2)$ with
$\Omega^*_{(2)}(M_1; {\cal V}_1)$ using the maps $h_i^*$, but there is the
technical problem that $h_i^*$, as originally defined on smooth forms,
need not be $L^2$-bounded if $h_i$ is not
a submersion. As in [HS, p. 90], 
we modify $\{h_i^*\}_{i=1}^2$ to obtain $L^2$-bounded
cochain homotopy equivalences between 
$\Omega^*_{(2)}(M_2; {\cal V}_2)$ and
$\Omega^*_{(2)}(M_1; {\cal V}_1)$ as follows. 
From [HS, p. 90], for suitably large $N$,
there is a submersion
$H_i : B^N \times M_i \rightarrow M_{3-i}$ such that
$H_i(0, m_i) = h_i(m_i)$. Here $B^N$ is an open ball
in an euclidean space of dimension $N$.
Furthermore, from the construction in 
[HS, p. 90], we may assume that $H_i$ is product-like near $\partial M_i$.
Fix $v \in \Omega^N_c(B^N)$ with $\int_{B^N} v = 1$.
Define a bounded cochain homotopy equivalence
$$t_i :
\Omega^*_{(2)}(M_i; {\cal V}_i) \rightarrow \Omega^*_{(2)}(M_{3-i}; 
{\cal V}_{3-i})$$
by
$t_i(\omega) = \int_{B^N} v \wedge H^*_{3-i}(\omega)$.
Let $\partial t_i$ be the analogous map from
$\Omega^*_{(2)}(\partial M_i; ({\cal V}_i)_0)$ to 
$\Omega^*_{(2)}(\partial M_{3-i}; ({\cal V}_{3-i})_0)$.

As $f_i$ and $g_i$ are homotopy inverses, there are bounded operators
$$A_i : \Omega^*(\partial M_i; ({\cal V}_i)_0) \rightarrow
\Omega^{*-1}(\partial M_i; ({\cal V}_i)_0)$$
and 
$$B_i : W^*_i \rightarrow W^{*-1}_i$$
 such that
$$I - g_i \circ f_i = D_{\partial M_i} A_i + A_i D_{\partial M_i}$$
and
$$I - f_i \circ g_i = D_{W_i} B_i + B_i D_{W_i}.$$
As $(f_i)^\prime = g_i$, $(D_{\partial M_i})^\prime = -D_{\partial M_i}$
and $(D_{W_i})^\prime = -D_{W_i}$, we can assume that
$(A_i)^\prime = - A_i$ and $(B_i)^\prime = - B_i$.
Let $\widehat{B}_i$ denote the obvious extension of $B_i$ to a map from
$\widehat{W}^*_i$ to $\widehat{W}^{*-1}_i$.
Put
$$C^*_i = \Omega^*_{(2)}(M_i; {\cal V}_i) \, \oplus \, 
\left( \Omega^*_{(2)}(0,2) \, \widehat{\otimes} \, \widehat{W}^*_i \right).$$

\nproclaim Theorem [conichominv]. 
The index class $\Ind({\cal D}^{C,+}_{{\rm sign,cone}})\in
 K_0(C^*_r(\Gamma))$ is the same for $M_1$ and $M_2$.

\noindent
{\it Proof.} We will show that [HS, Lemme 2.4] applies.
From [HS, p.90-91], there are bounded operators 
$y_i : \Omega^*(M_i; {\cal V}_i) \rightarrow \Omega^{*-1}(M_i; {\cal V}_i)$ 
such that
$$1 - t_i^\prime t_i = D_{M_i} y_i + y_i D_{M_i}.$$
Similarly, there are bounded operators
$z_i : \Omega^*(M_i; {\cal V}_i) \rightarrow \Omega^{*-1}(M_i; {\cal V}_i)$ 
such that
$$1 - t_i t_i^\prime = D_{M_i} z_{3-i} + z_{3-i} D_{M_i}.$$
(This follows from the proof of [HS, Proposition 2.5] in the case 
$\alpha = 0$, with
the operator $S$ being the analogue of [HS, p. 94 b -7].)
We can assume that $y_i$ and $z_i$ are product-like near
$\partial M_i$. Let $\partial y_i$ and $\partial z_i$ denote their
boundary restrictions.

Define $T_i : C^*_{i} \rightarrow C^*_{3-i}$ by
$$T_i = 
\cases{
\left(
\matrix{t_i & \phi(x) \, 
\widehat{\otimes} \, (A_{3-i} \circ \partial t_i \circ 
\widehat{g}_i)  \cr
0 & \id \, \widehat{\otimes} \,
(\widehat{f}_{3-i} \circ \partial t_i \circ \widehat{g}_i)}
\right) &if$\;* < m - \ha$, \cr
&\cr
&\cr
\left(
\matrix{t_i & 0 \cr
- \phi(x) \, \widehat{\otimes} \, (\widehat{f}_{3-i} \circ \partial t_i \circ 
{A}_i) &  
\id \, \widehat{\otimes} \,
(\widehat{f}_{3-i} \circ \partial t_i \circ \widehat{g}_i)} 
\right) & if$\;* > m - \ha$
.\cr
}
\Eq(T_i)$$ 
Define
$Y_i : C^*_i \rightarrow C^{*-1}_i$ by
$$Y_i = 
\cases{
\left(
\matrix{y_i & \phi(x) \, 
\widehat{\otimes} \, \widehat{\alpha}_i  \cr
0 & - \id \, \widehat{\otimes} \,
(\widehat{B}_i + \widehat{f}_i \circ \partial y_i \circ \widehat{g}_i
+ \widehat{f}_i \circ (\partial t_i)^\prime \circ A_{3-i} \circ 
\partial t_i \circ \widehat{g}_i)
}
\right) &if$\;* < m - \ha$, \cr
&\cr
&\cr
\left(
\matrix{y_i & 0 \cr
- \phi(x) \, \widehat{\otimes} \, \widehat{\beta}_i & 
- \id \, \widehat{\otimes} \,
(\widehat{B}_i + \widehat{f}_i \circ \partial y_i \circ \widehat{g}_i
+ \widehat{f}_i \circ (\partial t_i)^\prime \circ A_{3-i} \circ 
\partial t_i \circ \widehat{g}_i)
}
\right)
& if$\;* > m - \ha$
,\cr
}
\Eq(Y_i)$$
where $$\alpha_i : {W}^*_i \rightarrow \Omega^{*-2}(\partial M_i; 
({\cal V}_i)_0)$$
is a map which satisfies
$D_{\partial M_i} \alpha_i - \alpha_i D_{W_i} = g_i B_i - A_i g_i$,
$\widehat{\alpha}_i$ is its extension to $\widehat{W}^*_i$,
$$\beta_i : \Omega^{*}(\partial M_i; 
({\cal V}_i)_0) \rightarrow {W}^{*-2}_i$$
is a map which satisfies
$D_{W_i} \beta_i - \beta_i D_{\partial M_i} = - B_i f_i + f_i A_i$
and $\widehat{\beta}_i$ is the  extension of ${\beta}_i $ to $\widehat{W}^*_i$.
Here $\alpha_i$ and $\beta_i$ exist because of Lemma \thm[complex], 
and from the proof of Lemma \thm[complex],
we can take them to be continuous.
Define
$Z_i : C^*_i \rightarrow C^{*-1}_i$ by
$$Z_i = 
\cases{
\left(
\matrix{z_i & \phi(x) \, 
\widehat{\otimes} \, \widehat{\alpha}_i  \cr
0 & - \id \, \widehat{\otimes} \,
(\widehat{B}_i + \widehat{f}_i \circ \partial z_i \circ \widehat{g}_i
+ \widehat{f}_i \circ \partial t_{3-i} \circ A_{3-i} \circ 
(\partial t_{3-i})^\prime \circ \widehat{g}_i)
}
\right) &if$\;* < m - \ha$, \cr
&\cr
&\cr
\left(
\matrix{z_i & 0 \cr
- \phi(x) \, \widehat{\otimes} \, \widehat{\beta}_i & 
- \id \, \widehat{\otimes} \,
(\widehat{B}_i + \widehat{f}_i \circ \partial z_i \circ \widehat{g}_i
+ \widehat{f}_i \circ \partial t_{3-i} \circ A_{3-i} \circ 
(\partial t_{3-i})^\prime \circ \widehat{g}_i)
}
\right)
& if$\;* > m - \ha$
.\cr
}
\Eq(Z_i)$$
Define ${\cal E}_i : C^*_i \rightarrow C^{*}_i$ to be
$(-1)^{\rm degree}$.

We will think of the small positive number $t$ in the metric
\equ(metric) as a free parameter. One can check that the operators 
$T_i$, $Y_i$ and $Z_i$ are bounded and
have norms which are independent of $t$. In order to apply 
[HS, Lemme 2.4], it suffices to show that if $t$ is made small then
$(D_{C,i}^{\rm cone})^2$, 
$1 - T_i^\prime T_i - D_{C,i}^{\rm cone} Y_i - Y_i D_{C,i}^{\rm cone}$ and
$1 - T_{3-i} T_{3-i}^\prime - D_{C,i}^{\rm cone} Z_{i} - Z_{3} 
D_{C,i}^{\rm cone}$
can be made arbitrarily small in norm.

The formula for $(D_{C,i}^{\rm cone})^2$ was given in \equ(square).
One computes that
$$
\eqalign{
1 - T_i^\prime T_i - D_{C,i}^{\rm cone} Y_i - Y_i D_{C,i}^{\rm cone} & =
1 - T_{3-i}T_{3-i}^\prime - D_{C,i}^{\rm cone} Z_{i} - Z_{i} 
D_{C,i}^{\rm cone} \cr
& =
\cases{
\left(
\matrix{0 & \# dx \, (\partial_x \phi) \, 
\widehat{\otimes} \, \widehat{\alpha}_i  \cr
0 & 0
}
\right) &if$\;* < m - \ha$, \cr
&\cr
&\cr
\left(
\matrix{0 & 0 \cr
\# dx \, (\partial_x \phi) \, \widehat{\otimes} \, \widehat{\beta}_i & 
0
}
\right)
& if$\;* > m - \ha$
,\cr
}
}$$
where $\#$ denotes a power of $\sqrt{-1}$. (To do the calculations, it
is convenient to note that for fixed $x \in (0,2)$, the restriction of
$D_{C,i}^{\rm cone}$ to $\{x\} \times \partial M_i$ is of the form 
\equ(diffp) with $F = \partial M_i$ and $\epsilon = \phi(x)$.)

In each case, the result is of the form
$dx \otimes K(x)$ where $K(x)$ is a bounded operator which is only nonzero
when $x \in \left( {1 \over 4}, {1 \over 2} \right)$. 
From \equ(metric), the norm of $dx$ in this region is $t$, while the norm of 
$K(x)$ is independent of $t$. The theorem follows. \QED

%\fine
\expandafter\ifx\csname sezioniseparate\endcsname\relax%
   \input macro \fi

%--------------------------------------------------------
\numsec=7
\numfor=1
\numtheo=1
\pgn=1

\beginsection 7. Equality of the conic  and  APS-index classes

We first define the generalized APS-index. 
Fix numbers $t_1, t > 0$ and consider a 
Riemannian metric on $\interior(M)$ whose restriction to $(0,2) \times
\partial M$ is
$$g_M = t_1^{-2} dx^2 + g_{\partial M}. \Eq(defmetric)$$ 
In what follows, we will think of $t_1$ as a parameter associated to
$x = 0$ and $t$ as a parameter associated to $x = 2$.
Let $\chi \in C^\infty(0,2)$ be a positive function such that
$\chi(x) = t_1 $ if $x \in (0, 1/2)$ and 
$\chi(x) = t$ if $x \in (3/2, 2)$.
Let us go through the steps to define the signature operator as in
Section 5, with the differences that
$\tau$ is now given
on $\Omega^*_c(0,2) \, \widehat{\otimes} \, \Omega^*(\partial M)$ by
$$
\tau(1 \wedge \alpha) = dx \wedge t_1^{-1} 
\tau_{\partial M} \widehat{\alpha},
$$
$$
\tau(dx \wedge \alpha) = 1 \wedge i^{-(2m-1)} t_1 
\tau_{\partial M} \widehat{\alpha}.
$$ and $\tau_{\rm alg}$ is now given
on $\Omega^*_c(0,2) \, \widehat{\otimes} \, \widehat{W}^*$ by
$$
\tau_{\rm alg}(1 \wedge \alpha) = dx \wedge \chi(x)^{-1}  
(\varphi(2-x))^{2m-1-2|\alpha|}
\tau_{\widehat{W}} 
\widehat{\alpha},
$$
$$
\tau_{\rm alg}(dx \wedge \alpha) =  1 \wedge i^{-(2m-1)} \chi(x) 
 (\varphi(2-x))^{2m-1-2|\alpha|}
\tau_{\widehat{W}} \widehat{\alpha}.
$$
\noindent
That is, metrically speaking, we have a product structure near $x = 0$ 
and a cone on the
algebraic complex near $x = 2$.
Consider the corresponding perturbed signature operator $\SC$.
It has an invertible boundary operator
$\SC(1)_0$. We define $(H^1,\Pi_>)$ to be the $C^*_r(\Gamma)$-Hilbert
module of Sobolev-$H^1$ elements $(\alpha,v)$ of $\Omega^*_{(2)}(M;\Fl)
\oplus \left( \Omega^*_{(2)}[0,2)\widehat{\otimes}\whW^* \right)$ such that
$$\Pi_>(\alpha|_{\pa M}\oplus v(0)) = 0,$$
where $\Pi_>$ is the spectral projection onto the positive part
of the operator $H$ of Lemma \thm[ident.2].

The work of Wu [Wu] shows that 
$$\SCp:(H^1,\Pi_>)^+ \rightarrow L^{2,-}$$
is a $C^*_r(\Gamma)$-Fredholm operator which thus defines
a higher index class $\Ind(\SCAPSp)\in K_0(C^*_r(\Gamma))$.

\nproclaim Remark [remark]. The higher index class $\Ind(\SCAPSp)$ 
is independent of the choices of  $t_1$, $t$ and $\chi$.
In this section we will take 
$t=t_1$ and $\chi(x) = t$. In Section 10 we will take $t_1=1$. In Section 
9 we will consider a $b-$metric, in which case we effectively have
$\chi(x)=x$ for $x \in (0, 1/2)$.

\nproclaim Theorem [c=aps]. The following equality holds in 
$K_0(C^*_r(\Gamma))$:
$\Ind({\cal D}_C^{{\rm sign,cone},+})=\Ind(\SCAPSp)$.

\smallskip
\noindent
{\it Proof.} 
Using \equ(speccon), one can construct 
a continuous family 
of elliptic conic-operators $\{ {\cal D}_C^{{\rm sign,cone},+}(s) 
\}_{s \in [0,1]}$ 
with ${\cal D}_C^{{\rm sign,cone},+}(1) =
{\cal D}_C^{{\rm sign,cone},+} $,
such that for any $s\in [0,1] $, the restriction of 
${\cal D}_C^{{\rm sign,cone},+}(s)$ to the conic end $x \in (0, 1/4)$
is given by
$$\left( J^{-1} \, {\cal D}_{C}^{{\rm sign,cone}}(s) \, J \right)^+ =
\Theta \cdot \Phi
\left( 
t \left( \partial_x + s { m-{1 \over 2}- {\rm degree} \over
x} \right) + {{\cal D}_C^{\rm sign}(1)_0 \over x}
\right)
\Phi^{-1},$$
and 
with the principal symbols in the interior also varying continuously in $s$.
This implies that ${\cal D}_C^{{\rm sign,cone},+}(s)$
is $C^*_r(\Gamma)$-Fredholm for any $s\in [0,1] $
and that $\{ {\cal D}_C^{{\rm sign,cone},+}(s) \}_{s \in [0,1]}$ 
is a continuous path of $C^*_r(\Gamma)$-Fredholm operators.
Let  $\{ {\cal D}_C^{{\rm sign,APS},+}(s) 
\}_{s \in [0,1]}$ be the corresponding family of APS-type operators
on the manifold-with-boundary $M - ([0, 1/8) \times \partial M)$,
where $M - ([0, 1/8) \times \partial M)$ has the subspace metric
induced from the conical metric on $M$.
As the index is locally constant on the space of 
$C^*_r(\Gamma)$-Fredholm operators,
it follows that
$\Ind({\cal D}_C^{{\rm sign,cone},+}(s))$
and
$\Ind({\cal D}_C^{{\rm sign,APS},+}(s))$ are each independent of
$s \in [0,1]$.

At $s = 0$, the restriction of 
${\cal D}_C^{{\rm sign,cone},+}(0)$ to the conical end is conjugate to
the operator
$\Theta (t \partial_x + x^{-1} H)$, where $H$ is as in Lemma \thm[ident.2]. If 
$\psi$ lies in the kernel of this operator then
$\psi(x) = (8x)^{-H/t} \psi(1/8)$. 
If this is to lie in the Hilbert-module domain
then, in particular,
it must be square-integrable and so we must have $\Pi_> \psi(1/8) = 0$. 
One can make a similar statement for the cokernel.
Then an argument as in the proof of [BC, Theorem 1.5] shows that
$\Ind({\cal D}_C^{{\rm sign,cone},+}(0)) =  
\Ind({\cal D}_{C}^{{\rm sign,APS},+}(0))$.
Hence $$\Ind({\cal D}_C^{{\rm sign,cone},+}) = 
\Ind({\cal D}_C^{{\rm sign,cone},+}(1)) =  
\Ind({\cal D}_{C}^{{\rm sign,cone},+}(0)) = 
\Ind({\cal D}_C^{{\rm sign,APS},+}(0)) =  
\Ind({\cal D}_{C}^{{\rm sign,APS},+}(1)).$$ Also, as in
[BC, Theorem 1.5], one can show that
$\Ind({\cal D}_{C}^{{\rm sign,APS},+}(1)) =
\Ind({\cal D}_{C}^{{\rm sign,APS},+})$. The theorem follows. \QED

%\fine
\expandafter\ifx\csname sezioniseparate\endcsname\relax%
   \input macro \fi

%--------------------------------------------------------
\numsec=8
\numfor=1
\numtheo=1
\pgn=1

\beginsection 8. The enlarged $b$-calculus

In order to prove a suitable higher index formula
we shall now   change the perturbed signature operator
introduced above so as to get an element
of an appropriate $b$-calculus. 
%allowing the 
%construction of  superconnection heat kernels.
 Since the complex $C^*$
involves the additional piece $\Omega^*_c(0,2) \, \widehat{\otimes} \,
\whW^*$, this step is slightly more complicated than in [MP1]
and [LP2]. 

 Thus, let $M$ be a manifold with boundary $\pa M$.
We denote
by $u\in C^\infty(M)$ a boundary defining function
and we fix a Riemannian metric $g$ which is product-like
in a collar neighbourhood ${\cal U}$ of $\pa M$:
$g|_{{\cal U}}=du^2 + g_{\pa M}$.
As in [Me], we add to the manifold-with-boundary $M$
a cylindrical end $(-\infty,0]_u\times \pa M$. Similarly,
we add the half-line $(-\infty,0]_{u'}$ to the interval $[0,2)$
appearing in the definition of $C^*$.
The change of variables $u=\log x$, $u'=\log x'$ compactifies
these two manifolds and brings us into the framework of the $b$-geometry
of [Me]. The original manifold is contained in the $b$-manifold
so obtained, and the same is true for the interval.
We make an abuse of notation and call $(M,g)$ and $[0,2)$
the $b$-manifolds obtained above. Notice that $g$ is now
a $b$-metric which  is product-like near the boundary:
 $g=dx^2/x^2 + g_{\pa M}$, for $ 0\leq x \leq  1/2$.

We now define an appropriate {\it enlarged}  $\Bi$-$b$-calculus.
Besides  the usual $\Bi$-$b$-calculus on $M$ and the usual
$\Bi$-$b$-calculus on $[0,2)$, with values in $\whW^*$,
this enlarged version will have
 to involve operators of the following types:
$$ P:\ub\Omega^*([0,2);\whW^*)\rightarrow \ub\Omega^*(M;\Fli)
\quad\quad 
Q:\ub\Omega^*(M;\Fli)\rightarrow \ub\Omega^*([0,2);\whW^*).\Eq(2.1)
$$
For the definition of the $\Bi$-$b$-calculus we refer the
reader to [LP1, Sect. 12] and [LP3, Appendix].

We shall only need operators of order $-\infty$, i.e.
operators which are defined by {\it smooth} Schwartz kernels
on  suitable   blown-up spaces. The  blown-up space 
(see [Me, Sect. 4.2]) corresponding to
$P$  in \equ(2.1) is
$$[M\times [0,2)_{x'}\,;\,S] \quad\quad S=\{x=x'=0\}.$$
It comes with a blow-down map $\beta: [M\times [0,2);S]
\rightarrow M\times [0,2)$. There are three boundary hypersurfaces
in the manifold-with-corners $[M\times [0,2);S]$; the front face
${\rm bf}=\beta^{-1}(S)\equiv S_+N(S)$ and the left and 
right boundaries 
$${\rm lb}=\overline{\beta^{-1}(\pa M\times [0,2))}
\quad\quad
{\rm rb}=\overline{\beta^{-1}(M\times \{0\})}.
$$
We shall require  the Schwartz kernel $K_P$ of 
$P$  to lift to $[M\times [0,2);S]$ as a smooth
$C^\infty$ section of the bundle 
$\beta^*{\cal K}$
with
$${\cal K}_{(p,x')}\equiv {\rm Hom}_{\Bi}((\ub\Lambda^*[0,2)
\widehat{\otimes} \whW)_{x'}; (\ub\Lambda^*(M)\otimes\Fli)_p).$$
The kernel $\beta^* K_P$ is also required to vanish to infinite
order at ${\rm lb}$ and ${\rm rb}$.
We shall usually employ projective coordinates $(y,s,x')$
in a neighbourhood of ${\rm bf}\subset[M\times [0,2);S]$, 
with $s=x/x'$; see [Me, Ch. 4].
Notice that the front face is
diffeomorphic to $[-1,1]\times \pa M$. 
Restriction to the front face
followed by Mellin transform along $[-1,1]$
defines the indicial family of $P$ as an entire
family of $\Bi$-linear maps
$I(P,\lambda):\CC^2\otimes\whW^*\rightarrow \Omega^*(\pa M,\Fli_0)$,
$\lambda\in \CC$, with $\CC^2\cong \ub\Lambda[0,2)|_{x'=0}$;
see [Me, Sect. 5.2, formula (5.13)]. 
Using the Paley-Wiener theorem this construction can be reversed
(for smoothing operators)
as in [Me, Theorem 5.1 and Lemma 5.4].

Operators like $Q$ in \equ(2.1) are defined in a similar
way. They are integral operators 
with a Schwartz kernel  on $[0,2)\times M$
which lifts to become smooth  on the blown-up space 
$[[0,2)_{x'}\times M\,;\,\{x=x'=0\}]$ and vanishes
to infinite order at the left and right boundaries.
\bigskip

\noindent
{\bf Example.}
We shall now exibit two particular operators as in \equ(2.1).
These are $b$-version of the operators $\phi\widehat{f}$
and $\phi\widehat{g}$ already considered in
the previous section. 

Let $\phi\in C^{\infty}[0, \infty)$ be a nonincreasing function which is 
equal to $1$ on $[0,1/4)$ and equal to $0$ on $[1/2,\infty)$.
%Let $\alpha>>1$, consider $\phi_\alpha$
%defined by $\phi_\alpha(x)=\phi(x\alpha)$
%and its restriction to $[0,2)$. 
Let ${\cal U}\equiv [0,2)\times\pa M$ 
be a collar neighbourhood of $\pa M$.
As usual, we identify $\ub\Lambda^{j+1}_{\pa M}(M)$ with 
$\Lambda^{j+1}(\pa M)\oplus (\Lambda^j(\pa M)\wedge dx/x)$.
As already explained, using this identification and a trivialization
 $\Fli|_{\cal U}\cong [0,2)\times \Fli_0$,
we can write each element in $\ub\Omega^{j+1}({\cal U},\Fli)$ 
as 
$k^0\cdot \gamma^{j+1}(x) + (k^1 \wedge \gamma^j)$, with
$k^\ell\in\ub\Omega^\ell[0,2)$ and 
$\gamma^j\in\Omega^j(\pa M,\Fli_0)$.
Now let $\rho\in C^\infty_c(\RR)$ be an
even nonnegative test function such that $\int_\RR \rho(t)dt
=1$. Consider, as in [MP1, Lemma 9], the function   $\rho_\delta$
defined by $\rho_\delta(\lambda)=\delta^{-1}\cdot\rho(\lambda/\delta)$ for
$\delta >0$. For  $j<m$, consider
an element of $(\ub\Omega^*[0,2)\widehat{\otimes} \whW)^j$ of the form
$\omega= (h^0\otimes w^j) + (h^1\otimes w^{j-1})$.
Using the above identification, we define an operator
$$\widehat{g}_b
%^{\ell,j}
:(\ub\Omega [0,2)\widehat{\otimes} \whW)^j\rightarrow
\ub\Omega^{j+1}(M;\Fli)$$ 
by
$$\widehat{g}_b(\omega)
%(h^\ell\otimes w^j)
= \,\phi\cdot\int_0^\infty
\rho_\delta
(\log s) \phi(x/s) h^0(x/s) {ds\over s} \wedge \widehat{g}( w^{j}) \,
- \, i \, \phi\cdot\int_0^\infty
\rho_\delta
(\log s) \phi(x/s) h^1(x/s) {ds\over s} \wedge \widehat{g}( w^{j-1})
$$
with 
$\widehat{g} :
\widehat{W}^*
\rightarrow \Omega^*(\pa M;\Fli_0)$ as in Section 2.

Similarly, if $j> m$, we are going to define a $\Bi$-$b$-smoothing
operator
$$\widehat{f}_b:
\ub\Omega^j(M;\Fli)\rightarrow (\ub\Omega^*[0,2)
\widehat{\otimes} \whW^*)
^{j+1}
$$
as follows.
If $\omega\in\ub\Omega^j(M;\Fli)$ has
$\omega|_{{\cal U}} =
h^0  \omega^j + (h^1\wedge\omega^{j-1})$ with $h^\ell\in 
\ub\Omega^\ell[0,2)$ and $\omega^{j-\ell}
\in \Omega^{j-\ell}(\pa M;\Fli_0)$,  define
$$\widehat{f}_b (\omega)=
\phi\cdot
\int_0^\infty \rho_\delta(\log s)\phi(x/s) h^0 (x/s){ds\over s}\otimes
\widehat{f}(\omega^{j}) \, - \, i \, 
\phi\cdot
\int_0^\infty \rho_\delta(\log s)\phi(x/s) h^1 (x/s){ds\over s}\otimes
\widehat{f}(\omega^{j-1})$$
where  $\widehat{f} :
\Omega^*(\pa M;\Fli_0)\rightarrow \widehat{W}^*$ is as
in Section 2.

\medskip
\noindent
{\bf Remark.} The operators $\widehat{g}_b$ and $\widehat{f}_b$ also depend
on the choice of $\rho_\delta$. They are
 $\Bi$-$b$-smoothing operators of the type
described above, namely as
in \equ(2.1).

\medskip

\nproclaim Definition [enlarged]. 
The enlarged (small) $\Bi$-$b$-calculus
of order $m$, denoted $\widehat{\Psi}^m_{b,\Bi}$, is the space
of operators  $$P=\left( \matrix{P_{11} & P_{12} \cr P_{21} & P_{22}}
\right) \quad{\rm acting}\quad 
{\rm on}\quad \ub\Omega^*(M;\Fli)\oplus (\ub\Omega^*[0,2)
\widehat{\otimes}
\whW^*)$$
with $P_{11}$ and $P_{22}$ being $\Bi$-$b$-pseudodifferential
operators as in [LP1, Sect. 12] and $P_{12}$, $P_{21}$ as in \equ(2.1)
above. 

By construction, the  Schwartz kernels of $P$ will vanish
to infinite order  at the right and left boundaries.
Dropping this last condition but assuming conormal bounds
there, one can define, as in [Me, Sect. 5.16] [LP1, Sect. 12], the calculus with
bounds $\widehat{\Psi}^{m,\beta}_{b,\Bi}$ for $\beta > 0$.
Following the previous sections, we
shall  now consider a new differential $D_C$ on the perturbed
complex $C^*=\ub\Omega^*(M;\Fli)\oplus (\ub\Omega^* [0,2)
\widehat{\otimes}
\whW)$; on the degree $j$-subspace we put
$$D_C\equiv\left(\matrix{ D_M & 0 \cr
0 & D_{{\rm alg}}} \right)+\cases{\left(\matrix{0 & \widehat{g}_b \cr
0 & 0 } \right)&if$\;j<m$\cr \left(\matrix{0 & 0\cr
-\widehat{f}_b & 0 } \right)&if$\;j> m$\cr}\Eq(differ)$$
where $\widehat{g}_b$ and $\widehat{f}_b$ are defined as above.
By construction,
$D_C\in\widehat{\Psi}^1_{b,\Bi}$. Notice that $(D_C)^2$ is nonzero. 

Let $\SCb=D_C + (D_C)^*$
%# New macro \SCb = {\cal D}^C_{{\rm sign},b}
be the $b$-signature operator
associated to the $b$-complex  $(C^*,D_C)$.
Then  $\SCb=D_C- \tau_C D_C \tau_C$ 
is odd with respect to the $\ZZ_2$-grading
defined by the Hodge duality operator $\tau_C$ on $C^*$
(see Section 4). We shall call $\SCb$ the {\it perturbed}
$b$-signature operator.
More explicitly, on forms of degree
$m$, $\SCb$ is equal to 
$$\left( \matrix{ D_M -\tau_M D_M \tau_M& 0 \cr
0 & 
D_{{\rm alg}}-\tau_{{\rm alg}} D_{{\rm alg}} \tau_{{\rm alg}} }\right),
$$ 
whereas on forms of degree $j\not= m$, using Lemma \thm[antil],
$$\SCb=\left( \matrix{ D_M -\tau_M D_M \tau_M& 0 \cr
0 & D_{{\rm alg}}-\tau_{{\rm alg}} D_{{\rm alg}} \tau_{{\rm alg}} } \right)+
\cases{\left(\matrix{0 & \widehat{g}_b \cr
 \tau_{\whW} \, \widehat{f}_b \, \tau_{\pa M}& 0 }
\right)&if$\;j<m$\cr \left(\matrix{0 & -\tau_{\pa M} \, \widehat{g}_b
\,\tau_{\whW}\cr
-\widehat{f}_b & 0 } \right)&if$\;j>m$\cr}.
$$
The perturbed $b$-signature operator $\SCb$ is an element
of the enlarged $\Bi$-$b$-calculus defined above. 

\smallskip
\noindent
{\bf Notation.} The perturbed signature operator
depends both on  the choice of the functions $\rho$, $\phi$
and on the real number $\delta$.
For brevity, we shall write 
$$\SCb=\left( \matrix{ D_M -\tau_M D_M \tau_M& 0 \cr
0 & D_{{\rm alg}}-\tau_{{\rm alg}} D_{{\rm alg}} \tau_{{\rm alg}}} \right)+
\left(\matrix{0 & S(\delta) \cr T(\delta) & 0}
\right).
$$

\medskip
Notice that by employing the bundles $\Fl$,
$\Fl_0$ and by requiring the maps in the definition of $P_{12}$,
$P_{21}$ to be $C^*_r(\Gamma)$-linear and the operators $P_{11}$, $P_{22}$
to be in $\Psi^m_{b,C^*_r(\Gamma)}$ (see [LP1, Sect. 11]),
 one can define in a similar
way an enlarged $C^*_r(\Gamma)$-$b$-calculus, denoted
$ \widehat{\Psi}^*_{b,C^*_r(\Gamma)}$.

%\fine
\expandafter\ifx\csname sezioniseparate\endcsname\relax%
   \input macro \fi

%--------------------------------------------------------
\numsec=9
\numfor=1
\numtheo=1
\pgn=1

\beginsection 9. The  $b$-index class

We shall now show that under the present
assumptions, the operator $\SCb$ defines an index class
$\Ind(\SCbp)\in K_0(\Bi)$.

%#New macro : \SCbp={\rm D}^{C,+}_{{\rm sign},b} 

We first construct a parametrix for $\SCb$ with
$C^*_r(\Gamma)$-compact remainder.
\medskip
\noindent
{\bf Remark : } It does not seem to be mentioned in the literature that
the definition of a $C^*_r(\Gamma)$-compact operator in
[MF, Section 2] differs from that in [Ks, Definition 4].
The $C^*_r(\Gamma)$-compact operators of [MF, Section 2] form a left-ideal in
the $C^*_r(\Gamma)$-bounded operators, whereas the 
the $C^*_r(\Gamma)$-compact operators of [Ks, Definition 4] form a $2$-sided
ideal. In fact, the $C^*_r(\Gamma)$-compact operators of [Ks, Definition 4]
consist of the adjointable operators $K$ for which both $K$ and $K^*$ are
compact in the sense of [MF, Section 2]. 
In [MF, Theorem 2.4] it is implicitly assumed that the
$C^*_r(\Gamma)$-compact operators, as defined in [MF, Section 2], 
form a $2$-sided
ideal.  (The mistake is in the sentence ``Without loss of generality...'')
Hence there is a gap in the proof of [MF, Theorem 2.4].
However, it is easy to correct this by using the definition of
$C^*_r(\Gamma)$-compact operators from 
[Ks, Definition 4] throughout the paper [MF].  Then the results of 
[MF, Section 2] go through automatically and one can check that the
claims of [MF, Section 3] remain valid.  This definition of
$C^*_r(\Gamma)$-compact operators should also be used in
[LP1]-[LP4]. \medskip

%Proceeding in an abstract way let $P\in\widehat{\Psi}^m_{b,C^*_r(\Gamma)}$
%be an element in the enlarged $b$-calculus. 
The boundary behaviour of an element in the enlarged
$b$-calculus can be analyzed by looking separately at the
$b$-boundary $x=0=x'$ and the  boundary $x'=2$.
 Let ${\cal U}\cong [0,7/4)\times \pa M$ 
be a collar neighbourhhod
of $\pa M$. 
We shall consider the restriction
of our operators to $\ub\Omega^*({\cal U};\Fl|_{{\cal U}})\oplus
\ub\Omega^*([0,7/4))\widehat{\otimes}\whW^*$
and to $\Omega^*([3/2,2))\widehat{\otimes}\whW^*$ separately. 
In this section,  
we consider a ``conic 
metric'' at $x'=2$ on  $\Omega^*([3/2,2))\widehat{\otimes}\whW^*$ and we put
a $b-$metric at $x'=0$. That is, we
will use a formula as in Section 7 to define
$\tau_{{\rm alg}}$, but now with a function $\chi \in C^\infty(0,2)$ such that
$\chi(x) = x$ for $x \in  (0, 1/2]$ and $\chi(x) = t$ for $x \in 
(3/2, 2)$.

We shall construct a
parametrix by patching  a
$b$-boundary parametrix and a parametrix for the conic-signature operator
on $\Omega^*([3/2,2))\widehat{\otimes}\whW^*$.
We pass to the implementation of this program,
concentrating first and foremost on the $b$-parametrix near the 
$b$-boundary.

First, notice that  elements
in the enlarged $b$-calculus form an algebra. The proof of the
appropriate composition formulae proceeds as in [Me].
Next, we recall that each  $P\in\widehat{\Psi}^m_{b,C^*_r(\Gamma)}$
has a well defined indicial family
$$I(P,\lambda): \Omega^*(\pa M;\Fl_0)\oplus
\whW\rightarrow \Omega^*(\pa M;\Fl_0)\oplus
\whW,$$
where we have implicitly used suitable identifications in a neighbourhood
of the $b$-boundary.
If $P_{11}$ and $P_{22}$ are $b$-elliptic in the usual sense
(i.e. symbolically)
and if $I(P,\lambda)$ is uniformly invertible for each $\lambda\in\RR$
then, by inverse Mellin transform, we can construct a parametrix
$G\in\widehat{\Psi}^{-m,\beta}_{b,C^*_r(\Gamma)}$, $\beta >0$, with remainders 
$R_1, R_2 \in \widehat{\Psi}^{-\infty,\beta}_{b,C^*_r(\Gamma)}$ 
having vanishing indicial family or, equivalently,  a vanishing 
restriction to the front face. 
(For this construction see [Me, Prop. 5.28] and [LP1, Theorem 11.1].) 

These remainders define bounded maps between  $C^*_r(\Gamma)$-Hilbert
modules, from
$$L^2_b(M;\ub\Lambda^* M\otimes\Fl)\oplus L^2_{b,{\rm comp}}
([0,7/4);\ub\Lambda^* [0,7/4)\widehat{\otimes} \whW^*)$$ 
to  
$$x^\beta H^1_b(M;\ub\Lambda^* M\otimes\Fl)\oplus (x')^\beta H^1_{b,{\rm loc}}
([0,7/4);\ub\Lambda^* [0,7/4)\widehat{\otimes} \whW^*).$$
(Recall that in order to have a {\it compact} Sobolev embedding in 
the framework of $b$-Sobolev spaces it is necessary to have 
a gain both in the order of the Sobolev space and in the
weighting. See [LP1, Lemma 11.2].
For the definition of the $C^*_r(\Gamma)$-Hermitian
scalar product on $L^2_b(M,\ub\Lambda^* M\otimes\Fl)$ we refer
to [LP1, Sect. 11]).

Going back to our perturbed $b$-signature operator $\SCb$, we can compute
its indicial family as follows. First, using 
Lemma \thm[ident.2] and a harmless abuse of notation,
we can fix the identifications
$$\Psi^+\equiv\Phi^{-1}: \left(\ub\Omega^*(M;\Fl)|_{\pa M}\oplus 
(\ub\Omega^*[0,2)\widehat{\otimes}
\whW^*)|_{x'=0}\right)^+ \rightarrow \Omega^*(\pa M;\Fl_0)\oplus \whW^*
$$
$$\Psi^-
\equiv\Phi^{-1}\circ \Theta: \left(\ub\Omega^*(M;\Fl)|_{\pa M}\oplus 
(\ub\Omega^*[0,2)\widehat{\otimes}
\whW^*)|_{x'=0}\right)^- \rightarrow \Omega^*(\pa M;\Fl_0)\oplus \whW^*\;.
$$
We thus obtain an isomorphism 
$$\Psi=\Psi^+\oplus\Psi^-: \ub\Omega^*(M;\Fl)|_{\pa M}\oplus 
(\ub\Omega^*[0,2)\widehat{\otimes}
\whW^*)|_{x'=0}\rightarrow (\Omega(\pa M;\Fl_0)\oplus \whW)\otimes\CC^2\;.$$
Using {\it this} isomorphism,
the matrices
$$\gamma=\left( \matrix{0 & i \cr -i & 0} \right)
\quad\quad
\sigma=\left( \matrix{0 & 1 \cr 1 & 0} \right),
$$
and proceeding as in [MP1, Sect. 1],
we obtain the indicial
family $I(\SCb,\lambda): (\Omega(\pa M;\Fl_0)\oplus \whW)\otimes\CC^2
\rightarrow (\Omega(\pa M;\Fl_0)\oplus \whW)\otimes\CC^2$ to be
$$I(\SCb,\lambda)=\left( \matrix{ \gamma\lambda+\sigma 
{\cal D}_{{\rm sign},\pa M} & 0 \cr
0 & \gamma\lambda+\sigma\Dw  } \right)
+(-i)\cdot 
\cases{\left(\matrix{0 & \sigma\widehat{\rho_\delta}(\lambda)
\tau_{\pa M}
\widehat{g} \cr
-\sigma\widehat{\rho_\delta}(\lambda) \widehat{f}\tau_{\pa M} & 0 } 
\right)&if$\;j<m$\cr \left(\matrix{0 & 
\sigma\widehat{\rho_\delta}(\lambda) \widehat{g}\tau_{\widehat{W}}\cr
-\sigma\widehat{\rho_\delta}(\lambda) \tau_{\widehat{W}}\widehat{f} &
0 } 
\right)&if$\;j>m$\cr}.
\Eq(indicial)
$$
Following [MP1, Lemma 9], we shall
now show that $I(\SCb,\lambda)$ is invertible for all $\lambda\in\RR$.
Since $\widehat{\rho_\delta}(\lambda)\in\RR$ for $\lambda\in\RR$
we see that $I(\SCb,\lambda)=\gamma\lambda + A(\lambda)$
with $A(\lambda)$ self-adjoint. From the definition of $\gamma$ 
we only need to check the invertibility of $I(\SCb,\lambda)$
at $\lambda=0$. Since $\widehat{\rho_\delta}(0)=1$, we 
obtain immediately that $I(\SCb,0)=\SC(1)_0$ and the invertibility
thus follows
from the very definition of the perturbed differential complex given
in Section 2.
In summary, the perturbed
signature operator 
$\SCb\in\widehat{\Psi}^1_{b,C^*_r(\Gamma)}$ has an invertible indicial
family for $\lambda\in\RR$.

\noindent
We can therefore apply the above $b$-parametrix construction to $P=\SCb$,
obtaining a   $G\in\widehat{\Psi}^{-m,\beta}_{b,C^*_r(\Gamma)}$.
If we now patch this $b$-parametrix $G$ with
a  parametrix 
for the signature operator on $\Omega^*([3/2,2))\widehat{\otimes}
\whW^*$,  we obtain 
a parametrix $\widetilde{G}$ for $\SCb$ 
with $C^*_r(\Gamma)$-compact remainders.
We omit the standard details.

Thanks to the work of Mischenko and Fomenko, we infer that the operator
$\SCb$ has a well defined index class in $K_0(C^*_r(\Gamma))$.
Proceeding as in [L1, Section VI] and [LP3, Appendix], 
this index class can be sharpened into
a $K_0(\Bi)$-class, using an appropriate $\Bi$-$b$-Mischenko-Fomenko
decomposition theorem 
[LP1, Sect. 15]. 

\nproclaim Proposition [wd]. The index class $\Ind(\SCbp)
\in K_0(\Bi)$ only depends
on the signature operator on $M$ and not on the choice of the 
finitely generated
Hermitian complex $W^*$, the homotopy equivalence $f$ or the functions 
$\rho_\delta$ or $\phi$.

\smallskip
\noindent
{\it Proof.} The independence of the choices of $\rho_\delta$ and $\phi$
is proved as in [LP2, Prop. 6.4]. The fact that different choices of $W^*$
and 
$f$ do not affect the index class is proved using
the idea of the proof of [L5, Proposition 15]. \QED

\medskip
\nproclaim Definition [canon]. We shall call the higher index class
$\Ind(\SCbp)\in K_0(\Bi)$ constructed above the 
$b$-signature-index class associated to  a manifold-with-boundary
$M$ satisfying Assumption 1.

The $b$-signature-index class
depends neither on the choice of Riemannian
metric nor on the choice of the Hodge duality
operator $\tau_M$, 
as different choices give
operators that can be connected by a suitable $1$-parameter family of
operators.
In  Section 11  we shall
compute the Chern character of $\Ind(\SCbp)\in K_0(\Bi)$, with values in the
noncommutative de Rham homology of $\Bi$, in terms of the usual local integral
and the higher eta invariant defined in Section 2.

%\fine
\expandafter\ifx\csname sezioniseparate\endcsname\relax%
   \input macro \fi

%--------------------------------------------------------
\numsec=10
\numfor=1
\numtheo=1
\pgn=1

\beginsection 10. Equality of the APS  and 
$b$-index classes

We recall Remark \thm[remark] and that we assume $t_1=1$ in the definition
of the index class $\Ind(\SCAPSp)$.

\nproclaim Theorem [aps=b]. The following equality holds in
$K_0(C^*_r(\Gamma))$: $\Ind(\SCAPSp)= 
\Ind(\SCbp)$.

\smallskip
\noindent
{\it Proof.} The signature operator $ \SCp$, associated to
the odd operator $\SC(1)$ of Section 4, 
 induces an operator,
denoted $\SCcp$,
 acting on the same  $C^*_r(\Gamma)-$Hilbert modules
 as  $\SCbp$. This amounts to adding a cylindrical end
to the manifold with boundary $M$ of Section 4
and a half-line to the interval $[0,2)$.   The
extended operator $\SCcp$  is {\it not}
in the $b$-calculus; the role of the function $\rho_\delta$
in Section 8 was precisely that of providing  a perturbation
belonging to  the $b$-calculus. This will be crucial in order to 
prove the higher APS-index {\it formula} in Section 11. However
as far as {\it Fredholm properties} are concerned, the operator 
$\SCcp$, i.e. the operator 
$\SCp$  of Section 2 extended with cut-off functions
to the cylindrical
ends, can be proven directly to be $C^*_r(\Gamma)$-Fredholm. 
In order to show this fact we shall still employ
ideas from [Me, Sect 5.4, Sect. 5.5].
The Schwartz kernel of the perturbation in 
$\SCc$, constructed from 
$f$ and $g$ in Section 1, lifts to a distribution
on the $b$-stretched product which will be smooth outside
the $b$-diagonal
and vanishing to infinite
order at the left and right boundaries. 
The operator
$\SCc$ still 
admits an indicial family: 
$$I(\SCc,\lambda)=\left( \matrix{ \gamma\lambda+\sigma 
{\cal D}_{{\rm sign}}(\pa M) & 0 \cr
0 & \gamma\lambda+\sigma\Dw  } \right)
+(-i)
\cases{\left(\matrix{0 & \sigma\tau_{\pa M}\widehat{g} \cr
-\sigma\widehat{f}\tau_{\pa M}  & 0 } 
\right)&if$\;j<m$\cr \left(\matrix{0 & 
\sigma\widehat{g}\tau_{\whW}\cr
-\sigma\tau_{\whW}\widehat{f} & 0 } \right)&if$\;j>m$\cr}.
$$
Let us prove that there exists a bounded operator 
$G:H^{-1}_b\oplus H^{-1}_b \rightarrow L^2_b\oplus L^2_b$ and a positive
real number $s$  such that
$G\circ \SCc(1)-\id$ is bounded from
$L^2_b\oplus L^2_b$ into $x^s H^1_b\oplus (x')^s H^1_b $;
the two $L^2$ spaces here refer to $M$
and $[0,2)$, respectively.
$G$ will then be an inverse modulo $C^*_r(\Gamma)$-compacts; this will prove
the $C^*_r(\Gamma)$-Fredholm-property
for $\SCc$.
$G$ will be obtained as in Section 5
by patching a ``$b$-parametrix'' and
a parametrix for the (conic-)signature operator on 
$\Omega([3/2,2))\widehat{\otimes}\whW^*$. We shall only 
concentrate on the ``$b$-parametrix''.
Working symbolically first, we can find a $G_\sigma\in
\widehat{\Psi}^{-1}_{b,\Bi}$ such that
$G_\sigma \circ \SCc= \id + R$, with
$R$ sending $H^{p-1}_b\oplus H^{p-1}_b$ into $H^p_b\oplus H^p_b$ 
for any $p\in\ZZ$
and with a Schwartz kernel which lifts to the stretched
product as a distribution, smooth outside $\Delta_b$
and vanishing to infinite order at the left and right boundaries.

If $N>>1$ then there exists $A>>1$ such that
$$(\sum_{k\geq 0}^A (-1)^k R^k )\circ G_\sigma \circ \SCc
=
\id + K_A,$$
with the Schwartz kernel of $K_A$ lifting to a $C^N$-function
on the $b$-stretched product, smooth near 
${\rm lb}$ and ${\rm rb}$  and vanishing to infinite order
there.
If $N$ is big enough then, proceeding as in  [Me] (see 
Section 5.13),
we can find $T:H^{-1}_b\oplus H^{-1}_b \rightarrow L^2_b\oplus L^2_b$
such that the lift of the  Schwartz kernel of $T$
is $C^k$ in the interior, has conormal bounds 
at the right and left boundaries 
and satisfies
$$I(T,\lambda)\circ I(\SCc,\lambda)= -I(K_A,\lambda)$$
(we simply take the inverse Mellin
transform of $-I(K_A,\lambda)\circ (I(\SCc),\lambda)^{-1}$).
The operator $G=T+ \sum_{k\geq 0}^A (-1)^k R^k \circ G_\sigma$
provides a left $b$-parametrix with an 
(adjointable)
remainder which 
continuously maps $L^2_b\oplus L^2_b$ into $x^s H^1_b\oplus (x')^s H^1_b $
for 
a suitable $s > 0$.

Given $t\in [0,1]$,
each     $ F(t)=t\SCcp  + (1-t)\SCbp$ has 
 an invertible indicial family. Once again, 
$F(t)$ is not a $b-$pseudo-differential operator 
but its Schwartz kernel lifted to the $b-$stretched product vanishes to 
infinite order at ${\rm lb}$ and ${\rm rb}$.  
 One can then  construct  as above
a parametrix $G$ sending the 
$C^*_r(\Gamma)-$Hibert module $H^{-1}_b\oplus H^{-1}_b$ into
$L^2_b\oplus L^2_b$ such that $F(t)\circ G- \id$ and $G \circ F(t) - \id$ are 
$C^*_r(\Gamma)-$compact. 
Since the family $\{F(t)\}_{0\leq t \leq 1}$ is obviously 
continuous, we have $\Ind(\SCbp) = \Ind(\SCcp)$. 
Now we can find a finite number 
of elements $u_1,\ldots, u_N$ of $(\,\ub\Omega^*(M,\Fli) \oplus 
\ub\Omega^*([0,2))\widehat{\otimes} 
\whW^*\,)^-$ which vanish in a neighborhood of the $b$-boundary 
($x=0,\, x'=0$) (i.e. at $-\infty$ on the cylinders)
 such that if one denotes by $K$ the operator defined 
by $K(s_1,\ldots, s_N)= \sum_{j=1}^N s_j u_j$ 
for any 
$(s_1, \dots, s_N) \in (C^*_r(\Gamma))^N$ then the operator 
$\SCcp \oplus K$ is 
surjective from $(H^1_b)^+ \oplus (C^*_r(\Gamma))^N$ to $(L^2_b)^-$. 
Then, as is well 
known, $ \Ind(\SCcp)= 
[\Ker (\SCcp \oplus K ) - (\Bi)^N]$. Of course, we can assume 
that $u_1,\ldots,u_N$ are such that 
 $$\Ind(\SCAPSp)= 
 [\Ker_{{\rm APS}} (\SCp \oplus K ) - (\Bi)^N]
$$
Since the identification $\Ker_{{\rm APS}} (\SCp \oplus K )
\equiv \Ker (\SCcp \oplus K ) $
follows by standard arguments, the theorem is proved.

%\fine
\expandafter\ifx\csname sezioniseparate\endcsname\relax%
   \input macro \fi

%------------------------

\numsec=11
\numfor=1
\numtheo=1
\pgn=1

\beginsection 11. The higher index formula for the $b$-signature operator

We want to adapt the proof of the higher APS-index formula given
in [LP2] to the present situation. First we need to comment about the
existence of a heat kernel
for the perturbed signature Laplacian 
$(\SCb)^2$.
We shall concentrate on the $b$-boundary, 
since the heat kernel near $x'=2$ is a consequence of [Ch].
We  can write $(\SCb)^2$
as
$$(\SCb)^2=\Delta + P
\;\;{\rm with}\;\;P=\left(\matrix{ P_{11} & P_{12} \cr P_{21} & P_{22} }
\right),$$
where  $P$ is a smoothing operator in the 
enlarged $\Bi$-$b$-calculus. Moreover, on 
$\ub\Omega^*({\cal { U}};\, {\cal {V}}|_{{\cal {U}}} )\oplus 
(\ub\Omega^*([0,1/2]) \otimes
\whW^*)$, we have 
$$\Delta=\left( \matrix{ D_M-\tau_M D_M \tau_M & 0 \cr 0 & 
D_{{\rm alg}} - \tau_{{\rm alg}} D_{{\rm alg}} \tau_{{\rm alg}} } 
\right)^2 = \left( \matrix{ \Delta_M & 0 \cr 0 & \Delta_{[0,2)}\otimes \id +
\id \otimes \Delta_{\whW} } \right)\;.
$$

The heat kernel of $\Delta$ is certainly well-defined
as an element of a $\Bi$-$b$-heat calculus; see [Me, Sect. 7] and 
[LP1, Sect. 10]. Using exactly the
same technique as in the proof of [MP1, Proposition 8], we can construct the
heat kernel of $(\SCb)^2$ as follows. We set $H^{(0)}=\exp(-t\Delta)$ 
and consider
$$({d\over dt} + (\SCb)^2) H^{(0)}=R^{(0)},\quad\quad R^{(0)}=
 P \exp(-t\Delta).$$
Using the indicial family of $(\SCb)^2$ and the heat equation on the boundary,
we can inductively remove the whole Taylor series of $R^{(0)}$
at the front faces where it is not zero and thus define 
an $H$ in an enlarged   $\Bi$-$b$-heat calculus such that
$$({d\over dt} + (\SCb)^2) H=R
\quad {\rm with} \quad R\in 
C^\infty((0,\infty)_t;\rho_{{\rm bf}}^\infty
\cdot\widehat{\Psi}^{-\infty}_{b,\Bi}).$$
At this point, the heat kernel $\exp(-t(\SCb)^2)$
is obtained by summing the usual Duhamel's series:
$$e^{-t(\SCb)^2}=H+\sum_{k>1} \int_{t S^k} H(t-t_k) R(t_k-t_{k-1})
\cdots R(t_1)dt_k \ldots dt_1
\quad\in \widehat{\Psi}^{-\infty}_{b,\Bi}\quad\forall t>0,$$
with $t S^k=\{(t_1,t_2,\dots,t_k); 0\leq t_1
\leq t_2 \cdots \leq t_k\leq t\}$.
We refer the reader to [MP1, Proposition 8] for the details.

Next, we need to introduce a superconnection
$\AAA_s$ 
as in [L1, (51)] and [LP1], 
and define the associated superconnection
heat kernel.
We fix, as in \equ(sdconnect), a $\tau_{W}$-invariant
connection 
$$\nabla ^{\whW}: \whW^*\rightarrow\Omega_1(\Bi)\otimes_{\Bi}
\whW^*\, ,$$ 
define the superconnection
$$s(D_{{\rm alg}} - \tau_{{\rm alg}} D_{{\rm alg}} \tau_{{\rm alg}})
+ \nabla^{\whW}:
C^\infty([0,2);\ub\Lambda^* [0,2))
\widehat{\otimes} \whW^*\rightarrow C^\infty([0,2);\ub\Lambda^* [0,2))
\widehat{\otimes} (\Omega_*(\Bi) \otimes_{\Bi} \whW^*)$$
and consider the total superconnection
$$\AAA_s=\left( \matrix{s(D_M-\tau_M D_M \tau_M)+\nabla & 0 \cr
                        0 & s(D_{{\rm alg}} - \tau_{{\rm alg}} D_{{\rm alg}} \tau_{{\rm alg}})+\nabla^{\whW} }
          \right) :
C^*\rightarrow C^*\oplus (\Omega_1(\Bi)\otimes_{\Bi} C^*),
$$
with 
$$\nabla: C^\infty(M;\ub\Lambda^* M\otimes\Fli)
\rightarrow 
\Omega_1(\Bi)\otimes_{\Bi} C^\infty(M;\ub\Lambda^* M\otimes \Fli) \Eq(ref)$$ 
as in [L2, Proposition 9].
$\AAA_s$ extends to a map $ \Omega_*(\Bi) \otimes_{\Bi} C^*
\rightarrow \Omega_*(\Bi)\otimes_{\Bi} C^*$ which is odd
with respect to the total $\ZZ_2$-grading and satisfies
Leibniz' rule.

This is the unperturbed superconnection; we shall need
$$\pAAA_s=s\SCbs + \left( \matrix{\nabla & 0 \cr
                                  0 & \nabla^{\whW}}
\right)
$$
with 
$$\SCbs=\left(\matrix{D_M-\tau_M D_M \tau_M & 0 \cr
                     0 & D_{{\rm alg}} - \tau_{{\rm alg}} D_{{\rm alg}} \tau_{{\rm alg}} }\right)
+\eps(s)\left(\matrix{0 & S(\delta) \cr
                       T(\delta) & 0} \right)$$
where we have used the notation at  the end of the previous section
%\cases{\left(\matrix{0 & g_\alpha \cr
%\tau_{(0,2)}\otimes \tau_{\whW} \circ f_\alpha \circ \tau_B & 0 } 
%\right)&if$\;i<m$\cr \left(\matrix{0 & \tau_B \circ g_\alpha
%\circ \tau_{(0,2)}\otimes \tau_{\whW}\cr
%f_\alpha & 0 } \right)&if$\;i>m$\cr}.
%$$
and $\eps\in C^\infty (0,\infty)$ is a nondecreasing function such that
$\eps(s)=0$ for $s \in (0,2)$ and 
$\eps(s)=1$
for $s>4$.

\smallskip
\noindent
{\bf Remark.} The operator $\SCbs$, and thus the superconnection
$\pAAA_s$, depend on $\delta$.
Using Duhamel's expansion and the existence of the heat kernel 
$\exp(-(s\SCbs)^2)$ we can define the superconnection heat kernel
$\exp(-\pAAA^2_s)$. For each $s>0$, it is a smoothing operator
in the enlarged $b$-calculus with coefficients in
$\Omega_*(\Bi)$. We denote the latter space by
$\widehat{\Psi}_{b,\Omega_*(\Bi)}^{-\infty}$; thus
$\exp(-\pAAA^2_s)\in\widehat{\Psi}_{b,\Omega_*(\Bi)}^{-\infty}$. Then
$$\exp(-\pAAA_s^2)=\left(\matrix{E_{11} & E_{12} \cr
                                E_{21} & E_{22}} \right),$$
 with
$$ E_{11}\in
\Psi^{-\infty}_{b, \Omega_*(\Bi)}(M;
\Omega_*(\Bi)\otimes_{\Bi}\ub\Lambda^*M \otimes \Fli)
,\quad E_{22}\in
\Psi^{-\infty}_{b, \Omega_*(\Bi)}([0,2);\Omega_*(\Bi)
\otimes_{\Bi}
\ub\Lambda^* [0,2)\widehat{\otimes}
\whW^*) 
.$$ 
Each operator such as $E_{11}$ has a well defined
$b$-supertrace, with values in $\overline{\Omega}_*(\Bi)$; see [LP1, Section 
13].
Similarly each element such as $E_{22}$ will have a well-defined
$b$-supertrace, also with values in $\overline{\Omega}_*(\Bi)$.
Notice that the heat kernel in a neighbourhood
of $x'=2$ has a well-defined supertrace - 
there is no need for regularization
there.

We define the $b$-supertrace of $\exp(-\pAAA_s^2)$ as
$$\bSTR(\exp(-\AAA_s^2))=\bSTR(E_{11})+ \bSTR(E_{22}) \in
\overline{\Omega}_*(\Bi).$$
The same definition applies to any element
$$R=\left(\matrix{R_{11} & R_{12} \cr
                  R_{21} & R_{22} }\right) \in 
\widehat{\Psi}^{-\infty}_{b,\Omega_*(\Bi)}.$$

This $b$-supertrace is not necessarily zero on supercommutators. As in [Me]
(and then [MP1], [LP1]), one can write a formula for the $b$-supertrace
of a supercommutator of two elements 
$R, S \in 
\widehat{\Psi}^{-\infty}_{b,\Omega_*(\Bi)}$.

\nproclaim
Proposition [defect]. Given
$$R=\left(\matrix{R_{11} & R_{12} \cr
                  R_{21} & R_{22} }\right)\;,\quad
S=\left(\matrix{S_{11} & S_{12} \cr
                  S_{21} & S_{22} }\right) \in 
\widehat{\Psi}^{-\infty}_{b,\Omega_*(\Bi)}\;,$$
the following formula holds in $\overline{\Omega}_*(\Bi)$ :
$$\bSTR [R,S]=
{\sqrt{-1}\over 2\pi}\int_{-\infty}^{+\infty}
{\rm STR}({\partial I(R,\lambda) \over \partial\lambda}\circ 
I(S,\lambda)) d\lambda.
$$
Moreover, the same formula holds if 
$R_{11}$ and $R_{22}$ are $b$-differential.

\smallskip
\noindent
{\it Proof.} On applying some straightforward
linear algebra, the proof can be eventually reduced to the one
in [Me, Prop. 5.9]; the details are exactly
as in  [MP1, Prop. 9], [LP1, Sect. 13]
but with the additional (harmless) complication coming from the 
fact that we are dealing with the {\it enlarged} $b$-calculus.
Since the details are elementary but tedious, we omit them. \QED
\medskip

 Using the $b$-supercommutator formula 
and proceeding as in [MP1, Proposition 11], 
we can  now compute the $s$-derivative of
$\bSTR(\exp(-\pAAA_s^2)\,)$. To this  end,
we first need to analyze the boundary behaviour
of $\pAAA_s$.

Let $\nabla^{\pa M}: C^\infty(\pa M; \Lambda^* (\pa M)\otimes
\Fli_0\otimes\Cl(1))\rightarrow \Omega_1(\Bi) 
\otimes_{\Bi} C^\infty(\pa M; \Lambda^* (\pa M)\otimes
\Fli_0\otimes\Cl(1))
$ be  the
$\Cl(1)$ analog of the connection in \equ(ref). 
We consider
 $$\pBBB_s(\lambda)=s \left(\matrix{
\sigma {\cal D}_{{\rm sign}}(\pa M) & 0 \cr 0 & \sigma\Dw}
\right) +
s\eps(s)\left(\matrix{0 & I(S(\delta),\lambda) \cr
                      I(T(\delta),\lambda) & 0}\right)
+\left(\matrix{\nabla^{\pa B} & 0\cr
               0 & \nabla^{\whW}}\right)$$
where we recall, see \equ(indicial), that using our identifications at the
boundary, 
$$\left(\matrix{0 & I(S(\delta),\lambda) \cr
                      I(T(\delta),\lambda) & 0}\right)
=(-i)\cdot 
\cases{\left(\matrix{0 & \sigma\widehat{\rho_\delta}(\lambda)
                      \tau_{\pa M}
\widehat{g} \cr
-\sigma\widehat{\rho_\delta}(\lambda) \widehat{f}\tau_{\pa M} & 0 } 
\right)&if$\;j<m$\cr \left(\matrix{0 & 
\sigma\widehat{\rho_\delta}(\lambda) \widehat{g}\tau_{\widehat{W}}\cr
-\sigma\widehat{\rho_\delta}(\lambda) \tau_{\widehat{W}}\widehat{f} & 0 } \right)&if$\;j>m$\cr}.
\Eq(indicial2)$$
For each fixed $\lambda\in\RR,$ $\pBBB_s(\lambda)$
is a $\Cl(1)$-superconnection, mapping 
 $$(C^\infty(\pa M;\Lambda^*(\pa M)
\otimes \Fli_0)\oplus \whW^*) \otimes\Cl(1) \rightarrow
\Omega_*(\Bi) \otimes_{\Bi}
(C^\infty(\pa M;\Lambda^*(\pa M)
\otimes \Fli_0)\oplus \whW^*) \otimes\Cl(1).$$
It depends on $\delta$  through the indicial families
 \equ(indicial2). 

\nproclaim Proposition [tran].
The following formula holds in $\overline{\Omega}_*(\Bi)$:
$${d\over ds}\left(b-{\rm STR} \exp(-\pAAA_s^2)\right)=
-d\left({\rm b-STR}\left({d\pAAA_s\over ds}e^{-\pAAA_s^2}
\right)\right)- \widehat{\eta}_\delta (s)
\Eq(tranf)
$$
with
$$\widehat{\eta}_\delta (s)={i\over 2 \pi}\int_{\RR}
{\rm STR} \left( s\left(\matrix{\gamma & 0 \cr 0 & \gamma}\right)\cdot
{d\pBBB_s(\lambda)\over ds}e^{-\pBBB^2_s} \right) e^{-(s\lambda)^2} d\lambda$$
$$+{i\over 2 \pi}\int_{\RR} {\rm STR} \left( s\eps(s)\left(\matrix{\gamma & 0
\cr 0 & \gamma}\right)\lambda \left({d\over d\lambda} 
\left(\matrix{0 & I(S(\delta),\lambda) \cr
                      I(T(\delta),\lambda) & 0}\right)\right)
\cdot e^{-\pBBB^2_s}\right) e^{-(s\lambda)^2} d\lambda$$

\smallskip
\noindent
{\it Proof.} Using Proposition 
\thm[defect], the proof given in [MP1, Prop. 11], [LP2, Prop. 7.2]
can be easily adapted to the present situation. \QED

\smallskip
We define 
$\widetilde{\eta}_\delta (s) =
{\cal R}\, \widehat{\eta}_\delta (s)$,
with ${\cal R}$ as in Section 2, and we put
$\widetilde{\eta}(\delta)=
\int_0^\infty \widetilde{\eta}_\delta (s)ds$.
The convergence of the integral can be proven as
in Proposition 7.4 in [LP2].
Regarding the relationship of $\widetilde{\eta}(\delta)$
to the $\widetilde{\eta}_{\pa M}$ of Section 2 (see
equation \equ(heta)
with $F=\pa M$), we have the following proposition.

\nproclaim Proposition [equal]. The following
equality holds  for any $\delta > 0$:
$$\widetilde{\eta}(\delta)\, = 
\widetilde{\eta}_{\pa M}\;\;{\rm in}\;\;
\overline{\Omega}_*(\Bi)/d
\overline{\Omega}_*(\Bi). \Eq(equalf)$$

\smallskip
\noindent
{\it Proof.} The proof of Proposition 7.4 (2) in  [LP2]
applies {\it mutatis mutandis} to the more general case treated
here. \QED

\smallskip

These two propositions
are
crucial in the
proof of the following higher  index theorem for the $b$-signature operator :

\nproclaim Theorem [higher]. Let $M$ be an even-dimensional 
manifold-with-boundary. Let $\Gamma$ be a finitely-generated discrete group
and let $\nu : M \rightarrow B\Gamma$ be a continuous map, defined up to
homotopy.
We make  {\bf Assumption 1}
on $\pa M$. Let $g$ be a $b$-metric on $M$ which is 
product-like near the boundary and let $R^M$ be
the associated curvature $2$-form.  The following
formula holds for the Chern character of
the canonical higher $b$-signature-index class:
$$\ch(\Ind(\SCbp))=\int_M L (R^M/2 \pi) \wedge \omega
- \, \widetilde{\eta}_{\pa M}\;\;{\rm in}\;\;\overline{\HH}_*(\Bi)
\Eq(higherf)$$
with $\widetilde{\eta}_{\pa M}$ equal to the higher eta-invariant
of Section 2.

\smallskip
\noindent
{\it Proof.}
Integrating formula \equ(tranf), we obtain that for $u>t>0$
and modulo exact forms,
$${\cal R} \, \bSTR(e^{-\pAAA^2_u})= {\cal R} \, \bSTR(e^{-\pAAA_t^2})-
\int_t^u \widetilde{\eta}_\delta (s)ds\,.\Eq(transf2)$$
The limit of the right-hand-side as $t\rightarrow 0^+$ can be computed
as in [L1, Proposition 12] and [LP1] (we also use [Ch, Section 4] 
near $x'=2$).
Let us  consider the asymptotic expansion near $t=0$
of the first summand on the right-hand-side of \equ(transf2).
Using
[L1, Proposition 12], [LP1] and [Ch, Section 4], one sees that the coefficient
of $t^0$ will be the sum of three noncommutative differential
forms. 
The first term  is 
$$\int_M 
L(R^M /2 \pi) \wedge \omega \; .
$$

The second term is the $t^0$ term of the
$b$-integral over $[0,2)$ 
of the pointwise supertrace of the
heat kernel associated
to the superconnection
$$t(D_{{\rm alg}} - \tau_{{\rm alg}} D_{{\rm alg}} \tau_{{\rm alg}})
+ \nabla^{\whW}\;.$$
As we are effectively computing a heat kernel on the real line, the
local heat trace 
asymptotics will be of the form $t^{-1} \sum_{j=0}^\infty c_j t^{j}$.
(Note that the usual $t$ of the heat kernel expansion is, in our case,
$t^2$.)
Using the Duhamel formula, we see that the $t^0$ term is proportionate to 
the local supertrace of $[(D_{{\rm alg}} - 
\tau_{{\rm alg}} D_{{\rm alg}} \tau_{{\rm alg}}),
\nabla^{\whW}]$. As $\nabla^{\whW}$ is independent of $x$, this equals
${\rm Str}[{\cal D}^{\rm sign}_{\whW},\nabla^{\whW}]$. However, being
the supertrace of a supercommutator involving $\nabla^{\whW}$, this is exact 
as an element of
$\overline{\Omega}_*(\Bi)$. As we are working modulo exact forms,
it thus vanishes.

The third term is a eta-contribution coming from
$x'= 2$:
$$- \int_0^{\infty} \widetilde{\eta}_{\whW}(s)ds\;,$$
where we recall that 
$$\widetilde{\eta}_{\whW}(s) = 
{1 \over \sqrt{\pi}}\,{\cal R} \, 
{\rm STR}_{{\rm Cl}(1)} \sigma {\cal D}^{\rm sign}_{\whW}
\exp(-(s \sigma {\cal D}^{\rm sign}_{\whW} + \nabla^{\whW})^2 )\; .$$
Thus, modulo exact forms, 
we obtain:
$$\lim_{t\rightarrow 0^+} {\cal R} \, \bSTR(e^{-\pAAA_t^2})%-
%\int_t^u \widehat{\eta}_\delta (s)ds
%\right)
= \int_M 
L( R^M /2 \pi) \wedge \omega
- \int_0^{\infty} \widetilde{\eta}_{\whW}(s)ds, $$
%-\ha\int_0^u \widehat{\eta}_\delta (s) ds,$$
However, as in [L2, p. 227], a duality argument shows that for all $s > 0$,
modulo exact forms, 
$$\widetilde{\eta}_{\whW}(s)=0 \;.\Eq(wetavan)$$
Summarizing,
$${\cal R} \, \bSTR(e^{-\pAAA^2_u})= \int_M 
L( R^M/2 \pi )\wedge\omega
- \int_0^u \widetilde{\eta}_\delta (s)
ds\;.$$
Proceeding as in [LP1, Theorem 14.1]
and [LP3, Appendix], one can now complete the proof 
of the theorem.
We omit the details
as they are very similar to those
explained at length in the above references. 
(The Appendix of [LP3], which is 
based
on results in [L3, Section 6], extends
the higher APS-index theorem proven in [LP1, Section 14]
to any finitely-generated discrete group $\Gamma$, 
under a gap hypothesis for the boundary operator.) This proves the
theorem. 
\QED

%\fine
\expandafter\ifx\csname sezioniseparate\endcsname\relax%
   \input macro \fi

%--------------------------------------------------------
\numsec=12
\numfor=1
\numtheo=1
\pgn=1

\beginsection 12. Proofs of Theorem \thm[mainthm] and Corollaries
\thm[hominv]-\thm[cutandpaste]

Before proving the main results, let us comment about our normalization
of eta-invariants. In the case ${\cal B}^\infty = \CC$, our definition
of $\widetilde{\eta}_F$ in \equ(oddeta) 
gives {\it half} of the eta-invariant as defined in [APS, Theorem 3.10(iii)] 
for Dirac-type operators.
This normalization is more convenient for our
purposes, albeit
unconventional, and is also used in [BC]. In [APS, Theorem 4.14(iii)], 
the eta-invariant of
the (tangential) signature operator is defined in terms of an
operator on {\it even} forms, and so gives half of the APS eta-invariant
of the corresponding Dirac-type operator. The upshot is that when considering
the signature operator, our eta-invariant coincides with that considered
in [APS, Theorem 4.14(iii)].

\medskip
\noindent
{\bf Proof of Theorem \thm[mainthm] :} Suppose first that $M$ is
even-dimensional. By Theorem \thm[conichominv], the
conic index in $K_0(C^*_r(\Gamma))$ is an oriented-homotopy invariant.
By Theorems \thm[c=aps] and \thm[aps=b], 
the conic index equals the $b$-index. As 
${\cal B}^\infty$ is a dense subalgebra of $C^*_r(\Gamma)$ which is 
closed under the holomorphic functional calculus in $C^*_r(\Gamma)$, there
is an isomorphism 
$K_0({\cal B}^\infty) \cong K_0(C^*_r(\Gamma))$ [Co, Section IIIC]. By
Theorem \thm[higher], the Chern character of the index, as an element
of $\overline{\HH}_*({\cal B}^\infty)$, equals $\sigma_M$. The theorem
follows in this case.

If $M$ is odd-dimensional, say of dimension $n = 2m-1$, 
we can reduce to the even-dimensional case
by a standard trick, replacing $M$ by $M \times S^1$ and replacing $\Gamma$ by 
$\Gamma \times \ZZ$. Observe that by Fourier
transform, $C^*_r(\Gamma\times\ZZ)\cong
C^*_r(\Gamma)\otimes C^0 (S^1)$
and $\Bi_{\Gamma\times\ZZ^1}\cong \Bi_{\Gamma}\otimes C^{\infty}(S^1)$,
where $\otimes$ denotes a projective tensor
product.
Under these identifications, the signature operator of $M\times S^1$
can be identified with the suspension of the signature
operator on $M$, as in [Lu, p. 250] and [BF, p. 124].
%We take the smooth subalgebra of
%$C^*_r(\Gamma \times \ZZ)$ given by ${\cal B}^{\infty \prime} = 
%\Bi \otimes C^\infty(S^1)$.
Moreover, instead of the universal graded algebra 
$\Omega_*({\cal B}^{\infty}_{\Gamma\times\ZZ})$, 
it suffices
to deal 
with the smaller differential graded algebra 
$\Omega_*(\Bi_{\Gamma}) \, \widehat{\otimes} \,
\Omega^*(S^1)$. 
Let $\tau^{\prime}$ be the generator for $\HH^1(\ZZ; \ZZ) \subset
\HH^1(\ZZ; \CC)$. Put $\tau= \sqrt{-1}\,\tau ^{\prime}/2\pi
\in \HH^{1}(\ZZ; \CC)$. There is a natural desuspension map
$$
\langle \cdot\,,\,\tau \rangle:\overline{\HH}_*(\Bi_{\Gamma\times\ZZ})
\rightarrow \overline{\HH}_*(\Bi_{\Gamma})
$$
and one can check,
as in [MP2, Lemma 6], that
$$\langle \int_{M \times S^1} L(T(M \times S^1)) \wedge 
\omega_{\Gamma\times\ZZ}-\widetilde{\eta}_{\partial M \times S^1} \, , \,
\tau \rangle =
\int_{M} L(TM) \wedge \omega_{\Gamma} 
-
\widetilde{\eta}_{\partial M}
\;\;{\rm in}\;\;\overline{H}_*(\Bi_{\Gamma}).$$

In order to directly apply the even-dimensional results to $M \times S^1$, 
we would have to
know that if $\partial M$ satisfies Assumption 1 then 
$\partial M \times S^1$ satisfies
Assumption 1.  This is not quite true.  However, if we consider the
complex $\widehat{W}^*$ for $\partial M$, from \equ(newcomplexx),
and take its graded tensor product with 
$\Omega^*(S^1)$ then 
the terms in degrees $m-1$ and $m$ are $L \, \widehat{\otimes} \,
\Omega^1(S^1)$ and $L^\prime \, \widehat{\otimes} \,
\Omega^0(S^1)$, respectively,
with the differential between them being the zero map.
We can then go through all of the arguments in Sections 5-11 for $M \times 
S^1$, using
$\widehat{W}^* \, \widehat{\otimes} \,
\Omega^*(S^1)$ as the perturbing complex. 
The important point is that by duality, we again have
$\widehat{\eta}_{\widehat{W}} = 0$ [L2, p. 227].  
This implies Theorem 
\thm[mainthm] in the odd-dimensional case. \QED
\medskip
\noindent
{\bf Proof of Corollary \thm[hominv] :} 
This is an immediate consequence of
Theorem \thm[mainthm]. \QED

\medskip
\noindent
{\bf Example :}
Given $N > 0$ and $0 \le j \le 4N - 2$, put 
$M_1 = (\CC P^{2N} \# \CC P^{2N} \# D^{4N}) \times T^j$
and $M_2 = (\CC P^{2N} \# \overline{\CC P}^{2N} \# D^{4N}) \times T^j$. Then
$M_1$ and $M_2$ are homotopy equivalent as topological spaces, as they are
both homotopy equivalent to $(\CC P^{2N-1} \vee \CC P^{2N-1}) \times T^j$. 
Put $\Gamma = \ZZ^j$ and let $\nu_i : M_i \rightarrow
B\ZZ^j$ be the classifying maps for the universal covers.
Then Assumption 1 is satisfied.
Take $\tau = [B\ZZ^j] \in \HH^j(B\ZZ^j; \CC)$. 
Then $\langle \sigma_{M_1}, \tau \rangle = 
\sigma(\CC P^{2N} \# \CC P^{2N} \# D^{4N}) = 2$, while
$\langle \sigma_{M_2}, \tau \rangle = 
\sigma(\CC P^{2N} \# \overline{\CC P}^{2N} \# D^{4N}) = 0$. Thus by Corollary
\thm[hominv], $M_1$ and $M_2$ are not
homotopy equivalent as manifolds-with-boundary. 

\medskip
\noindent
{\bf Proof of Corollary \thm[novadd] :} From 
[L1, Corollary 2], the higher signature of
$M$ corresponding to $\tau \in \HH^*(\Gamma; \CC)$ is a nonzero constant
(which only depends on the degree of $\tau$) times
$\langle \int_M L(TM) \wedge \omega, Z_\tau>$. Let us choose a
Riemannian metric on $M$ so that a tubular neighborhood of $F$ is isometrically
a product. Then
$$\int_M L(TM) \wedge \omega = \int_A L(TM) \wedge \omega + 
\int_B L(TM) \wedge \omega =
\left( \int_A L(TM) \wedge \omega \right) - \widetilde{\eta}_{\partial A} + 
\left( \int_B L(TM) \wedge \omega \right) - \widetilde{\eta}_{\partial B},$$
as $\partial A$ and $\partial B$ differ in their orientations and
$\widetilde{\eta}$ is odd under a change of orientation.
Here we have
chosen a (stable) Lagrangian subspace $L$ for $\partial A$, if necessary,
and then taken the (stable) 
Lagrangian subspace $-L$ for $\partial B$.  
Thus 
$$\int_M L(TM) \wedge \omega = \sigma_A + \sigma_B. \Eq(topair)$$ 
The corollary follows from pairing both sides of \equ(topair) with 
the cyclic cocycle $Z_\tau$. \QED
\medskip
\noindent
{\bf Proof of Corollary \thm[cutandpaste] :} 
This is an immediate consequence
of Corollary \thm[novadd], along with the fact that $\sigma_A$ and
$\sigma_B$ are  
smooth topological invariants and only
depend on $\nu$ through its homotopy class. (If $\dim(M) = 2k+1$ then 
we use the 
(stable) Lagrangian subspace $L$ of $\overline{\HH}^k(F; {\cal V}_0)$,
assumed to be invariant under $(\phi_2 \circ \phi_1^{-1})^*$, to define
$\sigma_A$ and $-L$ to define $\sigma_B$.)
That $\sigma_A$ and $\sigma_B$ are 
smooth topological invariants follows from Theorem \thm[mainthm]. It also
follows more directly from
[L2, Proposition 27] and [L5, Theorem 6]. The papers [L2] and [L5] deal with
a slightly stronger
assumption than Assumption 1, but their proofs can be extended to the
present case, too. \QED

%\fine
\expandafter\ifx\csname sezioniseparate\endcsname\relax%
   \input macro \fi

%--------------------------------------------------------
\numsec=13
\numfor=1
\numtheo=1
\pgn=1

\beginsection 13. Appendix

Let $M$ be a compact oriented 
manifold-with-boundary and let $\nu:M\rightarrow B\Gamma$
be a continuous map. Let $M'$ be the associated
normal $\Gamma$-cover of $M$. If 
$\dim M=2m$ (resp. $\dim M=2m+1$) then, for the purposes of this
Appendix, we assume
that the differential form Laplacian has a strictly
positive spectrum on $\Omega^m(\partial M^\prime)$
(resp. $\Omega^m(\partial M^\prime)$).
This is a slightly stronger assumption
than Assumption 1; see Lemma \thm[Riem] and Lemma \thm[Riem2]. 

Under this assumption, a
higher $b$-signature-index class for manifolds with boundary 
was introduced in [LP4].
A higher signature
index formula was then proven in the virtually nilpotent case
using the higher APS index theorem proved in [LP1,2]. 
The regularization proposed there followed [LP2] and employed
the notion of a {\it symmetric} spectral section $\P$
for the boundary signature operator of $\cDi_{{\rm sign,b}}$.
The index class in [LP4] was denoted by $\Ind(\cDi_{{\rm sign,b}}^{\,+},\P)$
and was proven to be independent of the particular choice of symmetric
spectral section $\P$.
We shall now indicate how to prove that 
$\ch(\Ind(\cDi_{{\rm sign,b}}^{\,+},\P))=\sigma_M$, with $\sigma_M$
as \equ(highersig).
This will imply (for virtually nilpotent groups) that
the higher signatures considered in [L5]  and in [LP4]
are in fact the same. We shall only sketch the argument.

According to Theorem \thm[higher], it suffices to show
that 
$$\ch(\Ind(\cDi_C^{{\rm sign,b},+}))
=\ch(\Ind(\cDi_{{\rm sign,b}}^{\,+},\P)). \Eq(llp=lp)$$
Recall 
[LP2, Definition 6.3]
that $\Ind(\cDi_{{\rm sign,b}}^+,\P)$ is equal
to $\Ind((\cDi_{{\rm sign,b}} + A_{\P})^+)$, with $A_{\P}$ a regularizing
operator associated to $\P$. $A_{\P}$ is constructed as in Section 8
starting from a perturbation $A_{\P}^0$ on $\pa M$ which makes
the boundary operator $(\cDi_{{\rm sign,b}})_0$ invertible.
The symmetry of $\P$ corresponds to a vanishing
of $A_{\P}^0$ in middle degree plus a suitable $\ZZ_2$-grading
of $A_{\P}^0$ outside the middle degree, see [LP4, Definition 4.2].

We can extend the operator $\cDi_{{\rm sign,b}} + A_{\P}$
and make it act on 
$\ub\Omega^*(M;{\cal V}^\infty)\oplus
\ub\Omega^*[0,2)\widehat{\otimes}\whW^*$ without changing the
Chern character of the corresponding index class.
For this, it suffices to first consider
the operator 
$$\cDi^{\,\oplus}_{{\rm sign},b}=\left( \matrix{ {\cal D}_{{\rm sign},b}  & 0 \cr
0 & D_{{\rm alg}}-\tau_{{\rm alg}} D_{{\rm alg}} 
\tau_{{\rm alg}} }\right)                 \Eq(oplus) $$ 
and then the operator in the enlarged $b$-calculus given
by
$$\cDi^{\,\oplus}_{{\rm sign,b}}+\left( 
\matrix{ A_{\P} & 0 \cr 0 & 0 } \right)\equiv
\left( \matrix{ {\cal D}_{{\rm sign},b} + A_{\P} & 0 \cr
0 & D_{{\rm alg}}-\tau_{{\rm alg}} D_{{\rm alg}} 
\tau_{{\rm alg}} }\right). \Eq(oplus2)$$
Notice that the boundary operator of \equ(oplus2)
%, i.e.
%$$\left( \matrix{ {\cal D}^{\,0}_{{\rm sign},b}+A^0_{\P}  & 0 \cr
%0 & \cDi_{{\rm sign,\whW}}}\right) $$  
is not invertible; however we can add a perturbation
$A^0_{\whW}$
to $\cDi_{{\rm sign,\whW}}$ and make it invertible
on all of $\Omega^*(\pa M;{\cal V}^{\infty})\oplus \whW$.
The resulting operator is thus
$$\left( \matrix{ ({\cal D}_{{\rm sign},b})_0 + A^0_{\P}  & 0 \cr
0 & \cDi_{{\rm sign,\whW}}+A^0_{\whW}}\right)\,. \Eq(oplus2.5)$$ 
Let $A_{\whW}$ be the corresponding perturbation for 
$D_{{\rm alg}}-\tau_{{\rm alg}} D_{{\rm alg}} 
\tau_{{\rm alg}}$, constructed as in Section 8.
It is clear that the index class associated
to  the operator
$$\cDi^{\,\oplus}_{{\rm sign,b}}+\left( 
\matrix{ A_{\P}  & 0 \cr
0 & A_{\whW}}\right) 
\equiv
\left( \matrix{ {\cal D}_{{\rm sign},b} + A_{\P} & 0 \cr
0 & D_{{\rm alg}}-\tau_{{\rm alg}} D_{{\rm alg}} 
\tau_{{\rm alg}}+A_{\whW} }\right)\Eq(oplus3)$$
has the same Chern character as  
$\Ind(\cDi^{\, +}_{{\rm sign},b},\P)$.
We denote by $\cDi^{\,\oplus}_{{\rm sign},b}+A_{\P,\whW}$
the operator in \equ(oplus3). Thus
$$\ch(\Ind(\cDi^{\, +}_{{\rm sign},b},\P))=
\ch(\Ind(\cDi^{\,\oplus}_{{\rm sign},b}+A_{\P,\whW})^+)\,.\Eq(oplus4)$$
It remains to show that $\Ind((\cDi^{\,\oplus}_{{\rm sign},b}+A_{\P,\whW})^+)=
\Ind(\cDi_C^{{\rm sign,b},+})$ in $K_0(C^*_r(\Gamma))\otimes\QQ$.
To this end, we remark that the boundary operator of 
 $\cDi^{\,\oplus,+}_{{\rm sign,b}}$, denoted as usual by
$(\cDi^{\,\oplus}_{{\rm sign,b}})_0$, is {\it invertible}
in the middle degree of 
the Hermitian complex $\Omega(\pa M;{\cal V}^{\infty})\oplus \whW$.
The notion of spectral section and of symmetric spectral section
for $(\cDi^{\,\oplus}_{{\rm sign,b}})_0$
can be extended to the more general situation
considered in Sections 8 and 9. The APS-spectral
projections $\Pi_>$ for the boundary operator of \equ(oplus3) 
and for the boundary operator of  
$(\cDi_{{\rm sign,b}}^{C,+})$ in  Section 8
are both  examples of spectral sections  for
$(\cDi^{\,\oplus}_{{\rm sign,b}})_0$.
Here we have used the correspondence, explained in detail
in [MP1] and [Wu], between APS-spectral projections for
perturbed Dirac-type operators and spectral sections for
unperturbed Dirac-type operators. 
Let us denote by $\P^{\,\oplus}$
and $\Q^C$ these particular spectral sections.
It turns out that $\P^{\,\oplus}$
and $\Q^C$ are in fact {\it symmetric}
spectral sections; 
this follows from the structure of 
the above perturbation (see \equ(oplus2.5))
and of  the mapping-cone-perturbation of Section 4 (see \equ(psign)).
More precisely, it is implied by the vanishing in middle degree
plus a $\ZZ_2$-grading outside the middle degree for
$$\left( \matrix{ A^0_{\P} & 0 \cr
0 & A^0_{\whW}} \right)$$ 
and for 
$$\cases{\left(\matrix{0 & \epsilon\tau_{\pa M} \widehat{g} \cr
-\epsilon \widehat{f} \tau_{\pa M} & 0 } \right)&if$\;j < m- {1 \over 2}$\cr 
\left(\matrix{0 & \epsilon \widehat{g} \tw \cr
-\epsilon \tw \widehat{f} & 0 } \right)&if$\;j > m - {1 \over 2}$,\cr} 
$$
respectively.

Thus, by definition,
$$\Ind((\cDi^{\,\oplus}_{{\rm sign},b}+A_{\P,\whW})^+)=
\Ind(\cDi^{\,\oplus,+}_{{\rm sign},b},\P^{\,\oplus}),\quad\quad
\Ind(\cDi_C^{{\rm sign,b},+})=\Ind(\cDi^{\,\oplus,+}_{{\rm sign},b},\Q^C)
\,.\Eq(oplus6)
$$
 By the relative index theorem of [LP4, Prop. 6.2], suitably extended to this
more general setting,
we obtain
$$\Ind(\cDi^{\,\oplus,+}_{{\rm sign},b},\P^{\,\oplus})-
\Ind(\cDi^{\,\oplus,+}_{{\rm sign},b},\Q^C)=[\Q^C - \P^{\,\oplus}].
\Eq(oplus7)$$
However, both $\P^{\,\oplus}$ and $\Q^C$ are symmetric spectral sections.
Thus
by the symmetry argument in [LP4], see in particular
[LP4, Prop. 4.4], we have $[\Q^C - \P^{\,\oplus}]=0$
in $K_0(C^*_r(\Gamma))\otimes\QQ$. The claim then
follows from \equ(oplus4), \equ(oplus6) and \equ(oplus7). 
The odd dimensional case is similar, using ${\rm Cl(1)}$-symmetric
spectral sections.

This  proof shows that in the particular case $F=\pa M$
(and $\Gamma$ virtually nilpotent),
the two regularizations of the higher eta invariant
proposed in [L5, Definition 8] and [LP4, Definition 5.2] coincide. 
The regularization using symmetric spectral sections
can also be given for any closed oriented manifold $F$
and any normal $\Gamma$-cover $F^{\prime}$
satisfying the assumption that $\Delta_{F^{\prime}}$
is $L^2$-invertible in middle degree (but with $\Gamma$
still virtually nilpotent).
Using the 
above arguments and an extension of the jump-formula for 
higher eta-forms [MP1, Proposition 17], [LP2, Theorem 5.1], 
one can show that in
this general situation,  the   two 
definitions of the higher eta form given in [LP4, Definition 5.2] and 
[L5, Definition 8]
coincide.

%\fine
\expandafter\ifx\csname sezioniseparate\endcsname\relax%
   \input macro \fi

\beginsection References.

\frenchspacing
\item{[APS]} Atiyah M.F., Patodi V. and Singer I. 1975.  
Spectral asymmetry and Riemannian 
geometry I. {\it  Math. Proc. Cambridge Phil. Soc.}. {\bf 77},  p. 43-69
\item{[BaJ]} Baaj S. and Julg P. 1983. 
Th\'eorie bivariante de Kasparov et op\'erateurs non born\'es dans les
$C^*$-modules hilbertiens.
{\it C. R. Acad. Sc. Paris} {\bf 296}, p. 875-878  
\item{[BC]} Bismut J.-M. and Cheeger J. 1990. 
Families index for  manifolds with boundaries, superconnections and cones, 
I. {\it J. of Funct. Anal.} {\bf 89}, p. 313-363
\item{[BF]} Bismut J.-M. and Freed D. 1986. 
The analysis of elliptic families II.
Dirac operators, eta invariants and the holonomy theorem.
{\it Comm. Math. Phys.} {\bf 107}, p. 103-163  
\item{[BS]} Br\"uning J. and Seeley R. 1988. 
 An index theorem for first order regular
 singular operators. {\it  Amer. J. of Math.} {\bf 110} p. 659-714 
\item{[Bu]} Bunke U. 1995. 
On the gluing problem for the $\eta$-invariant.
{\it  J. Diff. Geom.} {\bf 41} p. 397-448 
\item{[Ch]} Cheeger J. 1983. Spectral geometry of singular Riemann spaces.
{\it J. Diff. Geom.} {\bf 18}, p. 575-657 
\item{[Co]} Connes A. 1994. {Noncommutative geometry}, 
Academic Press. San Diego
\item{[CM]} Connes A. and Moscovici H. 1990. Cyclic cohomology, the Novikov
Conjecture and hyperbolic groups. {\it Topology} {\bf 29}, p. 345-388
\item{[dlH]} de la Harpe. 1988. Groupes hyperboliques, alg\`ebres 
d'op\'erateur et un th\'eor\`eme de Jolissaint. {\it C. R. Acad. Sci. Paris}
S\'er. I Math. {\bf 327}, p. 771-774
\item{[D]} Dodziuk J. 1977. De Rham-Hodge theory for $L^2$-cohomology of
infinite coverings. {\it Topology} {\bf 16}, p. 157-165
\item{[Do]} Donnelly H. 1980.
The differential form spectrum of hyperbolic space. 
{\it Manuscripta Math} {\bf 33}, p. 365-385 
\item{[Ge]} Getzler E. 1993. The odd Chern character in cyclic homology and
spectral flow.
{\it Topology} {\bf 32}, p. 489-507 
\item{[H]} Hilsum M. 1989. Fonctoriali\'e en K-Th\'eorie bivariante
pour les vari\'et\'es lipschitziennes. 
{\it  K-Theory}. {\bf 3}, 
p. 401-440
\item{[HS]} Hilsum M. and Skandalis G.  1990. Invariance par homotopie de 
la signature  \`a coefficients dans un fibr\'e presque plat.
{\it  J. Reine Angew. Math}. {\bf 423}, 
p. 73-99
\item{[J]} Ji R. 1992.
Smooth dense subalgebras of reduced group $C^*$-algebras,
Schwartz cohomology of groups and cyclic cohomology.
{\it J. of Funct. Anal.} {\bf 107}, p. 1-33 
\item{[K]} Karoubi M. 1987. Homologie cyclique et K-th\'eorie, 
{\it Ast\'erisque}.
{\bf 149}
\item{[KKNO]} Karras U., Kreck M., Neumann W. and Ossa E. 1973. 
Cutting and pasting of manifolds; $SK$-groups. Publish or Perish. Boston
\item{[KM]} Kaminker J. and Miller J. 1985. Homotopy invariance of
the analytic index of signature operators over $C^*$-algebras.
{\it J. Operator Theory}. {\bf 14}, p. 113-127
\item{[Ks]} Kasparov G. 1980. Hilbert $C^*$-modules : theorems of
Stinespring and Voiculescu.
{\it J. Operator Theory}. {\bf 4}, p. 133-150
\item{[LP1]} Leichtnam E. and Piazza P. 1997.
The $b$-pseudodifferential calculus on Galois
coverings and a higher Atiyah-Patodi-Singer index
theorem.
{\it M\'emoires of the S. M. F.} {\bf 68} 
\item{[LP2]} Leichtnam E. and Piazza P. 1996. Spectral sections
and  higher Atiyah-Patodi-Singer index
theory on Galois coverings. {\it GAFA}. {\bf 8}, p. 17-58
%(http://www.dmi.ens.fr/dmi/preprints).
\item{[LP3]} Leichtnam E. and Piazza P. 1999. Homotopy
invariance of twisted higher signatures on manifolds with boundary. 
{\it  Bull. Soc. Math. France}. {\bf 127}, p. 307-331
\item{[LP4]} Leichtnam E. and Piazza P. 1998. A higher
APS index theorem for the signature operator on Galois coverings.
to appear, {\it Ann. Glob. Anal. and Geom.}
%Higher eta invariants and the 
%Novikov conjecture on manifolds with boundary.{\it CRAS} t. {\bf 327} 
%serie I, p. 497 (Geometrie differentielle).
\item{[L1]} Lott J. 1992. Superconnections and higher index theory. 
{\it GAFA}. {\bf 2}, p. 421-454  
\item{[L2]} Lott J. 1992. Higher eta invariants. {\it  K-Theory}.
{\bf 6}, p. 191-233
\item{[L3]} Lott J. 1997. Diffeomorphisms, analytic torsion 
and noncommutative 
geometry. {\it  Mem. of the Amer. Math. Soc} 141, no. 673. viii + 56 pp.
\item{[L4]} Lott J. 1996. 
The zero-in-the-spectrum question. {\it  Enseign. Math}.
{\bf 42}, p. 341-376
\item{[L5]} Lott J. 1998. Signatures and higher
 signatures of $S^1-$quotients. {\it Math. Annalen}, to appear
\item{[LL]} Lott J. and L\"uck W. 1995. $L^2$-topological invariants of
$3$-manifolds. {\it Invent. Math.} {\bf 120}, p. 15-60
\item{[Lu]} Lusztig G. 1971. Novikov's higher signature and families of
elliptic operators. {\it J. Diff.Geom.} {\bf 7}, p. 229-256 
\item{[Me]} Melrose R. 1993. {The 
Atiyah-Patodi-Singer index theorem}. 
A. and K. Peters. Wellesley, MA
\item{[MP1]} Melrose R. and Piazza P. 1997.
Families of Dirac operators, boundaries and the $b$-calculus.  
{\it J. Diff. Geom.} 
{\bf 46}, p. 99-180
\item{[MP2]} Melrose R. and Piazza P. 1997. 
An index theorem for families of Dirac operators
on odd dimensional manifolds with boundary. {\it J. Diff. Geom.} 
{\bf 46}, p. 287-334
\item{[MF]} Mischenko A. and Fomenko A. 1979. 
The index of elliptic operators over $C^*$-algebras.
{\it Izv. Akad. Nauk. SSSR, Ser. Mat.} 
{\bf 43}, p. 831-859
\item{[N]} Neumann W. 1975. Manifold cutting and pasting groups.
{\it Topology}. {\bf 14}, p. 237-244
\item{[R]} Ranicki, A. 1998. Higher-dimensional knot theory. Springer-Verlag.
\hfil\break
\noindent
Additions and errata at http://www.maths.ed.ac.uk/\~{ }aar/books
\item{[RW]} Rosenberg J. and Weinberger, S. 1993. Higher $G$-signatures for
Lipschitz manifolds.
{\it K-Theory}. {\bf 7}, p. 101-132
\item{[We]} Wegge-Olsen N. 1993. $K$-Theory and $C^*$-algebras.
Oxford University Press. New York
\item{[Wh]} Whitney H. 1957. Geometric integration theory.
Princeton University Press. Princeton, N.J.
\item{[Wu]} Wu F. 1997. The higher $\Gamma$-index for coverings of manifolds
with boundaries. {\it Fields Inst. Commun.} {\bf 17}, {\it 
Cyclic cohomology and
noncommutative geometry}. p. 169-183

\message{******** Run TeX twice to resolve cross-references}
\message{******** Run TeX twice to resolve cross-references}
\message{******** Run TeX twice to resolve cross-references}
\message{******** Run TeX twice to resolve cross-references}

\bye

ENDPAPER:

\end